\documentclass{amsart}
\usepackage{xypic}
\usepackage{url}
\begin{document}

\newtheorem{philos}[equation]{Philosophy}
\newtheorem{method}[equation]{Method}
\newtheorem{theorem}[equation]{Theorem}
\newtheorem{conjecture}[equation]{Conjecture}
\newtheorem{notation}[equation]{Notation}
\newtheorem{assumption}[equation]{Assumption}
\newtheorem{proposition}[equation]{Proposition}
\newtheorem{prop}[equation]{Proposition}
\newtheorem{lemma}[equation]{Lemma}
\newtheorem{corollary}[equation]{Corollary}
\theoremstyle{definition}
\newtheorem{definition}[equation]{Definition}
\newtheorem{remark}[equation]{Remark}
\newtheorem{example}[equation]{Example}
\newtheorem{wish}[equation]{Wish}
\newtheorem{situation}[equation]{Situation}
\theoremstyle{plain}
\numberwithin{equation}{section}
\numberwithin{figure}{section}

\title[On the $p$-adic Beilinson conjecture for number fields]{On the $\mathbf p$-adic Beilinson conjecture for number fields.}

\author{Amnon Besser}
\address{Department of Mathematics\\Ben-Gurion University of the Negev\\P.O.B. 653\\Be'er-Sheva 84105\\Israel}

\author{Paul Buckingham}
\address{Department of Pure Mathematics\\University of Sheffield\\Hicks Building\\Hounsfield Road\\Sheffield S3 7RH\\United Kingdom}

\author{Rob de Jeu}
\address{Faculteit Exacte Wetenschappen\\Afdeling Wiskunde\\Vrije Universiteit\\De Boelelaan 1081a\\1081 HV Amsterdam\\The Netherlands}

\author{Xavier-Fran\c cois Roblot}
\address{Institut Camille Jordan \\ Universit\'e de Lyon, Universit\'e
  Lyon~1, CNRS -- UMR 5208 \\ 43 blvd du 11 Novembre 1918, 69622 Villeurbanne Cedex \\ France}

\dedicatory{\hfill Dedicated to Jean-Pierre Serre on the occasion of his eightieth birthday.}

\begin{abstract}
We formulate a conjectural $p$-adic analogue of Borel's theorem relating regulators for higher $K$-groups
of number fields to special values of the corresponding $ \zeta $-functions, using syntomic regulators and
$p$-adic $L$-functions.
We also formulate a corresponding conjecture for Artin motives, and state a conjecture
about the precise relation between the $ p $-adic and classical situations.
Parts of the conjectures are proved when the number field (or Artin motive) is Abelian over the rationals,
and all conjectures are verified numerically in some other cases.
\end{abstract}



\maketitle

\def\spaceifletter{\futurelet\comingchar\dospaceifletter}
\def\dospaceifletter{\relax\ifmmode\else
  \ifcat A\noexpand\comingchar{} \fi
  \ifcat 0\noexpand\comingchar
  \ifx 0\noexpand\comingchar{} \fi
  \ifx 1\noexpand\comingchar{} \fi\ifx 2\noexpand\comingchar{} \fi
  \ifx 3\noexpand\comingchar{} \fi\ifx 4\noexpand\comingchar{} \fi
  \ifx 5\noexpand\comingchar{} \fi\ifx 6\noexpand\comingchar{} \fi
  \ifx 7\noexpand\comingchar{} \fi\ifx 8\noexpand\comingchar{} \fi
  \ifx 9\noexpand\comingchar{} \fi\fi
  \ifcat $\noexpand\comingchar{} \fi
  \ifcat \noexpand\relax\noexpand\comingchar{} \fi
\fi}

\def\art{\mathcal{M}_\pi}
\def\artR{\mathcal{M}_{\pi/\R}}
\def\a{\alpha}
\def\ab{\textup{ab}}
\def\Alt{\textup{Alt}}
\def\b{\beta}
\def\B #1 {B_{#1}(k)}
\def\Bei{\textup{Bei}}
\def\cl{\textup{Cl}}
\def\C{\mathbb C}
\def\Cp{{\mathbb C_p}}
\def\Drt{D_k^{1/2,\infty}}
\def\Drtp{D_k^{1/2,p}}
\def\Eul{\textup{Eul}}
\def\dd{{\textup{d}}}
\def\dr{\textup{dR}}
\def\f{\mathfrak f}
\def\Fr{\textup{Fr}}
\def\Gal{\operatorname{Gal}}
\def\hb{H_B}
\def\hd{H_{\mathcal D}}
\def\hdr{H_\dr}
\def\hm{H_M}
\def\hsyn{H_\syn}
\def\Hom{\operatorname{Hom}}
\def\spec{\operatorname{Spec}}
\def\isom{\cong}
\def\id{\textup{id}}
\def\Ind #1 #2 {\textup{Ind}_{#1}^{#2}}
\def\iso{\cong}
\def\kbar{\overline{k}}
\def\Ker{\operatorname{Ker}}
\def\khat{\widehat k}
\def\KMn #1 {{K_{2n-1}(#1)}}
\def\Kn #1 {K_{2n-1}(#1)_\Q}
\def\KnE #1 {K_{2n-1}(#1)_E}
\def\L #1 #2 #3 {L(#1,#2,#3)}
\def\Ls #1 #2 #3 {L^\ast(#1,#2,#3)}
\def\Lsharp #1 #2 #3 {L^\sharp(#1,#2,#3)}
\def\Lfunction{$L$-function\spaceifletter}
\def\Lfunctions{$L$-functions\spaceifletter}
\def\Li{\textup{Li}}
\def\Lmod #1 {P_{#1,p}}
\def\loc{\textup{loc}}
\def\Lp #1 #2 #3 #4 {L_{#1}(#2,#3,#4)}
\def\m{\mu}
\def\myatop #1 #2 {\genfrac..{0pt}{200pt}{#1}{#2}}
\def\mym{\mathbf{m}}
\def\Nm{\textup{Nm}}
\def\O{{\mathcal O}_k}
\def\OO{\mathcal O}
\def\padic{$p$-adic\spaceifletter}
\def\pint{\int^{\prime}}
\def\P{\mathfrak P}
\def\Pn{P}
\def\PP{{\P'}}
\def\ppl{\mathcal{P}_+}
\def\Q{\mathbb Q}
\def\Qp{{\mathbb Q_p}}
\def\Qpbar{{\overline {\mathbb Q}_p}}
\def\Qbar{\overline{\Q}}
\def\r{\rho}
\def\R{\mathbb R}
\def\Re{\textup{Re}}
\def\reg{\textup{reg}}
\def\rightiso{\buildrel{\sim}\over{\rightarrow}}
\def\s{\sigma}
\def\sgn{\sigma}
\def\smallsum{\textstyle\sum}
\def\syn{\textup{syn}}
\def\tB{\widetilde B}
\def\tensor{\otimes}
\def\tor{\textup{torsion}}
\def\t{\tau}
\def\tLmod #1 {{\widetilde P}_{#1,p}}
\def\V #1 #2 {V_{#1,#2}}
\def\W #1 #2 {W_{#1,#2}}
\def\z{\zeta}
\def\Z{\mathbb Z}
\def\Zp{\mathbb Z_p}

\section{Introduction}

The Beilinson conjectures about special values of \Lfunctions~\cite{bei85}
are a far reaching generalization of the class number formula for the
Dedekind zeta function. For every algebraic variety $X$ over the
rationals it predicts the leading term of the Taylor expansions of $L(H^i(X),s)$
at certain points, up to a
rational multiple, in terms of
arithmetic information associated with $X$, namely, its algebraic
K-groups $K_j(X)$~\cite{qui73b}. More generally, these conjectures can also be formulated for motives.

There have been several important steps taken towards verification of
these conjectures in various cases, although, strictly speaking they
have only been verified completely in the case where $X$ is the
spectrum of a number field, where they follow from famous theorems 
of Borel~\cite{Borel74, Borel77}.

To motivate what follows, let us briefly recall the conjecture that
interests us the most (for an introduction see~\cite{DenSch91,schn88}).
One associates with $X$ two cohomology 
groups. The first one is the Deligne cohomology $\hd^i(X_{/\R},\R(n)) $, which is an
$\R$-vector space. The second is \lq\lq integral\rq\rq\ motivic cohomology $\hm^i(X_{/\Z},\Q(n)) $,
which may be defined as a certain subspace of $K_{2n-i}(X)\tensor_\Z \Q$.
There is a regulator map defined by Beilinson,
\begin{equation}\label{regmap}
 \hm^i(X_{/\Z},\Q(n)) \to \hd^i(X_{/\R},\R(n))
\,.
\end{equation}
If $ 2n > i+1 $ then  $\det\hd^i(X_{/\R},\R(n))$ has a rational structure
coming from the relations between $ \hd^i(X_\R,\R(n)) $ and the de Rham and singular cohomology groups of~$X$
\cite[p.30]{schn88}.

The first part of the conjecture is that
the map in \eqref{regmap} induces an isomorphism between the lefthand side tensored with $\R$ and the
righthand side, and
consequently provides a second rational structure on $\det\hd^i(X_{/\R},\R(n))$.
The second part of the conjecture states that, assuming a suitable functional equation for
$ L(H^{i-1}(X),s) $, these two rational structures differ from each other by the
leading term in the Taylor expansion of $ L(H^{i-1}(X),s) $
at $ s=i-n $.  Because of the expected functional equation one can reformulate the conjecture in terms
of $ L(H^{i-1}(X),n) $ (see \cite[Corollary~3.6.2]{bei85} or \cite[4.12]{Jan88b}).

As mentioned before, this conjecture has only been verified in the case of number
fields, due to difficulties in the computation of motivic
cohomology. What has been verified in several other cases is a 
form of the conjecture in which one assumes the first part.
For this one finds $\dim_{\R} \hd^i(X_{/\R},\R(n)) $ elements of
$\hm^i(X_{/\Z},\Q(n))$, checks that their images under \eqref{regmap} are
independent, hence should form a basis of $ \hm^i(X_{/\Z},\Q(n))$ according to the
first part of the conjecture, and verifies the second part using these
elements.

The idea that there should be a $p$-adic analogue of Beilinson's
conjectures has been around since the late 80's.  Such a conjecture was
formulated and proved by Gros in~\cite{Gro90,Gro94} in the case of Artin
motives associated with Dirichlet characters.
In the weak sense mentioned before, it was proved for certain CM elliptic curves
in~\cite{Col-Gro89} (where the relation with the syntomic regulator is proved
in~\cite{Bes98b} and further elucidated in~\cite{Bes-Den99}), and for elliptic
modular forms it follows from Kato's work (see~\cite{Scho98}).

The book~\cite{Peri96} contains a very general
conjecture about the existence and properties of \padic
\Lfunctions, from which one can derive a \padic Beilinson conjecture.
Rather than explain this in full detail we shall give a 
sketch of this conjecture similar to the sketch above of the
Beilinson conjecture.

For the \padic Beilinson conjectures one has
to replace Deligne cohomology with syntomic cohomology~\cite{Gro94,Niz97,Bes98a}, 
the Beilinson regulator with the syntomic regulator, and \Lfunctions with \padic
\Lfunctions.

Syntomic cohomology $\hsyn^i(Y,n)$ is defined for a smooth
scheme $ Y $ of finite type over a complete discrete valuation ring of mixed
characteristic $(0,p)$ with perfect residue field.
For the $\Q$-variety $X$ we obtain, for all but finitely many primes, a map
\begin{equation}\label{regpmap}
 \hm^i(X_{/\Z},\Q(n)) \to \hsyn^i(Y,n)
\end{equation}
where $Y$ is a smooth model for $X$ over $\Z_p$. This cohomology group is a
$\Qp$-vector space.

The theory of \padic \Lfunctions starts with Kubota and Leopoldt's
\padic $ \z $-function \cite{Kub-Leo}, $ \z_p(s) $, which is defined by
interpolating special values of complex valued Dirichlet \Lfunctions.
This principle has been extended to $ \z $- and \Lfunctions in various situations, resulting in 
corresponding \padic functions for totally real number fields~\cite{Bar78,Cas-Nog79,Del-Rib},
CM fields~\cite{Kat76,Kat78} and modular forms~\cite{MTT86}.
(Given the occasion, let us note that the approach of Deligne and Ribet using modular forms was
initiated by Serre~\cite{Ser73}.)

The $p$-adic Beilinson conjecture therefore has many
similarities with its complex counterpart. However, there is a very
important difference. In general, when $  2n>i+1 $, there is no hope that~\eqref{regpmap}
induces an isomorphism after tensoring the lefthand side with $\Q_p$. 
To see this, consider a number field~$k$ with ring of algebraic integers $ \O $.
By Borel's theorem (see Theorem~\ref{boreltheorem}) we have,
in accordance with Beilinson's conjectures,
\begin{equation*}
  \dim_\Q \hm^1(\spec(k)_{/\Z},\Q(n))=
\left\{ 
\begin{aligned}
      r_2 & \text{ when $ n \ge 2 $ is even }
\\
r_1 + r_2 & \text{ when $ n \ge 2 $ is odd }
\end{aligned}
 \right.
\end{equation*}
with $r_1$ (resp.\ $ 2 r_2 $) the number of real (resp.\ complex)
embeddings of $k$. Thus, motivic cohomology ``knows'' about the number
of real and complex embeddings of $k$, as does, by its definition, 
Deligne cohomology, which in this case becomes
$ \hd^1(\spec(k)_{/\R},\R(n)) \iso \{ (x_\s)_{\s: k \to \C} \text{ in } \R(n-1)^{r_1+2r_2} \text{ such that } \overline{x_\s} = x_{\overline \s}\} $.
But syntomic cohomology, which depends only on the completion of~$k$  at~$p$, does
not.  In this case, we obtain
$ \hsyn^1(\spec(\O\tensor_\Z\Zp),n) \iso \Q_p^{r_1+2r_2} $.

The solution to this problem suggested in~\cite{Peri96} is to make the
\padic \Lfunction depend on a subspace of syntomic cohomology which
is complementary to the image of the regulator. While this might seem
artificial, there are other reasons for choosing
this solution. In most cases one chooses a particular subspace and
obtains a special case of the conjecture.

In the important special case of a totally real number field, or more generally an
Artin motive over $ \Q $, associated with a Galois
representation whose kernel is totally real 
(let us call these totally real Artin motives), no such subspace
is required when $ n\ge2 $ is odd
(see Proposition~\ref{whenconj}).  The same holds for Artin motives where the conjugacy class of complex
conjugation acts as multiplication by $ -1 $ and $ n \ge2 $ is even (we may think of those as the ``negative
part" of CM Artin motives).

If $ \chi $ denotes the character associated with either representation, then
the Beilinson conjecture relates the regulator of $ K_{2n-1} $ with the Artin \Lfunction of $ \chi $
at $ n $.
For the \padic \Lfunctions one has to consider $ \chi \tensor \omega_p^{1-n} $ with $ \omega_p $ the Teichm\"uller
character for the prime number $ p $.  Then the fixed field of
the kernel of the representation underlying
$ \chi \tensor \omega_p^{1-n} $ is totally real,
and it is perhaps no coincidence that in precisely this
case the existence of a  \padic \Lfunction that is not identicaly zero has been
established, by~\cite{Bar78,Cas-Nog79,Del-Rib} for the case of fields
and by~\cite{greenberg83} for Artin motives. 

The goal of the present work is to describe in detail the conjectures
for the cases of totally real fields, totally real Artin motives, as well as the ``negative part" of CM
Artin motives, and describe the (conjectural) relation between the classical and \padic conjectures.
We test everything numerically in several cases, and also deduce most of the conjectures for Abelian Artin
motives from work by Coleman~\cite{Col82}.

There have been several developments that allow us to carry out this
numerical verification. In~\cite{dJ95} de Jeu proved part of Zagier's
conjecture concerning a more explicit description of the $ K $-theory
(tensored with $ \Q $) of number fields  (see Section~\ref{ksection}).
While this conjecture is not known to give the K-theory of such fields in all cases,
it does in practice.  Thus it provides a way of computing them, and Paul Buckingham
wrote a computer implementation for this.  In~\cite{BdJ03} Besser and de Jeu computed the syntomic regulator for
(essentially) the part of the $ K $-theory of a number field described by Zagier's conjecture
 and showed that it is given by applying the \padic polylogarithm.
Those \padic polylogarithms were invented by Coleman~\cite{Col82}
using his theory of \padic integration but are not so easy to
compute. In~\cite{BdJ06} Besser and de Jeu devised an
algorithm for this computation. Taken together,
these developments allow us to compute~\eqref{regmap} and ~\eqref{regpmap}
for number fields.
Finally, building on earlier work in \cite{Sol-Rob}, Roblot has dealt with the computational aspects of computing
\padic \Lfunctions for Abelian characters over $ \Q $ or a real quadratic field \cite{Rob06}.

This paper is organized as follows.  In Section~\ref{conjecturesection} we recall Borel's theorem as
well as various facts about \Lfunctions and \padic \Lfunctions, and formulate a conjectural \padic analogue of
Borel's theorem.
In Section~\ref{motivicsection} we introduce Artin motives with coefficients in a number field $ E $
in terms of representations of the Galois group,
determine when the $ \Q $-dimension of the left-hand side of~\eqref{regpmap}
equals the $ \Qp $-dimension of the right-hand side (corresponding to equality in Proposition~\ref{whenconj}),
define both classical and \padic \Lfunctions with coefficients in $ E $, and formulate the motivic Beilinson
conjecture with coefficients in $ E $, Conjecture~\ref{motiveconjecture}, a small part of which is the
same as a conjecture by Gross.
In Section~\ref{ksection} we describe the set-up for finding elements in the $ K $-groups of number fields
using Zagier's conjecture, and the classical and \padic regulators on them.
We also prove most of Conjecture~\ref{motiveconjecture} for Abelian Artin motives over $ \Q $ (Proposition~\ref{abelprop}
and Remark~\ref{abelrem}).
In Section~\ref{practical} we discuss a few computational aspects of implementing Zagier's conjecture
and describe the Artin motives that we consider later for the numerical examples,
and in the process prove Gross's conjecture for Artin motives obtained from $ S_3 $ and $ D_8 $-extensions of $ \Q $. 
Then in Section~\ref{pLsection} we sketch how the \padic \Lfunctions can be computed in certain cases,
and make the required Brauer induction explicit for the Artin motives we want to consider.  Finally,
the last section is devoted to the results of the numerical calculations for examples.

Amnon Besser, Paul Buckingham and Rob de Jeu would like to thank the
EC network Arithmetic Algebraic Geometry for travel support.  Rob de
Jeu would like to thank the Tata Institute of Fundamental Research for
a productive stay during which this paper was worked on, and the
Nuffield Foundation for a grant under the Undergraduate Research
Bursary programme (NUF-URB03) that enabled Paul Buckingham to develop
a computer program for finding the necessary elements in $ K $-theory.
Amnon Besser would like to thank the Institute for Advanced Study in
Princeton for a stay during which part of the paper was worked on,
and the Bell companies Fellowship and the James D. Wolfensohn fund
for financial support while at the institute.  Finally, the authors
would like to thank Alfred Weiss for very useful explanations about
\padic \Lfunctions, and, in particular, for bringing to their attention
the paper~\cite{greenberg83}.

\medskip
\noindent
{\bf Notation}\hfil\break
Throughout the paper, for an Abelian group $ A $ we let $ A_\Q $ denote $ A \tensor_\Z \Q $.
If $ B $ is a subgroup of $ \C $ and $ n $ an integer then we let $ B(n) = (2 \pi i)^n B \subseteq \C $
We normalize the absolute value $ |\cdot|_p $ on the field of $p$-adic numbers $\Qp$ 
in such a way that $|p|_p= p^{-1}$, and use the same notation for its extensions to
an algebraic closure $ \Qpbar $ and $\Cp = \hat{\overline\Q}_p$.

\section{The \padic Beilinson conjecture for totally real fields}
\label{conjecturesection}

Let $ k $ be a number field with $ r_1 $ real embeddings, $ 2r_2 $ complex embeddings, 
ring of algebraic integers $ \O $, and discriminant $ D_k $.
As is well-known, $ \O^* $ is a finitely generated Abelian group
of rank $ r = r_1 + r_2 - 1 $, and its regulator $ R $ satisfies
$  w \sqrt{  |D_k| } \, \text{Res}_{s=1}\zeta_k(s) =  2^{ r_1 } (2\pi)^{ r_2 } \, |\text{Cl}(\O )| \, R $,
with $ \text{Cl}(\O) $ the class group of $ \O $, and $ w $ the number of roots of unity in $ k $.
Because  $ K_0(\O) \iso \text{Cl}(\O) \oplus \Z $ and $ K_1(\O) \iso \O^* $, 
so $ |\text{Cl}(\O)| = |K_0(\O)_{\text{tor}}| $
and $ w = | K_1(\O)_{\text{tor}} | $, this can be interpreted
as a statement about the
$ K $-theory of $ \O $, and it is from this point of view that
it can be generalized to $ \zeta_k(n) $ for $ n \geq 2 $.
Namely, in \cite{qui73b}, Quillen proved that $ K_m(\O) $ is a finitely generated Abelian group for all~$ m $.
Borel in \cite{Borel74} computed its rank when $ m \ge 2 $.  For $ m $ even this rank is zero,
but for odd $ m $ it is $ r_1+r_2 $ or $ r_2 $, and in \cite{Borel77}
he showed that a suitably defined regulator of $ K_{2n-1}(\O) $ is related to $ \zeta_k(n) $.
Since $ K_{2n-1}(\O) \rightiso K_{2n-1}(k) $ when $ n \ge 2 $
we can rephrase his results for $ K_{2n-1}(\O) $ in terms of $ K_{2n-1}(k) $.
Also, we replace Borel's  regulator map $ \reg_B : K_{2n-1}(\C) \to \R(n-1) $ ($ n \ge 2 $)
with Beilinson's regulator map $ \reg_\infty $ (see \cite[\S~4]{schn88}), which is half the Borel map
regulator by \cite[Theorem~10.9]{Bur02}.
Because $ k\tensor_\Q\C \iso \oplus_{\s: k \to \C} \C $,
and $ n > \dim \spec(k\otimes \C) $, we obtain by \cite[p.9]{schn88}
\begin{equation}\label{HDid}
\begin{aligned}
   \hd^1(\spec(k\otimes \C)_{/\R},\R(n))
& \iso
   H^0(\spec(k\otimes \C)_{/\R},\R(n-1))
\cr
& \iso
   \Bigl(\bigoplus_{\s: k\to \C} \R(n-1)\Bigr)^+
\,,
\end{aligned}
\end{equation}
which consists of those $ (x_\s)_\s $ with $ x_{\overline\s} = \overline{x_\s} $.
Finally, for any embedding $ \s : k \to \C $ we let $ \s_* : K_{2n-1}(k) \to K_{2n-1}(\C) $ be the induced map.

\begin{theorem}\label{boreltheorem}{\bf (Borel)}
Let $ k $ be a number field of degree $ d $, with $ r_1 $ real embeddings and $ 2 r_2 $ complex embeddings, and let $ n \ge 2 $.
Then the rank $ m_n $ of $ K_{2n-1}(k) $ equals $ r_2 $ if $ n $ is even and $ r_1+r_2 $ if $ n $ is odd.
Moreover, the map
\begin{equation}\label{Bregulator}
\begin{aligned}
 K_{2n-1}(k) & \to \bigoplus_{\s : k \to \C } \R(n-1)
\\
\a & \mapsto (\reg_\infty\circ\s_*(\a))_\s
\end{aligned}
\end{equation}
embeds $ K_{2n-1}(k)/\tor $ as a lattice in
$ (\oplus_\s \R(n-1))^+ \iso \R^{m_n} $,
and the volume $ V_n(k) $ of a fundamental domain of this lattice satisfies
\begin{equation}\label{borelq}
 \zeta_k(n)  \sqrt{|D_k|} = q \pi^{n(d-m_n)} V_n(k) 
\end{equation} 
for some $ q $ in $ \Q^* $.
\end{theorem}

\begin{remark}\label{klingensiegelremark}
In Theorem~\ref{boreltheorem} $ m_n = 0 $ precisely when $ k $ is totally real and $ n \ge 2 $
is odd.  In this case the given relation holds (with $ V_n(k) = 1 $)
by the Siegel-Klingen theorem \cite[Chapter~VII, Corollary~9.9]{Neu99}.
\end{remark}

This theorem is equivalent with Beilinson's conjecture for $ k $.
In order to deal with this in detail 
(see Remarks~\ref{BBeq} and~\ref{BBmotiverem}) and in order to introduce \padic \Lfunctions 
we recall some facts about Artin \Lfunctions \cite[Chapter VII, \S~10-12]{Neu99}.
Let $ k $ be a number field, $ d = [k:\Q] $, 
and let $ \chi $ be a $ \C $-valued Artin character of $ \Gal(\kbar/k) $.
For a prime number $ \ell $ we define 
\[
\Eul_\ell(s, \chi, k) =  \prod_{\mathfrak{l} \mid \ell} \Eul_{\mathfrak l}(s,\chi,k)
\]
where $ \Eul_{\mathfrak l}(s,\chi,k) $ is the reciprocal of the Euler factor for $ \mathfrak l $
and the product is over primes $ \mathfrak l $ of $ k $ lying above $ \ell $.  Then for $ s \in \C $
with $ \Re(s)>1 $ we can write the Artin \Lfunction of $ \chi $ as
\[
 \L s {\chi} k = \prod_\ell \Eul_\ell(s,\chi,k)^{-1} 
\,.
\]
For an infinite place $ v $ of $ k $ we let
\begin{equation*}
L_v(s, \chi, k) = \left\{
\begin{aligned}
& \hbox to 1.5 in {$ L_\C(s)^{\chi(1)} $\hfill}\hbox to 1.5 in{when $ v $ is complex,\hfill}
\\
& \hbox to 1.5 in {$ L_\R(s)^{n^+} L_\R(s+1)^{n^-} $\hfill}\hbox to 1.5 in{when $ v $ is real,\hfill}
\end{aligned}
\right.
\end{equation*}
where $ L_\R(s) = \pi^{-s/2} \Gamma(s/2) $, $ L_\C(s) = 2 (2\pi)^{-s}  \Gamma(s) $
and
$ n^\pm = \frac d2 \overline\chi(1) \pm \sum_v  \frac 12 \overline\chi(\Fr_w) $,
with the sum taken over the real places $ v $ of $ k $
and $ \Fr_w $ the generator of the image of $ \Gal(\C/\R) $ in $ \Gal(\kbar/k) $
corresponding to any extension $ w : \kbar \to \C $ of~$ v $.
Then 
\begin{equation*}
\Lambda(s, \chi, k) =  \L s {\chi} k \prod_{v|\infty} L_v(s, \chi, k)  
\end{equation*}
extends to a meromorphic function, satisfying the functional equation
(cf.\ \cite[p.\ 541]{Neu99})
\begin{equation*}
\Lambda(s, \chi, k) = W(\chi) C(\chi,k)^{1/2}  C(\chi,k)^{-s} \Lambda(1-s, \overline\chi, k) 
\end{equation*}
with $ W(\chi) $ a constant of absolute value~1
and $ C(\chi,k) = |D_k|^{\chi(1)} \Nm_{k/\Q}(\f(\chi,k)) $ for $ \f(\chi,k) $ the Artin conductor of $ \chi $.
Therefore
\begin{equation}\label{oldfuneq}
\begin{aligned}
\L 1-s {\overline\chi} k 
& =
W(\overline\chi) C(\overline\chi,k)^{s-\frac12} (2 (2\pi)^{-s}\Gamma(s))^{d\overline\chi(1)}
\\
& \qquad\times
 (\cos(\pi s/2))^{n^+}  (\sin(\pi s/2))^{n^-} \L s {\chi} k 
\,.
\end{aligned}
\end{equation}

Following \cite[pp.~80--81]{greenberg83} we shall now describe a \padic \Lfunction $ \Lp p s {\chi} k $ when
$ k $ is a totally real number field, $ p $ a prime and 
$ \chi : \Gal(\kbar/k) \to \Qpbar $ a suitable Artin character.
We begin with the case of 1-dimensional Artin characters.

If $\sigma:\Qpbar \to \C$ is any isomorphism then
$\sigma \circ \chi$ is a complex Artin character, so we have the Artin \Lfunction $  \L s {\sigma\circ\chi} k $.
We may also view $ \chi $ as a character of a suitable ray class group, so
by \cite[Corollary~9.9 and page~509]{Neu99} all $ \L m {\sigma\circ\chi} k $ for $m \in \Z_{\leq 0}$ 
are in $ \Q(\s\circ\chi) $ and the values
\begin{equation}\label{sigmaindependence}
  \Ls m {\chi} k = \sigma^{-1}(\L m {\sigma\circ\chi} k )
\end{equation}
are independent of the choice of $\sigma$.

In the same way, we define
\begin{equation}\label{sigmaeuler}
\Eul^\ast_\ell(m, \chi, k) = \sigma^{-1} (\Eul_\ell(m,\sigma\circ\chi, k)) \,,
\end{equation}
which is clearly independent of the choice of $ \s $.  To construct
the \padic \Lfunction one finds a \padic analytic or
meromorphic
function on an open ball around 0 that interpolates the values 
$\Ls m {\chi} k $.  For 1-dimensional $ \chi $ with
the fixed field $ k_\chi $ of its kernel totally
real this was achieved independently by Deligne and Ribet~\cite{Del-Rib},
Barsky~\cite{Bar78}, and Cassou-Nogu\`es \cite{Cas-Nog78, Cas-Nog79}.  We shall sketch a proof of
the following theorem and Remark~\ref{interpolrem} below in Section~\ref{pLsection}.
(When $ k_\chi $ is not totally real the \padic \Lfunction is
identically zero since the values interpolated are all zero by the
functional equation of the \Lfunction.)

\begin{theorem}\label{padicLfunction}
For $ p $ prime, let $\mathcal{B}$ in $ \Cp $ be the open ball with centre $0$ and radius $qp^{-1/(p-1)}$
where $q = p$ if $ p>2 $ and $ q=4 $ if $p=2$.
If $ k $ is a totally real number field 
and $ \chi : \Gal(\kbar/k) \to \Qpbar $ a $ 1 $-dimensional Artin character, then there exists a unique 
$ \Cp $-valued function $ \Lp p s {\chi} k $ on $ \mathcal{B} $
satisfying the following properties:
  \begin{enumerate}
  \item $\Lp p s {\chi} k $ is analytic if $\chi$ is non-trivial and
     meromorphic with at most a simple pole at $s = 1$ if $\chi$ is trival;
  \item if $m$ is a negative integer such that $m \equiv 1 $ modulo $\varphi(q)$ then
    \begin{equation*}
      \Lp p m {\chi} k = \Eul^\ast_p(m,\chi,k) \Ls m {\chi} k 
      \,. 
    \end{equation*}
  \end{enumerate}
\end{theorem}

If $\chi : \Gal(\kbar/k) \to \Qpbar $ is any Artin character then by Brauer's induction theorem
\cite[(10.3)]{Neu99}
there exist 1-dimensional Artin characters $\chi_1, \dots,\chi_t$ on subgroups
$ G_1,\dots,G_t $ of $ \Gal(\kbar/k) $
of finite index,
and integers $a_1, \dots, a_t$, such that
\begin{equation}\label{brauerinduction}
  \chi = \sum_{i=1}^t a_i \Ind {G_i} {\Gal(\kbar/k)} (\chi_i)
\,.
\end{equation}
If $ k_\chi $ is
totally real then we can assume that the same holds for the fixed fields $ k_i $ of the $ G_i $, and
we define the \padic \Lfunction of $\chi$ by
\begin{equation}\label{artinLp}
  \Lp p s {\chi} k = \prod_{i=1}^t \Lp p s {\chi_i} k_i  ^{a_i}
\,,
\end{equation}
which is a meromorphic function on $\mathcal{B}$ (see Section~\ref{pLsection}).  For any 
negative integer $m \equiv 1 $ modulo $ {\varphi(q)} $ the value
$ \Lp p m {\chi} k $ is defined and equals
$  \Eul^\ast_p(m,\chi,k) \Ls m {\chi} k $ by well-known properties of Artin
\Lfunctions~(see \cite[Prop.~10.4(iv)]{Neu99}),
showing the function is independent of how we express $\chi$ as a sum of induced 1-dimensional characters.

\begin{remark}
  In \cite{greenberg83}, Greenberg proves that the Main Conjecture of
  Iwasawa theory implies the $p$-adic Artin conjecture, that is, that
  the \padic \Lfunction of an Artin character $\chi$ is analytic on
  the open ball $\mathcal{B}$ if it does not contain the trivial
  character, and has at most a simple pole at $ s=1 $ otherwise.
  It therefore follows from the proof of the Main Conjecture by
  Mazur and Wiles \cite{MW84} for $ p $ odd that $\Lp p n {\chi} k $ for $ p\ne2 $ is
  defined for all integers $n\ne1$.  In particular, the values of the \padic \Lfunctions
  in Conjecture~\ref{pborel} below exist by Theorem~\ref{padicLfunction},
  and those in Conjecture~\ref{motiveconjecture} exist when $ p\ne2 $ but have to be assumed
  to exist when $p = 2$.
\end{remark}

\begin{remark}\label{interpolrem}
Let $W_p \subset \Z_p^*$ be the subgroup of $(p-1)$-th
roots of unity if $p$ is odd and let $W_2 = \{\pm 1\}$.  The
Teichm\"uller character on $\Gal(\kbar/k)$ is defined as the composition
\begin{equation}\label{TM}
\omega_p: \Gal(\kbar/k) \to \Gal(\Qbar/\Q) \to \Gal(\Q(\mu_q)/\Q) \rightiso (\Z/q\Z)^* \rightiso W_p
\,,
\end{equation}
where the last map sends an element of $(\Z/q\Z)^*$ to the unique element of $W_p$ to which it is congruent modulo $q$. 
For an Artin character $ \chi : \Gal(\kbar/k) \to \Qpbar $ and an integer $ l $,
$\chi\omega_p^l$  is also an Artin character.
If $ m \le 0 $ satisfies $ l+ m \equiv 1 $ modulo $ {\varphi(q)} $ and either $ \chi $ is 1-dimensional
or $ k_\chi $ is totally real then $ \Lp p m {\chi\omega_p^l} k = \Eul_p^*(m, \chi, k) \Ls m {\chi} k $.
\end{remark}

The \padic Beilinson conjecture is going to predict the special
values of \padic \Lfunctions in terms of a \padic regulator. The
required regulator is the syntomic regulator~\cite[Theorem~7.5]{Bes98a}. Let $F$
be a complete discretely valued field of characteristic $0$
with perfect residue field of characteristic
$p$ and let $X$ be a scheme that is smooth and of  finite type over 
the valuation ring ${\OO}_F$.  Then the above mentioned paper associates to $X$ its
rigid syntomic cohomologies $\hsyn^i(X,n)$, as well as syntomic
regulators (i.e.,\ Chern characters) $\reg_p: K_{2n-i}(X) \to
\hsyn^i(X,n)$. In this work, unlike in~\cite{Bes98a}, we shall
need to change the base field $F$. We therefore prefer to denote the
syntomic cohomology by $\hsyn^i(X/\OO_F,n)$. For the formulation of the
conjecture the following basic fact is required.

\begin{lemma}\label{synlemma}
  We have $\hsyn^1(\spec(\OO_F)/\OO_F,n) \isom F$, for all
  $n>0$ and consequently we have a syntomic regulator
  $\reg_p: K_{2n-1}(\OO_F)\to F$. 
  Furthermore, the map $\reg_p$ commutes in the obvious way with
  finite extensions of fields and with automorphisms of such fields,
  provided that their residue fields are 
  algebraic over the prime field.
\end{lemma}

\begin{proof}
  The first claim follows from part 3
  of~\cite[Proposition~8.6]{Bes98a}.  For the second claim we note
  that 
  by part 4 of the same proposition we have in this case an
  isomorphism between syntomic and modified syntomic cohomology (the
  latter only exists under the additional assumption on the residue
  field). The compatibility with finite base changes now follows from
  the same result for modified syntomic cohomology
  in~\cite[Proposition~8.8]{Bes98a}. The same holds for automorphisms,
  although not explicitly stated in the above reference, since the
  relevant base change results, e.g., Proposition~8.6.4, hold for this
  type of base change as well.
\end{proof}

As a consequence of the lemma we can also, by abuse of notation,
define for any complete discretely valued subfield $F\subset \Qpbar$ the
regulator $\reg_p: K_{2n-1}(\OO_F)\to F$. In~\cite{BdJ03}, two of the authors of the present
work showed how one can sometimes compute the map $\reg_p$ by using
\padic polylogarithms.

We now restrict our attention to a totally real number field $ k $ with $ [k:\Q] = d $.  Our goal will
be to formulate a conjecture that is the \padic analogue of Theorem~\ref{boreltheorem} for $ K_{2n-1}(k) $
with $ n \ge 2 $.  Since this $ K $-group is torsion when $ n $ is even but has rank $ d $ when $ n $ is odd,
we only consider odd $ n\ge 2 $.

In preparation for the more general construction that will follow
in Section~\ref{motivicsection} let us first reformulate Theorem~\ref{boreltheorem} in this special case.
Let $ a_1,\dots,a_d $ form a $ \Z $-basis of $ \O $ and let $ \s_1^\infty,\dots,\s_d^\infty $ be the embeddings of $ k $ into $ \C $.
Then we define $ \Drt = \det(\s_i^\infty(a_j)) $, a real root of the discriminant of $ k $.
Similarly, if $ \a_1,\dots,\a_d $ form a $ \Z $-basis of $ K_{2n-1}(k)/\tor $ then we let
$ R_{n,\infty}(k) = \det(\reg_\infty\circ\s_{i*}^\infty(\a_j)) $.
Then $ \Drt $ and $ R_{n,\infty}(k) $ are well-defined up to sign, and the relation in
Theorem~\ref{boreltheorem} in this case is equivalent with
\begin{equation}\label{rat}
\zeta_k(n)  D_k^{1/2,\infty} = q(n,k) R_{n,\infty}(k) 
\end{equation}
with $ q(n,k) $ in $ \Q^* $.

Now let $ F \subset \Qpbar $ be the topological closure of the Galois closure of $ k $ embedded in $ \Qpbar $ in any way.
If $ \s_1^p,\dots,\s_d^p $ are the embeddings of $ k $ into $ F $ 
then $ \Drtp = \det(\s_i^p(a_j)) $ is a root in $ F $ of the discriminant of $ k $.
For $ \s : k \to F $ an embedding we denote the induced map $ K_{2n-1}(k) \to K_{2n-1}(F) $ by $ \s_* $.
Then we define a \padic regulator in $ F $ by $  R_{n,p}(k) = \det(\reg_p\circ\s_{i*}^p(\a_j)) $.  Both
$ \Drtp $ and $ R_{n,p}(k) $ are well-defined up to sign.

Remark~\ref{interpolrem} suggests that the role of $ \zeta_k(n) $ in
a \padic analogue of Theorem~\ref{boreltheorem} should be played by
$ \Lp p n \omega_p^{1-n} k / \Eul_p(n,k) $ where $ \z_k(s) = \prod_l \Eul_l(s,k)^{-1} $ for $ \Re(s) > 1 $,
so we can hope that 
\[
\Lp p  n \omega_p^{1-n} k  \Drtp  = q_p(n,k) \, \Eul_p(n,k) R_{n,p}(k)
\]
for some $ q_p(n,k) $ in $ \Q^* $.

More precisely, because $ R_{n,\infty}(k) / \Drt $ and $ R_{n,p}(k) / \Drtp $
are invariant under reordering the $ \s_i^\infty $ or $ \s_i^p $, and
transform in the same way if we change the bases of $ \O $ and $ K_{2n-1}(k)/\tor $, 
we can make the following conjecture.

\begin{conjecture}\label{pborel}
For $ k $ a totally real number field, $ p $ prime, and $ n \ge 2 $ odd, we have, with notation as above:
\begin{enumerate}

\item in $ F $ the equality
\[
\Lp p  n \omega_p^{1-n} k   D_k^{1/2,p}  =   q_p(n,k) \, \Eul_p(n,k) R_{n,p}(k) 
\]
holds for some $ q_p(n,k) $ in $ \Q^* $;
\item in fact, $ q_p(n,k) = q(n,k) $;
\item $ \Lp p n \omega_p^{1-n} k $ and $ R_{n,p}(k) $ are non-zero.
\end{enumerate}
\end{conjecture}
As mentioned in the introduction, this, and the corresponding parts of Conjecture~\ref{motiveconjecture}
below, can be deduced (with some effort) from a much more general conjecture of Perrin-Riou~\cite[4.2.2]{Peri96}.

\begin{remark}
The conjecture is similar to the result that the residue of $ \z_p(s,k) $ at $ s=1 $
is related to the Leopoldt regulator of $ \O^* $ through exactly the same formula as for the residue of $ \z(s,k) $ and
the Dirichlet regulator \cite{Colm88}, with part(4) corresponding to the Leopoldt conjecture.  However,
we have not tried to determine if Colmez's normalization of the sign for the regulator is
the same as here, especially given the sign error in the proof of Lemma~4.3 in
loc.\ cit.\ (see Section~\ref{pLsection}).
\end{remark}

\begin{remark}\label{basisremark}
We can use a basis of a subgroup of finite index of $K_{2n-1}(k)/\tor $ in the definition of $ R_{n,\infty}(k) $ and $ R_{n,p}(k) $,
or even a $ \Q $-basis of $ \Kn k $, without affecting the rationality of $ q(n,k) $, $ q_p(n,k) $ or
their equality.
Similarly we can replace the $ \Z $-basis of $ \O $ with a $ \Q $-basis of~$ k $ in the definitions of $ \Drt $ and~$ \Drtp $.
\end{remark}

\begin{remark}\label{BBeq}
As is well-known, for $ k $ and $ n $ as in Conjecture~\ref{pborel},
by~\eqref{oldfuneq} $ \z_k(s) $ at $ s=1-n $ has a zero of order $ d $ and
the first non-zero coefficient in its Taylor expansion, $ \z_k^\sharp(1-n) $, equals
$ (2 \pi i )^{d(1-n)} ((n-1)!/2)^d  D_k^{n-1/2} \z_k(n) $.
If we take $ \Drt = D_k^{1/2} $ in~\eqref{rat} then we obtain
\begin{equation*}
\z_k^\sharp(1-n) =  ((n-1)!/2)^d q(n,k) D_k^{n-1} \widetilde R_{n,\infty}(k)
\end{equation*}
for Beilinson's renormalized regulator $ \widetilde R_{n,\infty}(k) = (2\pi i)^{d(1-n)} R_{n,\infty}(k) $.
In computer calculations $ ((n-1)!/2)^d q(n,k) D_k^{n-1} $ often has only 
relatively small prime factors, so the larger prime factors in $ q(n,k) $ correspond to $ D_k^{1-n} $.
This phenomenon also occurs in the calculations for Conjecture~\ref{motiveconjecture} below
(see Remark~\ref{constantremark}).
\end{remark}

\begin{remark}
(1)
The thought that $ \Lp p n \omega_p^{1-n} k $ is non-zero for $ n \ge 2 $ and
odd when $ k $ is a totally real Abelian extension of~$ \Q $
is mentioned by C.~Soul\'e in \cite[3.4]{sou81}.

(2)
F. Calegari \cite{Cal05} (see also \cite{Beu06}) proved that,
for $ p=2 $ and 3, $ \z_p(3) $, which in those cases equals  $ \Lp p 3 \omega_p^{-2} {\Q} $, is irrational. 
(More results along these lines are described in Remark~\ref{newrem}(3).)

(3)
Parts (1) and (2) of Conjecture~\ref{pborel} hold when $ k $ is a totally real Abelian number field,
and in fact the corresponding parts of a much stronger conjecture that we shall describe in Section~\ref{motivicsection}
hold for cyclotomic fields (see Proposition~\ref{abelprop}).

(4)
We numerically verified part (3) of Conjecture~\ref{pborel}, and its
more refined version Conjecture~\ref{motiveconjecture}(4) below, in
certain cases; see Remark~\ref{dirichletcalculations} and
Section~\ref{examplesection}.
\end{remark}

\section{A motivic version of the conjecture}
\label{motivicsection}

If $ E $ is any extension of $ \Q $, and $ k/\Q $ is finite and Galois with Galois group $ G $,
then we let $ M^E = E \tensor_\Q k $ and $ \KnE k = E \tensor_\Q \Kn k $, which  are $ E[G] $-modules.
The goal of this section is to refine Conjecture~\ref{pborel} to a conjecture for Artin motives with
coefficients in $ E $, or equivalently, idempotents in the group ring $ E[G] $, when $ E $ is a number field.

\begin{definition}
For an idempotent $ \pi $ in $ E[G] $ we let $ M_\pi^E = \pi M^E $ and $ K_{2n-1} (M_\pi^E) = \pi \KnE k $.
\end{definition}

Now fix an embedding $ \phi_\infty : k \to \C $.
The pairing $ G \times k \to \C $ mapping $ (\s,a) $ to $ \phi_\infty(\s(a)) $ leads to an $ E $-bilinear
pairing
\begin{equation}
\begin{aligned}\label{discpairing}
E[G] \times M^E & \to E \tensor_\Q \C
\cr
(e \s, e' \tensor a ) & \mapsto e e' \tensor \phi_\infty(\s(a))
\,.
\end{aligned}
\end{equation}
and we consider its restriction
\begin{equation*}
(\,\cdot\,,\,\cdot\,)_\infty : E[G]\pi \times M_\pi^E \to E \tensor_\Q \C
\,.
\end{equation*}
Similarly, replacing $ \phi_\infty $ with a fixed embedding $ \phi_p : k \to \Qpbar $ we obtain
\begin{equation*}
(\,\cdot\,,\,\cdot\,)_p : E[G]\pi \times M_\pi^E \to E \tensor_\Q F
\,.
\end{equation*}
where 
$ F \subset \Qpbar $ is the topological closure of $ \phi_p(k) $,
which is independent of $ \phi_p $ since $ k/\Q $ is Galois.

\begin{lemma}\label{dimlemma}
Let $ \pi $ in $ E[G] $ be an idempotent.  Then $ \dim_E(M_\pi^E) = \dim_E(E[G]\pi) $.
\end{lemma}

\begin{proof}
Since $ \pi^2 = \pi  $, \eqref{discpairing} is identically 0 on 
$ E[G]\pi \times  M_{1-\pi}^E $ and
$ E[G](1-\pi) \times  M_\pi^E $.
But the determinant of this pairing is, up to multiplication by
$ E^* $, equal to $ D_k^{1/2,\infty} $, hence non-zero.  This implies the lemma.
\end{proof}

We now introduce pairings similar to $ (\cdot,\cdot)_\infty $ and $
(\cdot,\cdot)_p $ but replacing $ M_\pi^E $ with $ \KMn M_\pi^E  $.
If we denote the map $ \Kn k \to \Kn {\C} $ induced by $ \phi_\infty $
by $ \phi_{\infty*} $ and let $ \Phi_\infty $
be the composition with the Beilinson regulator map $ \reg_\infty $ for $ \C $,
\begin{equation*}
\Phi_\infty = \reg_\infty \circ \phi_{\infty*} : \Kn k \to \Kn {\C} \to \R(n-1) \subset \C
\,,
\end{equation*}
then the pairing $ G \times \Kn k \to \C $ given by mapping $ (\s,\a) $ to $ \Phi_\infty(\s(\a)) $ gives
rise to an $ E $-bilinear pairing
\begin{equation}\label{EKpairing}
\begin{aligned}
E[G] \times \KMn M^E & \to E \otimes_\Q \C
\\
(e\s, e' \tensor \a ) & \mapsto  e e' \tensor \Phi_\infty(\s(\a))
\end{aligned}
\end{equation}
and we consider its restriction
\begin{equation*}
[\,\cdot\,,\,\cdot\,]_\infty : E[G]\pi \times \KMn M_\pi^E \to E \tensor_\Q \C
\,.
\end{equation*}

By Lemma~\ref{synlemma} and the definition of $ F $ the syntomic regulator gives us
\begin{equation*}
\reg_\syn : \Kn {\OO_F} \to F
\,.
\end{equation*}
If we write $ \phi_{p*} $ for the composition $ \Kn k \iso \Kn {\O} \to \Kn {\OO_F} $, with the second
map induced by $ \phi_p : k \to F $, and let $ \Phi_p $ be the composition
\begin{equation*}
\reg_\syn \circ \phi_{p*} : \Kn k \to \Kn {F} \to F
\,,
\end{equation*}
then we can similarly obtain a pairing
\begin{equation*}
[\,\cdot\,,\,\cdot\,]_p : E[G]\pi \times \KMn M_\pi^E \to E \tensor_\Q F
\,.
\end{equation*}

We now fix ordered $ E $-bases of $ E[G]\pi $, $ M_\pi^E $, and $ \KMn M_\pi^E  $.

\begin{definition}\label{discdef}
For $ *=\infty $ or $ p $
we let $ D(M_\pi^E)^{1/2,*} $ be the determinant of the pairing~$ (\cdot,\cdot)_* $,
computed with respect to our fixed bases of $ E[G]\pi $ and $ M_\pi^E $.
\end{definition}

Note that $ D(M_\pi^E)^{1/2,\infty} $ is non-zero by the proof of Lemma~\ref{dimlemma}, and the same
argument works for $ D(M_\pi^E)^{1/2,p} $.

\begin{definition}\label{regdef}
If $ \dim_E(E[G]\pi) $ equals $ \dim(\KMn M_\pi^E ) $ then
for $ *=\infty $ or $ p $
we let $ R_{n,\infty}(M_\pi^E) $ be the determinant of the pairing~$ [\cdot,\cdot]_* $,
computed with respect to our fixed bases of $ E[G]\pi $ and $ \KMn M_\pi^E $.
\end{definition}

For future use we prove the following.

\begin{lemma}\label{Qpquot}
If $ \dim_E(E[G]\pi) = \dim(\KMn M_\pi^E ) $ then
\begin{enumerate}
\item
$ R_{n,\infty}(M_\pi^E) / D(M_\pi^E)^{1/2,\infty} $
is independent of the basis of $ E[G]\pi $, of $ \phi_\infty $, and lies in $ E \tensor \R $;
\item
$ R_{n,p}(M_\pi^E) / D(M_\pi^E)^{1/2,p} $
is independent of the basis of $ E[G]\pi $, of $ \phi_p $, and lies in $ E \tensor \Qp $.
\end{enumerate}
\end{lemma}

\begin{proof}
We prove the second statement, the proof of the first being entirely similar.
Choosing a different $ E $-basis of $ E[G]\pi $ corresponds to letting an $ E $-linear transformation
act on $ E[G]\pi $, and in the given quotient the resulting determinant cancels.
Since $ k/\Q $ is Galois we get all possible embeddings of $ k $ into $ \Qpbar $ by replacing $ \phi_p $ with $ \phi_p \circ \s $
for $ \s $ in $ G $.   For both $ R_{n,p}(M^E) $ and $ D^{1/2,p}(M^E) $ this corresponds to letting $ \s $ act on $ E[G]\pi $,
and the resulting determinant cancels as before.  That the quotient lies in $ E \tensor \Qp $ follows
because it lies in $ E \tensor F $ and we have just proved that it is invariant under
$ \Gal(F/\Qp) $ by Lemma~\ref{synlemma}.
\end{proof}

We now investigate when the two dimensions in Definition~\ref{regdef} are equal.  The answer is given by Proposition~\ref{whenconj} below, but
we need some preliminary results.

\begin{proposition}\label{kerisreal}
If $ k/\Q $ is a finite Galois extension with Galois group $ G $ and $ \pi $ is an idempotent in $ E[G] $ then
$ \dim_E (\KMn M_\pi^E ) \le \dim_E (E[G] \pi) $.  Equality holds precisely when any $($hence every$)$ $ \t $
in the image in $ G $ of the conjugacy class of complex conjugation in $ \Gal(\Qbar/\Q) $ acts
by $ (-1)^{n-1} $ on $ E[G]\pi $.
\end{proposition}

\begin{proof}
For the statement we can first replace
$ E $ by a finitely generated subfield since $ \pi $ contains only finitely many elements of $ E $, next
embed $ E $ into $ \C $ and use this embedding to enlarge $ E $ to $ \C $.
So we may assume that~$ E=\C $.

According to Theorem~\ref{boreltheorem} the pairing~\eqref{EKpairing}
gives an injection
\begin{equation}\label{Rinjection}
\begin{aligned}
\R \tensor_\Q \Kn k  & \to \R[G]^\vee = \Hom_\R(\R[G], \R)
\cr
\a & \mapsto f_\a \text{ with } f_\a(\s) = \frac{1}{(2 \pi i)^{n-1}} \Phi_\infty(\s(\a))
\end{aligned}
\end{equation}
and by extending the coefficients we get an injection $ \C \tensor_\Q \Kn k  \to \C[G]^\vee $.
The image of $ \pi (\Kn k \tensor_\Q \C) $ under the last map vanishes on $ \C[G](1-\pi) $ since $ (1-\pi) \pi = 0 $
so that $ \pi \Kn k \tensor_\Q \C $ injects into $ (\C[G]/\C[G](1-\pi))^\vee \iso (\C[G]\pi)^\vee $, which proves our
first inequality.

As for equality, we know by Theorem~\ref{boreltheorem} that~\eqref{Rinjection} has as its image the subspace
of $ \R[G]^\vee $ where, under the action of $ \Gal(\Qbar/\Q) $ on
$ \R[G]^\vee $ via $ (\tilde \s f)(\s) = f(\tilde \s^{-1}\s) $, the conjucagy class of complex conjugation in $ \Gal(\Qbar/\Q) $
acts as multiplication by $ (-1)^{n-1} $.
The same will therefore hold with complex coefficients.  Since any element $ \t $ in the conjugacy class
of complex conjugation has order~1 or~2, $ \C[G] \pi $ decomposes
into eigenspaces for the eigenvalues $ \pm 1 $ and the desired equality can only hold if
$ \t $ acts as multiplication by $ (-1)^{n-1} $ on all of $ \C[G] \pi $.
\end{proof}

We now determine precisely when the equality of dimensions as in Proposition~\ref{kerisreal} can occur,
and for this we need a preliminary result.

\begin{prop}\label{CMprop}
Let $ \psi $ be a representation of $ \Gal(\Qbar/\Q) $ over $ \Qbar $ that factorizes through the Galois group of a finite Galois extension
of $ \Q $ and for which $ \psi(\t) $ acts as multiplication by $ (-1)^{n-1} $ for any $ \t $ in the conjugacy
class of complex conjugation in $ \Gal(\Qbar/\Q) $.
Then the fixed field of $ \Ker(\psi) $ is a finite Galois extension of $ \Q $ that is totally real if
$ n $ is odd, and CM if $ n $ is even.
\end{prop}

\begin{proof}
That the fixed field $ k $ of $ \Ker(\psi) $ is a finite Galois extension of $ \Q $ is clear.
When $ n $ is odd $ \psi(\t)=1 $ for any $ \t $ in the conjugacy class
of complex conjugation so that $ k $ is totally real.
For even $ n $ we let $ \omega $ be the composition
$ \Gal(\Qbar/\Q) \to \Gal(\Q(\mu_4)/\Q) \rightiso \{\pm1\} \subset \Qbar^* $.
Then $ \Ker(\psi\omega) $ contains $ \t $ and $ \Ker(\psi) \cap \Ker(\omega) $,
hence its fixed field is a totally real Galois extension $ k^+ $ of $ \Q $, contained in $ k(\mu_4) $.
Similarly $ k $ is contained in the CM field $ k^+(\mu_4) $ so, since $ k $ is not totally real, it must
be CM.
\end{proof}

\begin{remark}\label{CM}
For a CM field $ k $ with $ k/\Q $ Galois, its maximal totally real subfield is Galois
over $ \Q $ and is the fixed field of an element of order two in the centre of $ \Gal(k/\Q) $,
which we shall refer to as the complex conjugation of $ k $.
\end{remark}

\begin{proposition}\label{whenconj}
Let $ k/\Q $ be a finite Galois extension with Galois group $ G $, $ E $ any extension of $ \Q $, and
$ \pi $ an idempotent of $ E[G] $.  Let $ k' $ be the fixed field of the kernel of the
representation of $ G $ on $ E[G] \pi $.
Then for $ n \ge 2 $ the equality $ \dim_E(E[G]\pi) = \dim_E(\KMn M_\pi^E ) $ holds precisely in the following
cases:
\begin{enumerate}
\item
$ k' $ is totally real and $ n $ is odd;

\item
$ k' $ is a CM field, $ n $ is even, and the complex conjugation of $ k' $ acts on $ E[G] \pi $ as multiplication by $ -1 $.
\end{enumerate}
\end{proposition}

\begin{proof}
From Propositions~\ref{kerisreal} and~\ref{CMprop} we see that there cannot
be any other cases.  Conversely, by Proposition~\ref{kerisreal} equality holds in both.
\end{proof}

We now recall and introduce some terminology for later use.

\begin{definition}\label{mindef}
Let $ G $ be a finite group and $ E $ any extension of $ \Q $.
\begin{enumerate}
\item
If $ \pi $ is a central idempotent of $ E[G] $ such that $ E[G] \pi $ is a minimal (non-zero) 2-sided ideal of $ E[G] $
then $ \pi $ is a primitive central idempotent.

\item
If $ \pi $ is a primitive central idempotent of $ E[G] $
such that $ E[G] \pi \iso M_m(E) $ as $ E $-algebras for some $ m $ then we call a primitive idempotent corresponding
to $ \pi $ any element in $ E[G] \pi \subseteq E[G] $ that maps to a matrix in $ M_m(E) $ that is conjugate
to a matrix with~1 in the upper left corner and 0's elsewhere.  (Since all automorphisms of $ M_m(E) $
as $ E $-algebra are inner this is independent of the isomorphism $ E[G] \pi \iso M_m(E) $.)
\end{enumerate}
\end{definition}

\begin{remark}{\label{minrem}}
(1)
If $ E[G] \pi \iso M_m(E) $ then any idempotent in $ E[G] \pi $ can be written as a sum of 
orthogonal primitive idempotents corresponding to $ \pi $ as one sees immediately by diagonalizing the
matrix $ A $ that $ \pi $ maps to, which satisfies $ A^2=A $.

(2)
If $ G = \Gal(k/\Q) $  and $ \pi $ is a primitive central idempotent of $ E[G] $
then the dimensions of $ E[G] \tilde \pi $ and $ \tilde \pi \KnE k $ are equal for some non-zero idempotent
$ \tilde \pi $ in $ E[G] \pi $ if and only
if the same holds for one (hence any) primitive idempotent corresponding to $ \pi $.
Indeed, the dimensions for $ \pi $ do not change if we replace it with a conjugate in $ E[G] $, and they
add up for sums of orthogonal idempotents.
\end{remark}

We now introduce \Lfunctions, both classical and \padic, in the following context.  Let $ E $ be a finite extension
of $ \Q $, $ k $ a number field, and consider a character $ \psi $ of a representation of $ G = \Gal(\kbar/k) $ on a finite dimensional
$ E $-vector space~$ V $ that factorizes through $ \Gal(\kbar/k') $ for some finite extension $ k' $
of $ k $.   Then for every embedding $ \s : k \to \C $ we have the Artin \Lfunction $ \L s \s(\psi) k $,
with for every prime $ \P $  of $ k $ the reciprocal $ \Eul_\P(s, \s(\psi) ,k) $ of the Euler factor
corresponding to $ \P $.
Under the natural isomorphism $ E \tensor_\Q \C \rightiso \oplus_\s \C_\s  $
the $ \L s \s(\psi) k $ correspond to a canonical $ E\tensor_\Q \C $-valued \Lfunction that we
denote by $ \L s {\psi\tensor\id} k $.
Similarly, for every prime $ \P $ of $ k $ we have an $ E\tensor_\Q\C $-valued $ \Eul_\P( s, \psi\tensor\id, k) $
corresponding to the $ \Eul_\P(s,\s(\psi),k) $.

We now move on to the \padic \Lfunctions, and assume that $ k $ is totally real, $ p $ a
prime number, $ a $ any integer and $ \omega_p $ the Teichm\"uller character $ \Gal(\Qbar/k) \to \Q_p^* $
\eqref{TM}.
If $ \t : E \to \Qpbar $ is an embedding and $ k_{\t(\psi)\omega_p^a} $ is totally real
then from Section~\ref{conjecturesection} we have the \padic \Lfunction $ \Lp p s \t(\psi)\omega_p^a k $,
which is not identically zero.
In this case, using the natural isomorphism $ E\tensor_\Q\Qpbar \rightiso \oplus_\t \overline{\Q}_{p,\t} $,
they give us an $ E\tensor_\Q \Qpbar $-valued \padic \Lfunction on $ \Zp $ or $ \Zp\setminus\{1\} $
that we denote by $ \Lp p s \psi\tensor\omega_p^a k $.

\begin{lemma}\label{pLvalueslemma}
The values of $ \Lp p s \psi\tensor\omega_p^a k $ are in $ E\tensor_\Q \Qp $.
\end{lemma}

\begin{proof}
Using Brauer induction for $ \psi $ (cf.~\eqref{brauerinduction}) it suffices to do prove this when $ \psi $
is 1-dimensional.  But then for each $ \t : E \to \Qpbar $ the function $ \Lp p s \t(\psi)\omega_p^a k $ can be described
as in~\eqref{merpL} with $ l=-1 $, from which the result is clear.
\end{proof}

\begin{remark}
It follows from Remark~\ref{interpolrem} that $ \Lp p s \psi\tensor\omega_p^a k $ satisfies
\[
\Lp p m \psi\tensor\omega_p^a k = \Eul_p(m, \psi\tensor\id, k) \L m {\psi\tensor\id} k
\]
for integers $ m\le0 $ congruent to $ 1-a $ modulo $ \phi(q) $,
where both sides lie in $ E = E\tensor_\Q \Q $ inside $ E\tensor_\Q \Qp $ and $ E\tensor_\Q\C $ respectively,
and
\[
\Eul_p(s, \psi\tensor\id, k) = \prod_{\P|p} \Eul_\P(s, \psi\tensor\id, k)
\]
the product being over all primes of $ k $ dividing~$ p $.
\end{remark}

\begin{remark}\label{idenrem}
If $ E=\Q $ we shall identify $ E\tensor_\Q \C $ with $ \C $, and write $ \L s {\psi} k $, etc.,
instead of $ \L s {\psi\tensor\id} k $, etc.  Similarly we identify $ E\tensor_\Q\Qp $ with $ \Qp $
and write $ \Lp p s {\psi}\omega_p^a k $ instead of  $ \Lp p s {\psi}\tensor\omega_p^a k $.
\end{remark}

We now have all the ingredients for the generalization and refinement of Conjecture~\ref{pborel}.
Starting with a finite Galois extension $ k/\Q $ with Galois group $ G $, $ E $ a finite extension of
$ \Q $,
and $ \pi $ an idempotent in $ E[G] $, we let $ \psi_\pi $ be the natural representation
of $ \Gal(\Qbar/\Q) $ on $ E[G]\pi $ and $ \chi_\pi $ its associated character.  If
$ \dim_E(E[G]\pi) = \dim_E(\pi \KnE k ) $ for some $ n \ge 2 $ then,
for any prime $ p $ and any embedding $ \t : E \to \Qpbar $,
$ \t(\psi_\pi) \omega_p^{1-n} $ is trivial on the conjugacy class of complex conjugation
in $ \Gal(\Qbar/\Q) $ by Proposition~\ref{kerisreal}.
In particular, $ \Q_{\t(\psi_\pi)\omega_p^{1-n}} $ is totally real
and therefore $ \Lp p n {\chi_\pi\tensor\id} {\Q} $ is not identically zero.
With $ F \subset \Qpbar $ the topological closure of $ \phi_p(k) $ as before,
using ordered bases for $ E[G]\pi $, $ M_\pi^E $ and $ \KMn M_\pi^E $,
we have $ D(M_\pi^E)^{1/2,\infty} $ and the regulator $ R_{n,\infty}(M_\pi^E) $ in $ E \tensor \C $ as well
as $ D(M_\pi^E)^{1/2,p} $ and  $ R_{n,p}(M_\pi^E) $ in $ E \tensor F $.

\begin{conjecture}\label{motiveconjecture}
With notation as above, if $ \dim_E(E[G]\pi) = \dim_E(\pi \KnE k ) $ for some $ n \ge 2 $, then
\begin{enumerate}
\item
in $ E \tensor_\Q \C $ we have
\[
\qquad\qquad
\L n {\chi_\pi\tensor\id} {\Q}  D(M_\pi^E)^{1/2,\infty} = e(n, M_\pi^E) R_{n,\infty}(M_\pi^E)
\]
for some $ e(n, M_\pi^E) $ in $ (E\tensor_\Q \Q)^* $;
\item
in $ E \tensor_\Q F $ we have
\[
\qquad\qquad
\Lp p n {\chi_\pi}\tensor\,\omega_p^{1-n} {\Q} D(M_\pi^E)^{1/2,p}  = e_p(n,M_\pi^E) \, \Eul_p(n, \chi_\pi\tensor\,\id,\Q) R_{n,p}(M_\pi^E)
\]
for some $ e_p(n,M_\pi^E) $ in $ (E\tensor_\Q \Q)^* $;
\item
in fact, $ e_p(n, M_\pi^E) = e(n, M_\pi^E) $;
\item
$ \Lp p n {\chi_\pi}\tensor\omega_p^{1-n} {\Q} $ and $ R_{n,p}(M_\pi^E) $ are units in $ E \tensor_\Q \Qp $
and $ E \tensor_\Q F $ respectively.
\end{enumerate}
\end{conjecture}

\begin{remark}\label{newremark}
(1)
One sees as in the proof of Lemma~\ref{Qpquot} that the validity of each part of the conjecture is independent of
the chosen bases of $ E[G]\pi $, $ M_\pi^E $ and $ \KMn M_\pi^E $.

(2)
Since $ E\tensor_\Q \C \rightiso \oplus_\s E_\s $, where the sum is over all embeddings $ \s : E \to \C $,
an identity in $ E \tensor_\Q \C $ is equivalent to the corresponding identities for all such embeddings.
The same holds if we replace $ \C $ with $ \Qpbar $.
\end{remark}

In Remark~\ref{BBmotiverem} we shall make explicit how part (1) of this conjecture is equivalent with Beilinson's
conjecture \cite[Conjecture~3.4]{bei85} for an Artin motive associated with $ \pi $.
First, we make various remarks about its dependence
on $ E $, etc., and on its relation with Conjecture~\ref{pborel}.

\begin{remark}\label{newrem}
(1)
An equivalent for Conjecture~\ref{motiveconjecture}(1) can be formulated for any idempotent $ \pi $ if $ k/\Q $ is any
finite Galois extension (see the end of Section~\ref{practical}),
in which case it is a conjecture by Gross (see \cite[p.~210]{Neu88}).
In that case it was proved by Beilinson if the action of $ G=\Gal(k/\Q) $ on $ E[G]\pi $
is Abelian (see loc.\ cit.).  One can deduce that it holds for any $ \pi $ in $ E[G] $ if this action factors through
$ S_3 $ or $ D_8 $ (see Proposition~\ref{grossprop}).

(2)
If the action of $ G $ on $ E[G] $ is Abelian then parts (1)-(3) of Conjecture~\ref{motiveconjecture}
also hold for any $ \pi $ to which the conjecture applies (see Remark~\ref{abelrem}).

(3) Extending and simplifying earlier work by F. Calegari \cite{Cal05}, F. Beukers in \cite{Beu06}
proved that $  \z_p(2) $ is irrational when $ p=2 $ or 3.  This value equals $ \Lp p 2 \chi\omega_p^{-1} {\Q} $
with $ \chi $ the primitive character on $ (\Z/4\Z)^* $ for $ p=2 $, and the primitive character on $ (\Z/3\Z)^* $ for $ p=3 $.
Moreover, if $ \chi $ is the odd primitive character on $ (\Z/8\Z)^* $ then he also shows that $ \Lp 2 2 \chi\omega_2^{-1} {\Q} $ is irrational.
It follows that the conjecture holds in full for $ \pi $ corresponding to the non-trivial representation
of $ \Gal(k/\Q) $ and $ n=2 $ when
$ (k,p) $ is one of $ (\Q(\sqrt{-1}),2) $, $ (\Q(\sqrt{-2}),2) $ and $ (\Q(\sqrt{-3}),3) $.

(4) We have verified numerically that part~(4) of the conjecture holds in certain cases;
see Remark~\ref{dirichletcalculations} as well as Section~\ref{examplesection}.
\end{remark}

\begin{remark}\label{consequencesremark}
(1)
If $ \pi = \pi_1 + \cdots + \pi_m $ with $ \pi_i^2 = \pi_i $ and $ \pi_i \pi_j = 0 $
when $ i \neq j $, then the conjecture for $ \pi $ is implied by the conjecture for all $ \pi_i $ because
$ D(\oplus_i M_i^E)^{1,2/\infty} = \prod_i D(M_i^E)^{1/2,\infty} $, etc., as one easily sees by using bases.

(2)
If $ \pi $ is in $ E[G] $ and $ E' $ is an extension of $ E $, then we may view $ \pi $ as an element
of $ E'[G] $ as well, and the conjectures for $ M_i^E $ and $ M_i^{E'} $ are equivalent:
we can use the same bases over $ E' $ as over $ E $, so that $ D(M_i^E)^{1/2,\infty} =  D(M_i^{E'})^{1/2,\infty} $
in $ E \tensor \C \subseteq E' \tensor \C $, and the same holds for all the other ingredients (including
the $ L $-functions).

(3)
By comparing bases one sees immediately that
if $ \pi $ and $ \pi' $ are conjugate under the action of $ E[G]^* $ then
the conjectures for $ \pi $ and for $ \pi' $ are equivalent.

(4)
If, for a primitive central idempotent $ \pi_i $ of $ E[G] $, $ E[G] \pi_i \iso M_m(E) $ for some $ m $,
and $ \pi $ is a primitive idempotent corresponding to $ \pi_i $, then the conjecture for $ \pi_i $ is implied
by the conjecture for $ \pi $.
Indeed, we can decompose $ \pi_i $ into a sum of orthogonal primitive idempotents
as in Remark~\ref{minrem}(1), and the truth of the conjecture for $ \pi_i $ is implied by its truth
for each such primitive idempotent.  But all such primitive idempotents are conjugate to $ \pi $
hence part (3) above applies.
\end{remark}

\begin{remark}\label{galoisdescent}
If $ k/\Q $ is a Galois extension with Galois group $ G $ and $ H $ a subgroup of $ G $ with fixed field $ k^H $,  then
$ K_m(k)_\Q^H = K_m(k^H)_\Q  $, a result known as Galois descent.
With $ \pi_H = |H|^{-1} \sum_{h \in H} h $, an idempotent in $ E[G] $, this implies that $ \pi_H K_m(k) = K_m(k^H) $.
\end{remark}

\begin{remark}\label{sameconj}
(1)
Let $ k/\Q $ be a Galois extension with Galois group $ G $.  If $ N $ is a normal subgroup of $ G $ corresponding to
$ k' = k^N $ and $ \pi_N = |N|^{-1} \sum_{h \in N} h $, then $ \pi_N^2 = \pi_N $ and the natural
map $ \phi : E[G] \to E[G/N] $ has kernel $ E[G] (1-\pi_N) $ and induces an isomorphism $ E[G] \pi_N \to E[G/N] $.
Indeed, it is clear that $ \pi_N $ is central in $ E[G] $ since $ N $ is normal in $ G $ and that $ 1-\pi_N $
is in the kernel of $ \phi $.  Also, since $ N $ acts trivially on $ E[G] \pi_N $ this is an $ E[G/N] $-module
generated by one element, so its dimension over $ E $ cannot be bigger than $ |G/N| $.  Since $ \phi $
is obviously surjective our claims follow.  Therefore in this situation for any idempotent $ \pi $ in $ E[G/N] $ there
is a canonical idempotent $ \tilde \pi $ lifting $ \pi $ to $  E[G] \pi_N $, and the natural map $ E[G] \tilde \pi \to E[G/N] \pi $
is an isomorphism of $ E[G] $- and $ E[G/N] $-modules.

Then the statements of Conjecture~\ref{motiveconjecture} for $ \pi $ in $ E[G/N] $ or for $ \tilde \pi $ in $ E[G] $ are equivalent.  Namely,
$$
\tilde \pi \KnE k = \tilde \pi \pi_N \KnE k =  \tilde \pi \KnE k' = \pi \KnE k'
$$
inside $ \KnE k $ so that we can use the same bases for either side.   The same holds for $ \tilde \pi (E\tensor k) $
and $ \pi (E \tensor k') $ inside $ E \tensor k $.
Moreover, $ E[G]\tilde\pi $ is the pullback to $ G $ of the $ G/N $-representation $ E[G/N] \pi $,
so that $ \L s  {E[G]\tilde\pi} {\Q} = \L s {E[G/N]\pi} {\Q} $ and similarly for the \padic \Lfunctions.

(2)
If $ k $ is a totally real number field $ k $ let $ \tilde k $ be its (totally real) Galois closure over $ \Q $.
Then $ k = \tilde k ^H $ for some subgroup $ H $ of $ G = \Gal(\tilde k/\Q) $, $ \pi = |H|^{-1} \sum_{h \in H} h $
is an idempotent in $ \Q[G] $, and for $ n \ge 2 $ odd Conjecture~\ref{motiveconjecture} for $ \pi $ is equivalent to Conjecture~\ref{pborel}
for $ k $.  Namely,  $ \Kn k = \pi \Kn \tilde k \subseteq \Kn \tilde k $ by Remark~\ref{galoisdescent} and $ \pi \tilde k = k $
so that we can use the same $ \Q $-bases in each case.  Moreover, $ \Q[G] \pi \iso \Ind H G (1_H) = \Q[G] \tensor_{\Q[H]} 1_H $
with $ 1_H $ the trivial 1-dimensional representation of $ H $, as one easily sees by mapping $ \b $
in $ \Q[G] \pi $ to $ \b \tensor v $ and $ \sum_\s a_\s \s \tensor(\lambda v) $ to $ \lambda \sum_\s a_\s \s \pi $, $ \{v\} $ being a basis of $ 1_H $.
By well-known properties of Artin \Lfunctions \cite[Prop.~10.4(iv)]{Neu99}
this implies that $ \z_k(s) = \L n {\Q[G]\pi} {\Q} $ and similarly for the
\padic \Lfunctions.
\end{remark}

\begin{remark}\label{BBmotiverem}
We make the relation between Conjecture~\ref{motiveconjecture}(1) and Beilinson's conjecture 
\cite[Conjecture~3.4]{bei85} for (Artin)
motives  explicit since the relation between the two elements of $ E^* $ involved
also shows up very explicitly in our computer calculations, suggesting that the element for the formulation
at $ s=1-n $ is simpler than for our formulation at $ s=n $ (see Remark~\ref{constantremark}).
We provide some details since we could not find a detailed enough reference in the literature.

With notation as in Conjecture~\ref{motiveconjecture} and~\eqref{oldfuneq} we have
\begin{equation}\label{newfuneq}
\begin{aligned}
\L 1-s {\chi_\pi^\vee\tensor\id} {\Q} 
& =
W(\chi_\pi^\vee\tensor\id) C(\chi_\pi^\vee\tensor\id)^{s-\frac12} (2 (2\pi)^{-s}\Gamma(s))^{\chi^\vee(1)}
\\
& \qquad\times
 (\cos(\pi s/2))^{n^+}  (\sin(\pi s/2))^{n^-} \L s {\chi_\pi\tensor\id} {\Q} 
\,,
\end{aligned}
\end{equation}
with $ \chi_\pi^\vee $ the dual character of $ \chi_\pi $,
$ W(\chi_\pi^\vee\tensor\id) $ in $ E\tensor\C $, and
$ C(\chi_\pi^\vee\tensor\id) = C(\chi_\pi\tensor\id)$ in $ \Q^* $.
If $ m=\dim_E(E[G]\pi) $, then $ n^+ = m $ and $ n^-=0 $ for $ n $ odd, and $ n^-=m $ and $ n^+=0 $ for
$ n $ even.  In either case $ \L s {\chi_\pi^\vee\tensor\id} {\Q} $ has a zero of order $ m $ at $ s=1-n $.
Moreover
\begin{equation}\label{edef}
W(\chi_\pi^\vee\tensor\id) C(\chi_\pi^\vee\tensor\id)^{1/2} = e_n i^{m(n-1)} D(M_\pi^E)^{1/2,\infty}
\end{equation}
for some $ e_n $ in $ E^* = (E\tensor_\Q \Q)^* $ by \cite[Propositions~5.5 and 6.5]{Del77} since 
$ D^{1/2,\infty} $ can be taken to be the same for $ M_\pi^E $ and the associated determinant representation
(cf.\ \cite[p.360]{Jan88b}).
Hence the first non-vanishing coefficient in its Taylor expansion
around $ s=1-n $, $ \Lsharp {1-n} {\chi_\pi^\vee\tensor\id} {\Q} $, equals
\begin{equation*}
\delta_n e_n (2\pi i)^{m(1-n)} ((n-1)!/2)^m C(\chi_\pi^\vee\tensor\id)^{n-1} D(M_\pi^E)^{1/2,\infty} \L n {\chi_\pi\tensor\id} {\Q} 
\end{equation*}
with $ \delta_n = (-1)^{m(n-1)/2} $ when $ n $ is odd, and $ \delta_n = (-1)^{m(n+2)/2} $ when $ n $ is even.
In particular, Conjecture~\ref{motiveconjecture}(1) is equivalent with
\begin{equation}\label{newe}
\Lsharp {1-n} {\chi_\pi^\vee\tensor\id} {\Q} 
=
\delta_n e_n ((n-1)!/2)^m C(\chi_\pi^\vee\tensor\id)^{n-1} e(n,M_\pi^E) \widetilde R_{n,\infty}(M_\pi^E)
\end{equation}
with the renormalized regulator $ \widetilde R_{n,\infty}(M_\pi^E) = (2\pi i)^{m(1-n)}R_{n,\infty}(M_\pi^E) $.

Let us compare this with Beilinson's conjecture for a motive associated with $ \pi $.
We associate motives covariantly to smooth projective varieties over $ \Q $ 
as in \cite[\S~2.4]{bei85}.
The Galois group $ G $ acts on the left on $ k $, hence on the right on $ \spec(k) $, and we let 
$ \art $ be the motive corresponding to $ \pi $ under this action.
Then $ G $ acts on the left on the cohomology theories on $ \spec(k) $
as well as its  $ K $-theory, and after tensoring with $ E $
the corresponding groups for $\art$ are the images under $\pi$.
Thus, the relevant motivic cohomology of $\art$ is $\hm^1(\mathcal{M}_{\pi/\Z},\Q(n))= K_{2n-1} (M_\pi^E) $.

We shall need the non-degenerate $ E $-bilinear pairing
\begin{equation*}
\begin{aligned}
E[G] \times \bigoplus_{\t: k \to \C} E & \to E
\\
(\smallsum_\s a_\s \s,  (b_\t)_\t ) & \mapsto \smallsum_\s a_\s b_{\phi_\infty\circ\s^{-1}}
\,.
\end{aligned}
\end{equation*}
It factors through the left $ E[G] $-action on $ \oplus_\t E $
(given by $ \s((b_\t)_\t) = (b_{\t\circ\s^{-1}})_\t $),
hence is trivial on $ E[G](1-\pi)\times \pi(\oplus_\t E) $
and $ E[G]\pi\times (1-\pi)(\oplus_\t E) $.
We therefore obtain a non-degenerate $ E $-bilinear pairing
\begin{equation}\label{<>}
\langle \,\cdot\, , \,\cdot\, \rangle : E[G]\pi \times \pi(\oplus_\t E) \to E
\end{equation} 
that identifies $ E[G]\pi $ and $ \pi(\oplus_\t E) $ as dual $ E[G] $-modules.

Tensoring~\eqref{HDid} with $ E $ and applying $\pi$ we get
\begin{equation*}
    \hd^1(\artR,\R(n)) \isom  \pi \bigl( \bigoplus_{\tau: k\to \C} E\otimes \R(n-1) \bigr)
\,.
\end{equation*}
Note that the left hand side is a subspace of the right hand side by~\eqref{HDid}, but because
the regulator map~\eqref{Bregulator} tensored with $ \R $ is injective,
and $ \dim_E(K_{2n-1}(M_\pi^E)) = \dim_E(E[G]\pi) = \dim_E(\pi(\oplus_\t E)) $ by our assumption on $ \pi $
and since $  D(M_\pi^E)^{1/2,\infty} \ne 0 $ (as was noticed right after Definition~\ref{discdef}),
equality must hold.
For the Beilinson regulator we therefore have to compare
the two $ E $-structures on $ \det\hd^1(\artR,\R(n)) $ coming from
Betti cohomology,
\begin{equation*}
\hb^0(\artR,\Q(n-1)) \isom  \pi (\bigoplus_{\tau: k\to \C} E\otimes \Q(n-1)) 
\subseteq     \pi \bigl(\bigoplus_{\tau: k\to \C} E\otimes \R(n-1) \bigr) 
\,,
\end{equation*}
and that induced by the Beilinson regulator map~\eqref{Bregulator},
\begin{equation*}
\hm^1(\artR,\Q(n)) \to \hd^1(\artR,\R(n)) 
= \pi \bigl(\bigoplus_{\tau: k\to \C} E\otimes \R(n-1)\bigr) 
\,.
\end{equation*}
Choosing an $ E $-basis of $ \pi(\oplus_\t E) $, and multiplying it by $ (2\pi i)^{n-1} $ to obtain an $ E $-basis
of $ \pi (\oplus_\t E\tensor \Q(n-1)) $, it is easy to see that
Beilinson's regulator $ R_\Bei $ for $ \art $ satisfies
\begin{equation*}
R_\Bei = (2\pi i)^{m(1-n)} \det [\,\cdot\,,\,\cdot\,]_\infty / \det \langle\,\cdot\,,\,\cdot\,\rangle 
= \widetilde  R_{n,\infty}(M_\pi^E) / \det \langle\,\cdot\,,\,\cdot\,\rangle 
\end{equation*}
where all determinants are computed using the chosen $ E $-bases.

Finally, we compare \Lfunctions.  We have 
\begin{equation*}
 H_{\textup{et}}^0(\spec(k) \otimes \Qbar,\Q_l) \iso \bigoplus_{\t:k\to\Qbar} \Q_l
\iso \Q_l \tensor_\Q \biggl( \bigoplus_{\t:k\to\Qbar} \Q \biggr)
\end{equation*}
so that 
$H_{\textup{et}}^0(\art,\Q_l) \iso \Q_l \tensor \pi (\oplus_\t E)  $ as $ \Q_l\tensor E[G] $-modules.
The $ \Q_l $ plays no role, and as in~\eqref{<>} we see that
$ \pi(\oplus_{\t:k\to\Qbar} E ) $ is dual to $ E[G] \pi $ as $ E[G] $-module.
As the motivic \Lfunction uses
the geometric rather than the arithmetic Frobenius (cf.~\cite[p.26]{Kat94}),
we obtain that $ L(s,\art) $ of \cite[\S~3]{bei85} is equal to $ \L s {\chi_\pi\tensor\id} {\Q} $.
If $ \theta : E[G] \to E[G] $ is the $ E $-linear involution obtained by replacing each $ \s $ in $ G $ with $ \s^{-1} $,
then we need to consider $ \mathcal{M}_{\pi}^0 = \mathcal{M}_{\theta(\pi)} $ instead of $ \art $.
But the map $ E[G]\pi \times E[G]\theta(\pi) \to E $ mapping $ (\a,\b) $ to the coefficient of the neutral
element of $ G $
in $ \theta(\b)\a $
is easily seen to identify $ E[G]\theta(\pi) $ and $ E[G]\pi $ as dual $ E[G] $-modules, so that
$ L(s,\art^0) = \L s {\chi_\pi^\vee\tensor\id} {\Q} $.
\end{remark}

\section{More explicit $ K $-groups and regulators}
\label{ksection}

In this section we first describe an inductive procedure that conjecturally gives $ \Kn k $
($ n \ge 2 $) for any number field $ k $.  It is originally due to Zagier \cite{Zag91}, but we essentially
give a reformulation by Deligne \cite{del:imd}.
We also describe results concerning Conjecture~\ref{motiveconjecture} when the action of $ G $ on $ E[G]\pi $
is Abelian.

In order to describe Zagier's conjecture we need the functions
\begin{equation}\label{Lindef}
\Li_n(z) = \sum_{k \geq 1} \frac{z^k}{k^n}
\qquad (n \geq 0)
\end{equation}
for $ z $ in $ \C $ with $ |z| < 1 $ if $ n = 0 $ or 1, and $ |z|\le1 $ if $ n \ge 2 $.
In particular, $ \Li_1(z) $ is the main branch of $ -\log(1-z) $.
Using that $ \dd \, \Li_{n+1}(z) = \Li_n(z) \, \dd \log(z) $ they extend to multi-valued analytic functions on
$ \C\setminus\{0,1\} $.
By simultaneously continuing all $ \Li_n $ along the same path one obtains single-valued functions on $ \C\setminus\{0,1\} $
(see \cite[\S~7]{Zag91} or \cite[Remark~5.2]{dJ95}) by putting
\begin{equation}\label{Pndef}
\Pn_n(z)  = \pi_{n-1}\biggl( \sum_{j=0}^{n-1} \frac{b_j}{j!} (2\log|z|)^j \Li_{n-j}(z) \biggr)
\quad\quad (n \ge 1) \,,
\end{equation}
with $ b_j $ the $ j $-th Bernoulli number and $ \pi_{n-1} $ the projection of
$ \C=\R(n-1)\oplus\R(n) $ onto $ \R(n-1) $.
These functions satisfy $  \Pn_n(\bar z) = \overline{\Pn_n(z)}  $ as well as
\begin{equation}\label{poleq}
\Pn_n(z) + (-1)^n \Pn_n (1/z) = 0
\end{equation}
and
\begin{equation}\label{distrieq}
\Pn_n(z^m) = m^{n-1} \sum_{\z^m=1} \Pn_n(\z z)
\end{equation}
when $ m \ge 1 $ and $ z^m \ne 1 $.

We can now describe the conjecture.  Let $ d=[k:\Q] $, and for $ n \ge 2 $ let $ \tB_n(k) $ be a free
Abelian group on generators $ [x]_n^\sim $ with $ x\ne 0,1 $ in $ k $.
Define
\begin{alignat*}{1}
\widetilde P_n : \tB_n(k) & \to \left(\R(n-1)^d\right)^+
\\
[x]_n^\sim & \mapsto (P_n(\s(x)))_{\s: k \to \C}
\,,
\end{alignat*}
with
$ \left(\R(n-1)^d\right)^+ = \{ (a_\s)_\s \text{ in } \R(n-1)^d \text{ such that } a_{\bar\s} = \overline{a_\s} \} $.
Then we define inductively, for $ n \ge 2 $,
\begin{alignat*}{1}
\dd_n : \tB_n(k) & \to
 \left\{\begin{aligned}\bigwedge^2\lower 4pt\hbox{$ {}_\Z $} k^* \hbox to 30 pt{} & \text{ if $ n =2 $} \\ B_{n-1}(k) \tensor_\Z k^* & \text{ if $n>2$}\end{aligned} \right.
\\
\intertext{by}
[x]_n^\sim & \mapsto
   \left\{\begin{aligned}(1-x)\wedge x & \text{ if $ n =2 $} \\ [x]_{n-1} \tensor x & \text{ if $  n > 2 $\,,}\end{aligned} \right.
\end{alignat*}
where $ [x]_{n-1} $ denotes the class of $ [x]_{n-1}^\sim $ in $ B_{n-1}(k) $, which is defined as
$$
B_n(k) = \tB_n(k)/ \Ker(\dd_n) \cap \Ker(\widetilde P_n)
\,.
$$
There are some universal relations, one of which is that $ [x]_n + (-1)^n [1/x]_n = 0 $, a consequence
of~\eqref{poleq}.

\begin{conjecture}\label{zagierconjecture}
If $ n \ge 2 $ then
\begin{enumerate}
\item
there is an injection
$$
\dfrac{\Ker(\dd_n)}{\Ker(\dd_n) \hskip-.5pt \cap \Ker(\widetilde P_n)} \to K_{2n-1}(k)_\Q
$$
with image a finitely generated group of rank equal to $ \dim_\Q (\Kn k ) $;
\item
Beilinson's regulator map is given by $ \widetilde P_n $:
$$
\xymatrix{
\dfrac{\Ker(\dd_n)}{\Ker(\dd_n) \hskip-.5pt \cap \Ker(\widetilde P_n)} \ar[r]\ar[rd]_-{(n-1)!\widetilde P_n} & \Kn
k \ar[d]^-{\prod_{\s: k \to \C} \reg_\infty\circ\s_*} \cr
  & \left(\R(n-1)^d\right)^+
}
$$
commutes.
\end{enumerate}
\end{conjecture}

\begin{remark}
For $ n=2 $ the corresponding results were already known before Zagier made his conjecture
(\cite{bl00,susXkof}; see also \cite[\S~2]{gon:goc}).
\end{remark}

From the results in \cite[\S~5]{dJ95} one obtains the following.

\begin{theorem}\label{zagierthm}
Let $ k $ be a number field and let $ n \geq 2 $ be an integer.  Then there is an injection
$$
\Psi_n : \dfrac{\Ker(\dd_n)}{\Ker(\dd_n) \cap \Ker(\widetilde P_n)} \to \Kn k
$$
with finitely generated image, such that the diagram in Conjecture~\ref{zagierconjecture}(2)
commutes.
\end{theorem}

\begin{remark}
The part of Zagier's conjecture that remains open is the question if the rank of the image of
$ \Psi_n $ equals $ \dim_\Q(\Kn k ) $.
For a cyclotomic field $ k $ this holds for any $ n \geq 2 $ as we shall recall in Example~\ref{cyclo} below,
but it is also known for arbitrary number fields for $ n=2 $, as mentioned above, or $ n=3 $
(see \cite[\S~3]{gonXpam} and \cite[Appendix]{G3}).
\end{remark}

\begin{remark}\label{equivariance}
Although not stated explicitly in \cite{dJ95}, it is clear from the construction there that the map $ \Psi_n $
in Theorem~\ref{zagierthm}
is natural in $ k $.  In particular, if $ k/\Q $ is Galois then $ \Gal(k/\Q) $ acts on
$ \Ker(\dd_n) /\Ker(\dd_n) \cap \Ker(\widetilde P_n) $ (through the obvious action on the generators $ [x]_n^\sim $)
as well as on $ \Kn k $, and $ \Psi_n $ is $ \Gal(k/\Q) $-equivariant.
\end{remark}

\begin{example}\label{cyclo}
If $ k $ is a cyclotomic field then Zagier's conjecture is known in full.  By Theorem~\ref{zagierthm}
it suffices to see that the rank of the image of $ \Psi_n $ equals $ \dim_\Q \Kn k $.
For an $ N $-th root of unity $ \z\neq1 $, $ N[\zeta]_n $ is in
$ \Ker(\dd_n) / \Ker(\dd_n) \cap \Ker(\widetilde P_n) $.
If $ k = \Q(\z) $ for a primitive $ N $-th root of unity $ \z $
with $ N > 2 $ then $ \dim_\Q \Kn k = [k:\Q]/2 $ for $ n \geq 2 $, and 
the $ N [\z^j]_n $ with $ 1 < j < N/2 $ and $ \gcd(j,N)=1 $
have $ \R $-independent images under $ \widetilde P_n $ \cite[pp.420--422]{Zag91}
so that they are $ \Z $-independent in $ \Ker(\dd_n) / \Ker(\dd_n) \cap \Ker(\widetilde P_n) $,
and the same holds for their images under the injective map $ \Psi_n $.
For $ k=\Q $ one can easily check directly from Theorem~\ref{zagierthm} and~\eqref{Pndef}
that for odd $ n \geq 2 $ the element $ 2 \Psi([-1]_n) $ is non-zero.
\end{example}

For the \padic regulator we need to describe the \padic polylogarithms
introduced in~\cite{Col82}. We first make a choice of a branch of
the \padic logarithm. Recall that a homomorphism $\log_p: \OO_{\Cp}^*
\to \Cp$ is uniquely determined by the requirement that for $|x|_p<1$
it is given by the usual power series for $\log(1+x)$. To extend it to a
homomorphism $\log_p: \Cp^* \to \Cp$ it suffices to make an arbitrary choice
of $\log_p(p)$. Any such extension will be called a branch of the
$p$-adic polylogarithm. In what follows we fix one such choice.

Coleman first produced the $p$-adic polylogarithm out of a more
extensive theory of what is now called Coleman integration. For the
$p$-adic polylogarithm it is, however, possible to
give a more elementary and explicit theory.
To do this, consider the
class of functions $f: \Cp\setminus\{1\} \to \Cp$ that satisfy the following properties:
\begin{enumerate}
\item for any $a\in \Cp$ with $|a-1|_p=1$ we have a power series expansion
  for $f(z)$ in $z-a$ that converges for $|z-a|_p<1$;
\item for $0<|z-1|_p<1$ (resp.\ $|z|_p>1$) $ f(z) $ is  given by a polynomial in $\log_p(z-1)$
  (resp.\ $\log_p(1/z)$) with coefficients that are Laurent series in
  $z-1$ (resp.\ $1/z$), convergent for $0<|z-1|_p<1$ (resp.\ $|z|_p>1$).
\end{enumerate}
It is easy to see that differentiation sends this class surjectively onto itself.
The $p$-adic polylogarithms are functions $\Li_{n,p}(z) $ ($ n \ge 0 $)
in this class with the following properties:
\begin{enumerate}
\item $\Li_{0,p}(z) = \frac{z}{1-z}$;
\item $\Li_{n,p}(0) = 0 $ for $n\ge 0 $;
\item $ \dd \Li_{n+1,p}(z) = \Li_{n,p}(z) \frac{\dd  z}{z}$ for $n\ge 0 $;
\item for every $ n \ge 0 $ there is a power series $ g_n(v) $, convergent for $ |v|_p < p^{1/(p-1)} $,
such that $ g_n(1/(1-z)) = \Li_{n,p}(z)  - \frac{1}{p^n} \Li_{n,p}(z^p)$ when $ |z-1|_p>  p^{-1/(p-1)}$.
\end{enumerate}
Note that $1/z$ has a singularity at $0$ but $ \Li_{n,p}(z) /z$ does not,
because of the assumption that $\Li_{n,p}(0) = 0 $. In fact, it is easy to
see that $\Li_{n,p}(z)$ is given by the $p$-adically convergent power
series \eqref{Lindef} for $|z|_p<1$.

In~\cite{BdJ06} two of the authors of the present work described an
efficient algorithm for the computation of $p$-adic polylogarithms up to a prescribed precision.
This will be used in Section~\ref{examplesection} to compute the \padic regulator for $ n \ge 2 $ that we now proceed
to describe in terms of suitable combinations of $ \log_p^{n-m}(z) \Li_{m}(z) $.

By Proposition~6.4 and the correct version of Proposition~6.1 of \cite{Col82} the $ \Li_{n,p}(z) $
for $ n \ge 0 $ and $ z $ in $ \Cp\setminus\{1\} $ satisfy
\begin{equation*}
\Li_{n,p}(z) + (-1)^n \Li_{n,p}(1/z) = -\frac{1}{n!} \log_p^n(z)
\end{equation*}
when $ z \ne 0 $, and
\begin{equation}\label{distrieqp}
\Li_{n,p}(z^m) = m^{n-1} \sum_{\z^m=1} \Li_{n,p}(\z z)
\end{equation}
when $ m \ge 1 $ and $ z^m \ne 1 $.

For $ z $ in $ \Cp\setminus\{0,1\} $ and fixed $ n \ge 2 $, we define
\begin{equation}\label{Pnpdef}
\Lmod n (z) = \sum_{j=0}^{n-1} c_j \log_p^j(z) \Li_{n-j}(z)
\qquad (n \ge 2)
\end{equation}
for any $ c_j $ in $ \Cp $ satisfying $ c_0=1 $ and $ \sum_{j=0}^{n-1} c_j/(n-j)! =0 $, so 
that $ \Lmod n (z) + (-1)^n \Lmod n (1/z) = 0 $.
This is the case for $ c_j=b_j $
with $ b_j $ the $ j $-th Bernoulli numbers as before, although
the resulting formula is different from~\eqref{Pndef}.
Another possible natural candidate for the function $\Lmod n (z)$ is $
L_n(z) + L_{n-1}(z) \log_p(z)/n$, where
\begin{equation*}
  L_n(z) = \sum_{m=0}^{n-1} \frac{(-1)^m}{m!} \log_p^m(z) \Li_{n-m}(z) 
\end{equation*}
(see \cite[(1.3) and Remark~1.5]{BdJ03}).

The relations corresponding to~\eqref{poleq} and ~\eqref{distrieqp} for $ \Lmod n (z) $ ($ n\ge2 $)
are then
\begin{equation*}
\Lmod n (z) + (-1)^n \Lmod n (1/z) = 0
\end{equation*}
when $ z\ne 0,1 $, and
\begin{equation}\label{distrieqpLmod}
\Lmod n (z^m) = m^{n-1} \sum_{\z^m=1} \Lmod n (\z z)
\end{equation}
when $ m \ge 1 $ and $ z^m \ne 0,1 $.

For any branch of the \padic logarithm one then has the following result
for the \padic regulator map (see \cite[Theorem~1.10 or page 909]{BdJ03}).
Since its formulation depends on the embedding of $ k $ we deal with one such embedding at a time.

\begin{theorem}\label{synreg}
Let $ k $ be a number field, and let $ F \subset \Qpbar $ be the topological closure of the Galois
closure of any embedding $ k \to \Qpbar $.
For $ \s : k \to F $ let
$$
B_n^\s(k) = \langle  [x]_n\,\, |\text{ $ \s(x) $, $ 1-\s(x) $ are in $ \OO_\Cp^* $}\rangle
\subseteq B_n(k) = \dfrac{\tB_n(k)}{\Ker(\dd_n) \hskip-.5pt \cap \Ker(\widetilde P_n)}
\,.
$$
Then
\begin{alignat*}{1}
\tB_n(k) & \to F \cr
[x]_n & \mapsto \Lmod n (\s(x))
\end{alignat*}
induces a map
$$
\Lmod n ^\s : B_n^\s(k) \to F
$$
and the solid arrows in
$$
\xymatrix{
B_n^\s(k) \cap {\Ker(\dd_n)} \,\, \ar@{^{(}->}[r]\ar[rrd]_-{(n-1)!\Lmod n ^\s}
&
\dfrac{\Ker(\dd_n)}{\Ker(\dd_n) \hskip-.5pt \cap \Ker(\widetilde P_n)} \ar[r]^-{\Psi_n}\ar@{.>}[rd]^-{\,\,(n-1)!\Lmod n ^\s}
&
\Kn k \ar[d]^-{\reg_p\circ\s_*}
\\
&
&
\,F
}
$$
form a commutative diagram.
\end{theorem}

\begin{remark}\label{synregremark}
It was conjectured in \cite[Conjecture~1.14]{BdJ03} that the dotted arrow exists and that the full diagram commutes.
This is known to  hold for $ N [\z]_n $ if $ \z $ is any $ N $-th root of unity other than 1 (see Theorem~1.12 of loc.\ cit.).
\end{remark}

\begin{remark}
The formulae for the complex and \padic regulators involved in Theorems~\ref{zagierthm} and~\ref{synreg} are stated up to sign
in \cite[Proposition~4.1]{dJ95} (with a normalizing factor) and \cite[Theorem~1.10(2) or Proposition~7.14]{BdJ03} respectively, since they
depend on the choice of $ \Psi_n $, which is natural only up to sign.
In the \padic case it is clear from \cite[Proposition~7.10]{BdJ03} that the sign is $ (-1)^n $ if the
\lq\lq relativity isomorphism\rq\rq\ in the \padic case is normalized as in Proposition~5.7 of loc.\
cit.\ and the $ K $-theoretic relativity isomorphism as in (3.1) of loc.\ cit.\ 
corresponds to this under the naturality of the regulator map.
In order to see that the sign for the complex regulator is the same with the same choice of $ \Psi_n $,
we observe that one can use the techniques of \cite[Appendix~A and Section~5]{BdJ03} to describe the target
for the complex regulator of $ [z]_n $ with $ z $ in $ \C\setminus\{0,1\} $ in \cite{dJ95}
by means of a complex with the same formal structure as in \cite[Definition~5.2]{BdJ03}, but using
$ \R(n-1) $-valued $ C^\infty $-forms $ \psi $ with logarithmic poles on copies of $ (\mathbb{P}_\C^1\!\setminus\!\{1,z\})^j $
indexed by increasing functions $ f : [1,\dots,j] \to [1,\dots,n-1] $.  Using the maps
$ \psi \mapsto \int_{ (\mathbb{P}_\C^1)^j} \dd \arg(t_{f(1)}) \wedge \cdots \wedge \dd \arg(t_{f(j)}) \wedge \psi $
for a $ j $-form $ \psi $ with index $ f : [1,\dots,j] \to [1,\dots,n-1] $, and suitable normalizing factors, one writes down an $ \R $-linear map $ \Pi_\infty $ on such forms, with values in
$ \R(n-1) $, that vanishes on exact forms and satisfies the analogous normalization condition as $ \Pi=\Pi_p $
in \cite[Proposition~5.7]{BdJ03}.  Moreover, as in \cite[Section~2.5]{dJ95} one sees that the subcomplex
of forms on the component where $ j=n-1 $ that vanish if some $ t_l = 0 $ or $ \infty $ gives a quasi-isomorphism
under inclusion, and that the image of the regulator of $ [z]_n $ is given by $ \varepsilon_n(z) $ 
as in Section~4 of loc.\ cit.\ (up to sign since one has to use
the correct formula for $ \omega_n $, which can be obtained as in \cite[Section~6]{BdJ03}).
Using the computations on pages 236 and 237] of \cite{dJ95} one then sees that $ \Pi_\infty $ maps
$ \varepsilon_n(z) $ to $ (-1)^n (n-1)! P_n^\flat (z) $, where
$ P_n^\flat(z) =  \pi_{n-1}\left( \sum_{l=0}^{n-1} \frac{1}{l!} (-\log|z|)^l \Li_{n-l}(z) \right) $,
which induces the same map on $ \Ker(\dd_n)/\Ker(\dd_n) \cap \Ker(\widetilde P_n) $  as $ (-1)^n (n-1)! \, P_n(z) $.
Since the relativity isomorphism in $ K $-theory only depends on a choice of $ t_l =0 $ or $ \infty $ for
$ l=1,\dots,n-1 $ this shows that the signs in the formulae for the complex and \padic regulators of
$ [z]_n $ match for the same choice of $ \Psi_n $, and we adjust $ \Psi_n $ so that those signs disappear.
\end{remark}

We conclude this section by considering Conjecture~\ref{motiveconjecture} when the action of $ G $
on $ E[G]\pi $ is Abelian.

\begin{prop}\label{abelprop}
Let $ N \ge 2 $, $ k = \Q(\mu_N) $, $ G = \Gal(k/\Q) = (\Z/N\Z)^* $ where $ a $ in $ (\Z/N\Z)^* $
corresponds to $ \psi_a $ in $ G $ satisfying $ \psi_a(\z) = \z^a $ for all roots of unity $ \z $ in $ k $.
Assume that  $ E $ contains a root of unity of order equal to the exponent of $ G $
and let $ \pi $ in $ E[G] $ be the idempotent corresponding to an irreducible character $ \chi $ of $ G $.
Then parts (1), (2) and~(3) of Conjecture~\ref{motiveconjecture} hold
for $ \pi $ and $ n \ge 2 $ if $ \chi(-1) = (-1)^{n-1} $.
\end{prop}

\begin{proof}
That Conjecture~\ref{motiveconjecture}(1) holds in this case is well-known (cf.\ \cite{Neu88})
but we make it explicit here for part~(3) of the conjecture
(cf. \cite[Proposition~3.1]{Gro94}).

We deal with the case $ N=2 $ later so by Remark~\ref{sameconj}(1) we may assume $ N > 2 $
and $ \chi $ primitive. 
The idempotent $ \pi $ corresponding to $ \chi $ is
$ |G|^{-1} \sum_{a=1}^{\prime N} \chi^{-1}(a) \tensor \psi_a $, where the prime
indicates that we sum over $ a $ that are relatively prime to $ N $.

Note that by our assumption on $ \chi(-1) $ and Proposition~\ref{kerisreal} the dimension of $ K_{2n-1}(M_\pi^E) $ equals~1.
Fixing a generator $ \z $ of $ \mu_N \subset k^* $ we
use $ \pi $, $ \pi(1\tensor\z) $ and $ \pi(1\tensor [\z]_n) $
as basis vectors of $ E[G]\pi $, $ M_\pi^E $ and $ K_{2n-1}(M_\pi^E) $ respectively, where we simplify
notation by writing $ [\z]_n $ instead of $ \Psi_n([\z]_n) $ as in Theorem~\ref{zagierthm}.

In order to verify~(1) in $ E\tensor_\Q\C $ it suffices
by Remark~\ref{newremark}(2) to consider all embeddings $ \s: E \to \C $.
If $ \eta=\phi_\infty(\z) $ with $ \phi_\infty $ as in Section~\ref{motivicsection}
then the $ \s $-component  of $ D(M_\pi^E)^{1/2,\infty} $ becomes
\begin{equation*}
D(M_\pi^E)^{1/2,\infty}_\s
=
|G|^{-1} \sum_{a=1}^N \raise7pt\hbox{$\prime$} \chi_\s^{-1}(a) \eta^a
=
N |G|^{-1} \left(\sum_{a=1}^N \raise7pt\hbox{$\prime$} \chi_\s(a) \eta^{-a} \right) ^{-1}
\end{equation*}
by \cite[(3) on p.84]{Lan70}.
Similarly, by Theorem~\ref{zagierthm} the $ \s $-component of $ R_{n,\infty}(M_\pi^E) $ is
\begin{equation*}
R_{n,\infty}(M_\pi^E)_\s
=
(n-1)!\, |G|^{-1} \sum_{a=1}^N  \raise7pt\hbox{$\prime$} \chi_\s^{-1}(a) \Pn_n(\eta^a)
=
(n-1)! \, |G|^{-1} \sum_{a=1}^N  \raise7pt\hbox{$\prime$} \chi_\s^{-1}(a) \Li_n(\eta^a)
\end{equation*}
because $ \Pn_n(\eta^a) = (\Li_n(\eta^a)+\chi_\s^{-1}(-1)\Li_n(\eta^{-a}))/2 $ by our assumption on $ \chi(-1) $,
where $ \Li_n $ is computed using the power series in~\eqref{Lindef}.
They relate to the $ \s $-component of $ \L n {\chi} {\Q} $ via
$ \L n {\chi_\s} {\Q} = (n-1)!^{-1} R_{n,\infty}(M_\pi^E)_\s / D(M_\pi^E)^{1/2,\infty}_\s $
according to (2) on page 172 of \cite{Col82}
with the correct sign in the exponent of the Gauss sum as used here
(cf.\ \cite[p.~421]{Zag91}).  Since $ \L n {\chi_\s} {\Q} \ne 0 $ this also shows that
$ \pi(1\tensor [\z]_n) $ gives a basis of $ K_{2n-1}(M_\pi^E) $
as claimed.

Similarly, for~(2) and~(3) we consider all embeddings $ \t : E \to  \Qpbar  $.  If $ \phi_p(\z)=\eta_p $ then the $ \t $-component
of $ D(M_\pi^E)^{1/2,p} $ is $ N |G|^{-1} ( \sum_{a=1}^{\prime N} \chi_\t(a) \eta_p^{-a} )^{-1} $ and the $ \t $-component
of $ R_{n,p}(M_\pi^E) $ equals $ (n-1)! \, |G|^{-1} \sum_{a=1}^{\prime N} \chi_\t^{-1}(a) \Lmod n (\eta_p^a) $
by Remark~\ref{synregremark}, independent of our assumption on $ \chi(-1) $.   By~(3) on page~172 of
\cite{Col82} we have
\[
\Lp p n {\chi_\t\omega_p^{1-n}} {\Q} =  (n-1)!^{-1} \Eul_p(k,\chi_\t,\Q) R_{n,p}(M_\pi^E)_\t / D(M_\pi^E)_\t^{1/2,p}
\,.
\]
Therefore Conjecture~\ref{motiveconjecture}(1)-(3) hold with $ e(n,M_\pi^E) = e_p(n,M_\pi^E) $
in $ \Q^* \subseteq E^* $.

Now assume that $ N=2 $ so that $ \chi $ is the trivial character 1, $ \pi = 1 $ and $ n \ge2 $ is odd.
Taking $ 1 $, $ 1\tensor 1 $ and $ 1 \tensor [-1]_n $ as bases of $ E[G] = E $, $ M_\pi^E = E \tensor \Q $ and
$ K_{2n-1}(\Q)_E $ we find from Theorem~\ref{zagierthm},~\eqref{Lindef},~\eqref{Pndef} and~\eqref{distrieq}
that
\begin{equation*}
\L n 1 {\Q} = \z(n) = \frac{2^{n-1}}{1-2^{n-1}} P_n(-1) = \frac{2^{n-1}}{(n-1)!(1-2^{n-1})} R_{n,\infty}(M_1^E)
\,,
\end{equation*}
again proving that $ [-1]_n $ gives a basis for $ \Kn {\Q} $.
Coleman proves that
\[
\Lp p n \omega_p^{1-n} {\Q} = (1-p^{-n}) \lim_{x \to 1}{}'\, \Li_{n,p}(x)
\]
where the limit is taken in any subfield of $ \Cp $  that is of finite ramification degree over $ \Qp $.
But by~\eqref{Pnpdef} and~\eqref{distrieqpLmod} we can rewrite this as
\begin{equation*}
\Lp p n \omega_p^{1-n} {\Q}
= \frac{ (1-p^{-n}) 2^{n-1}}{1-2^{n-1}} \Lmod n (-1)
= \frac{(1-p^{-n}) 2^{n-1}}{(n-1)! (1-2^{n-1})} R_{n,\infty}^p(M_1)
\,.
\end{equation*}
\end{proof}

\begin{remark}\label{abelrem}
The $ \pi $ in the proof of Proposition~\ref{abelprop} are primitive in the sense of Definition~\ref{mindef}.
So, by Remarks~\ref{consequencesremark}(4) and~\ref{sameconj}(1),
if Conjecture~\ref{motiveconjecture} applies to $ \pi $ in $ E[G] $, then 
parts~(1), (2) and~(3) of it hold
for $ \pi $ for any $ E $ if the action on $ E[G]\pi $ is Abelian.
In particular, parts (1) and~(2) of Conjecture~\ref{pborel}
hold for any totally real Abelian number field 
by Remark~\ref{consequencesremark}.
\end{remark}

\begin{remark}\label{dirichletcalculations}
Computer calculations show that $ \Lp p n {\omega_p^{1-n}} {\Q} $ is in $ \Q_p^* $
when $ p=2,\dots,19 $ and $ n=2,\dots,20 $ is odd, verifying Conjecture~\ref{pborel}(3) in those cases.
Similarly, with notation and assumptions as in Proposition~\ref{abelprop},
for the 470 primitive characters $ \chi $ of $ \Gal(\Q(\m_N)/\Q) = (\Z/N\Z)^* $ with $ 2\le N \le 50 $, 
$ \Lp p n {\chi \tensor \omega_p^{1-n}} {\Q} $ lies in $ (E\tensor \Qp)^* $
for those values of $ p $ and $ n=2,\dots,20 $ whenever  $ \chi(-1) = (-1)^{n-1} $.
Thus Conjecture~\ref{motiveconjecture}(4) also holds in those cases.
\end{remark}

\section{Computing $ K $-groups in practice}
\label{practical}

Let $ n \ge 2 $.
It follows from Theorem~\ref{zagierthm} that Zagier's approach as described in Section~\ref{ksection} can be used to obtain at least
some part of $ \Kn k $.
 Using the notation introduced in Section~\ref{ksection}, in order to carry this approach out in practice
one starts with two sets $ S \subseteq S' $ of primes of $ \O $ and considers only 
$ [x]_n^\sim $ where $ x $ is an $ S $-unit and $ 1-x $ is an $ S' $-unit.
When we consider several
hundred of such $ [x]_n^\sim $'s (and avoid using both $ [x]_n^\sim $ and $ [1/x]_n^\sim $ since
$ [x]_n + (-1)^n [1/x]_n = 0 $ as mentioned in Section~\ref{ksection}) this method is
well suited to computer calculations.
The only point in Zagier's approach that cannot be carried out algebraically
is to determine which elements in $ \Ker(\tilde\dd_n) $ are actually
in $ \Ker(\tilde\dd_n) \cap \Ker(\widetilde P_n) $.  This is done using standard
methods for finding linear relations between the $ \widetilde P_n(\a_j) $'s for a basis $ \{\a_j\} $ of $ \Ker(\tilde\dd_n) $
with small integral coefficients.  
Typically, we used
about the first 50 decimals after the decimal point of all the polylogarithms
of all elements of $ k $ involved in the $ \a_j $ embedded into $ \C $ in all possible ways.
The relations that were found this way were then verified to hold up to at least 30 additional
decimal places.
As a final check, when $ \Ker(\dd_n) $ had the same rank as $ K_{2n-1}(k) $ the number $ q $
in~\eqref{borelq} (but now using the regulator $ \widetilde V_n(k) $ for the subgroup of $ \Kn k $ we obtained rather
than $ V_n(k) $; cf.~Remark~\ref{basisremark}) was computed with a working precision of up to 120
decimals.
It turned out that its reciprocal looked rational rather convincingly.

For example, for the cubic field $ k = \Q[x]/(x^3 - x^2 - 3 x + 1) $ with discriminant 148, one of the fields used in Example~\ref{S3computation},
we find for $ \z_k(n) |D_k|^{1/2} / \widetilde V_n(k) $ the values
\begin{alignat*}{1}
& -1.7531044558071585098612125639152666179693206722 \times 10^{-2}      \qquad (n=3)
\cr
&  -4.1170685884062518549452525064367732455835167754 \times 10^{-9}     \qquad (n=5)
\end{alignat*}
with reciprocals
\begin{alignat*}{1}
& -57.041666666666666666666666666666666666666666662            \qquad (n=3)
\cr
& -242891265.59999999999999999999999999999999999998            \qquad (n=5)
\,.
\end{alignat*}

As the calculation of a basis of $ \Kn k $, using the method outlined above,
tends to take (substantially) longer when $ [k:\Q] $, the absolute value of the discriminant of $ k $, or 
$ n $ increase, we describe two straightforward methods to reduce this.   They, and most of our arguments
below, will rely on Remark~\ref{galoisdescent}.

\begin{method}\label{firstmethod}
If $ k $ is an Abelian extension of $ \Q $, hence is contained in a cyclotomic field $ k' $,
then we can find $ \Kn k $ by computing $ \Kn k' $ as in Example~\ref{cyclo} and applying Remark~\ref{galoisdescent}
$($using Remark~\ref{equivariance}$)$.
\end{method}

Since parts (1), (2) and (3) of Conjecture~\ref{motiveconjecture} are known for all idempotents
of $ E[G] $ to which it applies if $ G $ is Abelian by Remark~\ref{consequencesremark}
and Proposition~\ref{abelprop}, this will not be used for
verifying Conjecture~\ref{motiveconjecture}.  Instead, in the cases that we shall consider 
we rely on the following method in order to find the $ K $-groups.

\begin{method}\label{subfields}
Let $ k/\Q $ be Galois with Galois group $ G $, $ E/\Q $ any extension, and
$ M $ an irreducible $ E[G] $-module in $ \KnE k $.
If $ H $ is a subgroup of $ G $ with $ M^H \ne 0 $ then 
$ M \subseteq E[G] \cdot \Kn k^H $.
If for every irreducible $ M $ we can take $ H \ne \{e\} $, then we reduce to
finding $ \Kn l $ together with the action of $ G $ on it inside $ \Kn k $ for proper subfields $ l $ of $ k $.
\end{method}

We now discuss in Examples~\ref{S3example},~\ref{D8example} and~\ref{sexticexample} how
to use Method~\ref{subfields} for certain Galois extensions $ k/\Q $.
In the decomposition of $ \KnE k $ according to central idempotents $ \pi_i $ of $ E[G] $
as in Section~\ref{motivicsection} we concentrate on
those primitive central $ \pi_i $ for which the action on $ E[G]\pi_i $ is not Abelian and
describe corresponding primitive idempotents $ \pi $ (as in Definition~\ref{mindef}(2)) for later
use.
Due to restrictions on when we can calculate the \padic \Lfunctions numerically~(see Section~\ref{pLsection})
we only consider $ k $ that are Abelian over a quadratic subfield.
Of course Proposition~\ref{whenconj} describes when Conjecture~\ref{motiveconjecture}
applies to $ \pi $, but we work out the structure of the $ K $-groups and the Galois action in
more detail.

\begin{example}\label{S3example}
Let $ k/\Q $ be an $ S_3 $-extension with quadratic subextension $ q $.
With $ E=\Q $ and $ \s $ in $ G $ of order three,
$ \pi_1 = (e+\s+\s^2)/3 $ and $ \pi_2 = 1-\pi_1 $ are orthogonal central idempotents with $ \pi_2 $ primitive
(as in Definition~\ref{mindef}(2)).  In $ \Q[G] = \Q[G] \pi_1  \oplus  \Q[G] \pi_2 $
the first summand consists of Abelian representations of $ G $ and the last of two copies of
the irreducible 2-dimensional representation $ V $ of~$ G $. 
A corresponding primitive idempotent for $ \pi_2 $ is $ \pi = \pi_2 \pi_H $ for
any subgroup $ H = \langle \t \rangle $ of $ G $ of order~2 where $ \pi_H = (e+\t)/2 $.

By Remark~\ref{galoisdescent} we have $ \pi_1 \Kn k = \Kn q $ so that
\begin{equation}
\label{S3iso}
\Kn k = \Kn q \oplus \pi_2 \Kn k
\,.
\end{equation}
The last summand is a $ \Q[G] $-module isomorphic to $ V^{t_n} $ for some $ t_n \geq 0 $.
As $ V = \pi V \oplus \t'\pi V $ for any $ \t'\ne\t $ of order~2, we find, if $ c = k^H  $,
\[
\pi_2 \Kn k = \pi_2 \Kn c \oplus \t' \pi_2 \Kn c
\]
again by Remark~\ref{galoisdescent}.
So we reduce to the calculation of $ \Kn q $ and $ \Kn c $ together with the action of~$ G $ on the latter.

A dimension count in~\eqref{S3iso} using Theorem~\ref{boreltheorem} determines~$ t_n $.
We distinguish two cases.

{\it Case 1.\/}
If $ q $ is real then $ t_n = 0 $ for $ n $ even and 2 for $ n $  odd.

{\it Case 2.\/}
If $ q $ is imaginary then $ t_n = 1 $.

\noindent
Conjecture~\ref{motiveconjecture} applies to $ \pi $ when $ t_n=2 $.
\end{example}

\begin{example}\label{D8example}
Let $ k/\Q $ be a $ D_8 $-extension where $ D_8 = \langle \s,\t | \s^4=\t^2=(\s\t)^2 = e \rangle $.
We fix an isomorphism $ G = \Gal(k/\Q) \iso D_8 $ and use notation as in Figure~\ref{D8figure}.

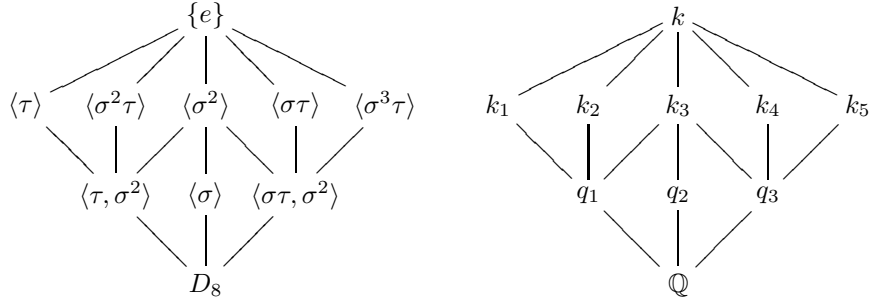
\begin{figure}[htb]
\hbox{
\quad
\xymatrix
@!=10pt
{
&&\{e\}\ar@{-}[lld]\ar@{-}[ld]\ar@{-}[d]\ar@{-}[rd]\ar@{-}[rrd]&&\cr
\langle\t\rangle& \langle\s^2\t\rangle & \langle\s^2\rangle & \langle\s\t\rangle & \langle\s^3\t\rangle \cr
& \langle\t,\s^2\rangle \ar@{-}[lu] \ar@{-}[u] \ar@{-}[ru] & \langle\s\rangle \ar@{-}[u] & \langle\s\t,\s^2\rangle \ar@{-}[lu] \ar@{-}[u] \ar@{-}[ru] & \cr
&& D_8 \ar@{-}[lu]\ar@{-}[u]\ar@{-}[ru] && \cr
}
\quad\hfill\quad
\xymatrix
@!=10pt
{
&& k \ar@{-}[lld]\ar@{-}[ld]\ar@{-}[d]\ar@{-}[rd]\ar@{-}[rrd]&&\cr
k_1& k_2 & k_3 & k_4 & k_5 \cr
& q_1 \ar@{-}[lu] \ar@{-}[u] \ar@{-}[ru] & q_2 \ar@{-}[u] & q_3 \ar@{-}[lu] \ar@{-}[u] \ar@{-}[ru] & \cr
&& \Q \ar@{-}[lu]\ar@{-}[u]\ar@{-}[ru] && \cr
}
\quad
}
\caption{The subgroups of $ G =  D_8 $ and their fixed fields in $ k/\Q $.}
\label{D8figure}
\end{figure}

\noindent 
Then $ \Q[G] = \Q[G] \pi_1 \oplus  \Q[G] \pi_2 $ where  $ \pi_1 = (1+\s^2)/2 $ and $ \pi_2 = 1-\pi_1 $ are
orthogonal central idempotents.  The first term is an Abelian representation of $ G $ and
the second as a $ \Q[G] $-module is isomorphic to $ V^2 $ for $ V $ the irreducible 2-dimensional representation
of $ G $ over~$ \Q $.
Then $ \pi_2 $ is a primitive central idempotent and  $ \pi = \pi_2 \pi_H $ a corresponding primitive idempotent 
if we let $ H = \langle\t \rangle $ and $ \pi_H = (1+\t)/2 $.

Since $ \pi_1 \Kn k = \Kn k_3 $ by Remark~\ref{galoisdescent}, for $ \Kn k $ we find
\begin{equation}\label{D8iso}
\Kn k =  \Kn k_3 \oplus \pi_2 \Kn k
\end{equation}
with the second term a $ \Q[G] $-module isomorphic with $ V^{t_n} $ for some $ t_n \geq 0 $.
If $ \t' $ in $ G $ has order two and does not commute with $ \t $ then $ V = \pi V \oplus \t'\pi V $,
so by Remark~\ref{galoisdescent} we obtain 
\[
\Kn k =  \Kn k_3 \oplus  \pi_2 \Kn k_1 \oplus \t' \pi_2 \Kn k_1 
\]
and we reduce to the calculation of $ \Kn k_3 $ and $ \Kn k_1 $ together with the action of $ G $ on
the latter.

We can find $ t_n $ by counting dimensions in~\eqref{D8iso}, and distinguish three cases.

{\it Case 1.\/}
For $ k_3 $ and $ k $ totally real $ t_n =0 $ for $ n $ even and $ t_n=2 $ for $ n $ odd.

{\it Case 2.\/}
For $ k_3 $ totally real but $ k $ not $ t_n=2 $ for $ n $ even and $ t_n=0 $ for $ n $ odd.

{\it Case 3.\/}
For $ k_3 $ not totally real $ t_n=1 $.

\noindent
Again Conjecture~\ref{motiveconjecture} applies to $ \pi $ when $ t_n = 2 $.
\end{example}

\begin{example}\label{sexticexample}
Let $ k/\Q $ be an $ S_3 \times \Z/3\Z $-extension, identify $ G = \Gal(k/\Q) $ and $ S_3 \times \Z/3\Z $
via a fixed isomorphism, and choose a generator $ \s $ of $ A_3 \subset S_3 $.
Let $ q $ be the quadratic subfield of $ k $, $ s' $ the fixed field of $ A_3 \times \{0\} $,
$ s'' $ the fixed field of $ \{e\} \times \Z/3\Z $ and $ s_i $ for $ i=1 $, 2 the fixed field of
$ H_i = \langle (\s,i) \rangle $.

If $ E $ is an extension of $ \Q $ that contains a primitive third root of unity $ \z_3 $ then
\begin{equation}\label{sexticEGdecomp}
\begin{aligned}
   E[G] & = E[G] \pi_1 \oplus E[G] \pi_2 \oplus E[G] \pi_3 \oplus E[G] \pi_4
\cr
        & \iso E[\Z/6\Z] \oplus W_1^2 \oplus W_{\z_3}^2 \oplus W_{\z_3^2}^2 
\end{aligned}
\end{equation}
where
\begin{equation*}
\begin{aligned}
\pi_1 & = ((e,0) + (\s,0) + (\s^2,0))/3 
\cr
\pi_2 & = (1-\pi_1)((e,0)+(e,1)+(e,2))/3
\cr
\pi_3 & = (1-\pi_1)((e,0)+\z_3^2(e,1)+\z_3(e,2))/3
\cr
\pi_4 & = (1-\pi_1)((e,0)+\z_3(e,1)+\z_3^2(e,2))/3 
\end{aligned}
\end{equation*}
are orthogonal central idempotents of $ E[G] $ satisfying $ \pi_1 + \pi_2 + \pi_3 + \pi_4 = (e,0) $, and
$ W_a  $ is the irreducible representation of $ G $ given by $ (\r,m)(v) = a^m \r(v) $ for
$ v $ in the irreducible 2-dimensional representation of $ S_3 $.
With $ \pi_{H_i} = \frac13 \sum_{h \in H_i} h $, a primitive idempotent corresponding to the
primitive central idempotent~$ \pi_j $ ($ j=3 $,~4) is $ \pi = \pi_j \pi_{H_i} $.

Using Remark~\ref{galoisdescent} for $ A_3 \times \{0\} $ and $ \{e\} \times \Z/3\Z $ we see that 
\[
\begin{aligned}
\KnE k & = \pi_1 \KnE k \oplus \pi_2 \KnE k \oplus \pi_3 \KnE k \oplus \pi_4 \KnE k 
\cr
      & = \KnE s' \oplus (1-\pi_1) \KnE s'' \oplus \pi_3 \KnE k \oplus \pi_4 \KnE k 
\,.
\end{aligned}
\]
Because $ W_a = W_a^{H_1} \oplus W_a^{H_2} $ when $ a \ne 1 $ we have
\[
\pi_j \KnE k = \pi_j \KnE s_1 \oplus \pi_j \KnE s_2 
\qquad (j=3, 4)
\]
so that we reduce to the calculation of $ \Kn s $ for the sextic subfields $ s $ of $ k $,
together with the action of $ G $ when $ s = s'' $, $ s_1 $ or~$ s_2 $.

We have $ \pi_3 \KnE k \iso W_{\z_3}^{t_n} $ and $ \pi_4 \KnE k \iso W_{\z_3^2}^{t_n} $ for the same $ t_n $
as one sees by taking $  E = \Q(\z_3) $ and considering the obvious action of $ \Gal(E/\Q) $ on $ E[G] $.
Since the dimension of $ (1-\pi_1) \KnE s'' $ was obtained in Example~\ref{S3example} we get the following
results for~$ t_n $.

{\it Case 1.\/}
For $ q $ real $ t_n = 0 $ when $ n $ is even and $ t_n=2 $ when $ n $ is odd.

{\it Case 2.\/}
For $ q $ imaginary $ t_n = 1 $.

\noindent
Conjecture~\ref{motiveconjecture} applies to $ \pi $ corresponding to $ \pi_j $ ($ j=3, $~4)
when $ t_n=2 $.
\end{example}

We finish this section with a result that was promised in Remark~\ref{newrem}(1).  Note that the assumption 
$ \dim_E(E[G]\pi) = \dim_E(\pi \KnE k ) $ is not needed here: we consider the regulator map with
values in $ \pi(\oplus_{\t:k\to\C} E\tensor \R(n-1))^+ $ obtained from~\eqref{Bregulator}, and use Beilinson's
regulator $ R_\Bei $ comparing $ E $-structures as in Remark~\ref{BBmotiverem}.  
Beilinson's conjecture \cite[Conjecture~3.4]{bei85} states that the order of vanishing at $ s=1-n $ of
$ \L s {\chi_\pi^\vee\tensor\id} {\Q} $ should be $ \dim_E(\pi \KnE k ) $ and that 
the first non-vanishing coefficient at $ s=1-n $ equals
\begin{equation*}
\Lsharp {1-n} {\chi_\pi^\vee\tensor\id} {\Q} = e R_\Bei
\end{equation*}
for some $ e $ in $ E^* $.

\begin{proposition}\label{grossprop}
Let $ k $ be a $ G=S_3 $-extension $($resp.\  $ G= D_8 $-extension$)$ of $ \Q $, and let $ n \ge 2 $.
Then this conjecture holds for any primitive idempotent $ \pi $ corresponding to $ \pi_2 $ in
Example~\ref{S3example} $($resp.\ Example~\ref{D8example}$)$.
\end{proposition}

\begin{proof}
Clearly it suffices to prove this when $ E=\Q $.
If $ G=S_3 $ then by Remark~\ref{consequencesremark}(3) we may assume
$ \pi_H = \pi_1 \pi_H + \pi $ with $ \pi_1 \pi_H $ and $ \pi = \pi_2 \pi_H $
orthogonal idempotents, where we use notation as in Example~\ref{S3example}.  The conjecture for $ \pi_H $
follows from Theorem~\ref{boreltheorem} for $ c = k^H $ by the arguments of Remark~\ref{sameconj} and
the functional equation~\eqref{oldfuneq} of $ \z_c(s) $.
Similarly, it holds for $ \pi_1 \pi_H = \frac16 \sum_{g \in G} g $,
where it corresponds to that theorem for $ k^G = \Q $,
and it follows that it must hold for $ \pi $ as well.
If $ G=D_8 $ the argument is identical, using
Theorem~\ref{boreltheorem} for $ k_1 $ and $ q_1 $.
\end{proof}

\section{Computing \padic \Lfunctions}
\label{pLsection}

We briefly sketch the method used to compute values of \padic
$L$-functions of 1-dimensional Artin characters.  Values of \padic
\Lfunctions of characters of higher dimension can then be deduced
using \eqref{artinLp}. The method we employ is somewhat similar to, but more technical than
the construction in \cite{Colm88}, and generalizes the one used in \cite{Sol-Rob}
to compute values of $p$-adic $L$-functions of real quadratic fields
at $s = 1$.   
Note that several modifications are necessary compared to \cite{Colm88} since the
construction there is for partial $ \z $-functions of the ray
class modulo $p$ whereas we need to work with partial $ \z $-functions of
arbitrary class groups. Also, we need to replace $p = 2$ in a certain
number of results by $q = 4$. We leave it to the careful reader to check
that these changes are indeed possible, or refer to the forthcoming
article \cite{Rob06} where the details of the actual method used will be
given.

Let $k$ be a totally real number field and let $d=[k:\Q]$. Let $\chi$ be a
$\Qpbar$-valued 1-dimensional Artin character of $\Gal(\kbar/k)$. As
in Section~\ref{conjecturesection}, we let $q = p$ if $p$ is odd, and
$q = 4$ if $p = 2$.

By class field theory, $\chi$ can be viewed as a character of a
suitable ray class group of~$k$.  We denote by $\mathfrak{f}$ the
conductor of $\chi$, which is an integral ideal of $k$, and we let
$\mathfrak{g}$ be the modulus with infinite part $\mathfrak{g}_\infty$
consisting of all infinite places of $k$, and finite part
$\mathfrak{g}_0$ equal to the least common multiple of $\mathfrak{f}$
and $q$.  Let $\mathfrak{a}_1, \dots, \mathfrak{a}_h$ be
representatives of the classes of the ray class group of $k$ modulo
$\mathfrak{g}$, and let $\zeta_{\mathfrak{a}_i}(s)$, $i = 1, \dots,h$,
be the corresponding partial $ \z $-functions.  According to
\cite[Chapter~VII, p.526]{Neu99} we have, for any isomorphism
$\sigma:\Qpbar \to \C$,
\begin{equation*}
\sum_{i=1}^h \s\circ\chi(\mathfrak{a}_i) \zeta_{\mathfrak{a}_i}(s) 
=
\Eul_p(s,\s\circ\chi,k) \L s {\s\circ\chi} k 
\,.
\end{equation*}
Since all $ \zeta_{\mathfrak{a}_i}(m) $  are rational for integers $ m \leq 0 $
by loc.\ cit.\ Chapter~VII, Corollary~9.9, we may identify them with
$ \sigma^{-1}(\zeta_{\mathfrak{a}_i}(m)) $  in $ \Qp $.
Thus for such $ m $ we have
\begin{equation}\label{partialzeta}
\sum_{i=1}^h  \chi(\mathfrak{a}_i) \zeta_{\mathfrak{a}_i}(m) 
= 
\Eul^\ast_p(m, \chi, k) L^\ast(m, \chi, k) 
\end{equation}
where $L^\ast(m, \chi, k)$ and $\Eul^\ast_p(m, \chi, k)$ are as
in \eqref{sigmaindependence} and \eqref{sigmaeuler}.

Now let $\beta\ne1 $ in $ \mathcal{O}_k$ be such that 
\begin{enumerate}
\item 
$\beta \equiv 1 \pmod{\mathfrak{g}_0}$ and $v(\beta) > 0$ for every infinite place $v$
of $k$;

\item
$\O/(\beta) \simeq \Z/b\Z$, where $b = \mathcal{N}(\beta)$ is the norm of the principal
  ideal $(\beta)$.
\end{enumerate}
Scaling the measures obtained in \cite[Lemme~4.4]{Colm88} for the $ \mathfrak{a}_i $ we obtain $ \Qp $-valued measures
$\widetilde\lambda_i$ ($ i=1,\dots,d $) on $\Zp^d$, depending on $\beta$ and $\mathfrak{a}_i$, such that
for all integers $m\le0$,
\[
(b^{m-1}-1) \zeta_{\mathfrak{a}_i}(m) = (-1)^{dm}
\mathcal{N}(\mathfrak{a}_i)^{-m} \int (x_1\cdots x_d)^{-m} \, d\widetilde\lambda_i
\,.
\]
The proof of Lemma~4.3 of loc.\ cit.\ shows, contrary to the statement
of that lemma, that the measures $\widetilde\lambda_i$ are supported
on $ (-1+q\Zp)^d $.  Therefore pulling back over multiplication
by $ -1 $ on $ \Zp^d $ we obtain measures $ \lambda_i $ on $ \Zp^d $,
supported on $(1+q\Zp)^d$, such that, for the same $ m $,
\begin{equation}\label{mnonpositive}
(b^{m-1}-1) \zeta_{\mathfrak{a}_i}(m) = 
\mathcal{N}(\mathfrak{a}_i)^{-m} \int (x_1\cdots x_d)^{-m} \, d\lambda_i
\,.
\end{equation}

For a fixed $x$ in $1 + q\Zp$ the function $s \mapsto x^s$ is defined
and analytic on $ \Zp $.  Since, in general, 
$\mathcal{N}(\mathfrak{a}_i)$ will not be congruent to $1$ modulo $q$,
we do the following.  For any integer $l$ and any
integral ideal $\mathfrak{a}$ of $k$ that is coprime to $p$, the 
function $s \mapsto \omega_p(\mathfrak{a})^l \langle
\mathcal{N}\mathfrak{a} \rangle^{-s}$ is an analytic $\Zp$-valued
function on $ \Zp $, whose value equals $ \mathcal{N}(\mathfrak{a})^{-m}$ at any
integer $m$ with $m + l \equiv 0 \pmod{\varphi(q)}$. 
We define the $p$-adic partial $ \zeta $-function of the class of $\mathfrak{a}_i$
by mapping $ s $ in $\Zp \setminus \{1\}$ to
\begin{equation*}
  \zeta_{p, \mathfrak{a}_i}(s) = 
  \left(b^{s-1} - 1\right)^{-1}
  \langle \mathcal{N}(\mathfrak{a}_i) \rangle^{-s} 
  \pint (x_1 \cdots x_d)^{-s} \, d\lambda_i \,,
\end{equation*}
where $ ' $ indicates that we restricted the domain of integration to $ (1+q\Zp)^d $.

Let $\Psi_i$ be the function defined by the above integral, that is, for $ s $ in $ \Zp $,
$$
\Psi_i(s) = \pint (x_1 \cdots x_d)^{-s} \, d\lambda_i\,.
$$
For $a$ in $ \Zp$ with $|a - 1|_p \leq q^{-1}$, and $ s $ in $ \Zp $, one can write $a^{-s} =
\exp_p(-s \log_p(a))$, where $\exp_p$ and $\log_p$ are the
\padic exponential and logarithm functions respectively; and thus the function $s
\mapsto a^{-s}$ can be expressed as a power series with coefficients
in $\Zp$, where the coefficient
of $s^m$ has absolute value at most $ q^{-m}p^{m/(p-1)}$. Developing $
(x_1\cdots x_d)^{-s} $ in this way as a power series of $s$ and using
the fact that the measures $\lambda_i$ have bounded norm 
by \cite[Lemme 4.2 bis]{Colm88},  we see that
$\Psi_i(s)$ can also be expressed as a power series in
$\Qp[[s]]$ where the absolute value of the coefficient of
$s^m$ is at most $ C_i q^{-m}p^{m/(p-1)}$ for some $ C_i>0 $.
In the same way, the
function $s \mapsto \langle \mathcal{N}(\mathfrak{a}_i) \rangle^{-s}$
can be expressed as a power series in $\Zp[[s]]$ whose
coefficients satisfy the same bounds.  
Similarly, since the only solution for $ s $ of $b^{s-1} = 1$ with $ |s| < q p^{-1/(p-1)} $
is $s = 1$, it follows from the \padic Weierstrass
preparation theorem that the function $s \mapsto b^{s-1}-1$ can be
expressed as the product of a power of $ p $, the polynomial $s-1$, and an invertible
power series in $\Zp[[s]]$ with the absolute value of the coefficient of $ s^m $ bounded by 
$ q^{-m} p^{m/(p-1)} $.
Therefore $\zeta_{p,\mathfrak{a}_i}(s)$ can be expressed as the quotient of a
power series in $\Qp[[s]]$ with the absolute value of the coefficient of $ s^m $ 
bounded by $ C_i' q^{-m} p^{m/(p-1)} $ with $ C_i'>0 $, and $s-1$.

If we define $p$-adic functions of $ s $,
\begin{equation}\label{merpL}
 L_p^{(l)}(s, \chi, k) =
\sum_{i=1}^h   \omega_p(\mathfrak{a}_i)^l \chi(\mathfrak{a}_i) \zeta_{p, \mathfrak{a}_i}(s) \,
\qquad (l \text{ modulo } \phi(q))
\,,
\end{equation}
then it follows from \eqref{partialzeta}, \eqref{mnonpositive}, and the equality
$ \zeta_{p,\mathfrak{a}_i}(m) = \omega_p(\mathfrak{a}_i)^m \zeta_{\mathfrak{a}_i}(m) $ for any integer $ m 
\le 0 $, that 
\begin{equation}\label{interpol1}
 L_p^{(l)}(m, \chi, k) = \Eul^\ast_p(m, \chi, k) L^\ast(m, \chi, k)
\end{equation}
for all non-positive $m$ such that $l + m \equiv 0 $ modulo $ \varphi(q) $.
In particular, the $p$-adic $L$-function of $\chi$ is given by $\Lp p s {\chi} k = L_p^{(-1)}(s, \chi, k)$.
Moreover, its residue at $s=1$ is zero if $\chi$ is non-trivial
because the functions $\omega_p(\mathfrak{a}_i)^{-1}\z_{p,\mathfrak{a}_i}(s)$
all have the same residue there by \cite[Corollaire on page 388]{Colm88}.

The estimates on the coefficients of the power series expansions above imply the
claims in Theorem~\ref{padicLfunction}(1).  Through the \padic Weierstrass preparation theorem
they also imply that the function in~\eqref{artinLp} can be written as the quotient of a power series 
in $ \Qp(\chi)[[s]] $ that converges
for $ s $ in $ \Cp $ with $ |s| < q p^{-1/(p-1)} $, and a polynomial in $ \Qp(\chi)[s] $.

We note in passing that if $l$ is any integer then
\begin{equation*}
  \Lp p s {\chi\omega_p^l} k 
=
 \sum_{i=1}^h \omega_p(\mathfrak{a}_i)^{l-1} \chi(\mathfrak{a}_i)  \zeta_{p, \mathfrak{a}_i}(s)
=
L_p^{(l-1)}(s, \chi, k)
\,, 
\end{equation*}
so that we obtain from \eqref{interpol1} that
\begin{equation*}
 \Lp p m {\chi\omega_p^l} k = \Eul^\ast_p(m, \chi, k) L^\ast(m, \chi, k) 
\end{equation*}
for all $m\le0$ with $l + m \equiv 1 $ modulo $ \varphi(q) $, as stated in Remark~\ref{interpolrem}.

In order to compute $L(s, \chi, k)$ for a fixed  $s$ in
$\Zp$, it suffices to compute $\zeta_{p, \mathfrak{a}_i}(s)$,
for which one mainly has to compute the $ \Psi_i(s) $.
The measures $\lambda_i$ can be represented as power
series in $ d $ variables,
$$
F_i(w_1, \dots, w_d) = \sum_{n_1, \dots, n_d \geq 0} a_{i, n_1, \dots, n_d} w_1^{n_1} \cdots w_d^{n_d} \,,
$$
with (bounded) coefficients
\begin{equation*}
 a_{i, n_1, \dots, n_d} = \int {x_1 \choose n_1} \cdots {x_d \choose n_d} d \lambda_i
\end{equation*}
that can be computed using Shintani's method~\cite{Shin76},
which is practical for calculations if the base field $k$ is $\Q$ or real
quadratic.  On the other hand, the function $ h $ on $\Zp$ defined by
$$
x \mapsto
\begin{cases}
  x^{-s} & \text{ if } x \text{ is in } 1+q\Zp, 
\\
  0 & \text{otherwise},
\end{cases} 
$$
is continuous and thus admits a Mahler expansion $\sum_{n \geq 0} c_{s, n} \binom{x}{n}$, where the
coefficients satisfy $c_{s, n} \to 0$ when $n \to \infty$ and can be easily computed recursively.  Then
we have 
\begin{equation*}
\Psi_i(s)
=
\int  h(x_1) \cdots h(x_d) \, d\lambda_i 
=
\sum_{n_1, \dots, n_d \geq 0} a_{i, n_1, \dots, n_d} c_{s,
  n_1} \cdots c_{s, n_d} \,. 
\end{equation*}
Thus, by computing enough terms in the above sum, we can get a good
approximation of the integral, and therefore, putting 
everything together, of the $p$-adic $L$-function.
 
\smallskip

In order to apply this we make~\eqref{brauerinduction} explicit in the cases we are interested in,
writing \Lfunctions of representations for the \Lfunctions of the corresponding characters.
In particular, we write the classical \Lfunctions of the irreducible 2-dimensional representations
of $ S_3 $ and $ D_8 $ in Examples~\ref{S3example} and~\ref{D8example}, and the irreducible 2-dimensional
representations $ W_a $ with $ a \ne 1 $ of $ S_3 \times \Z/3\Z$ in Example~\ref{sexticexample}, as
the \Lfunctions of Abelian characters over quadratic fields. 
For later use we also note how to write some of them in terms of $ \z $-functions of certain fields.

We take the field of coefficients $ E $ to be $ \C $ in the remainder of this section, and 
denote the trivial 1-dimensional representation of any group $ G $ by $ 1_G $.

\begin{example}\label{S3L}
Let  $ \Alt_{S_3} $ be the 1-dimensional representation of $ S_3 $ through the sign character, and $ V $ 
the irreducible 2-dimensional representation.  
Since 
\[
 \Ind A_3 S_3 \circ \Ind \{e\} A_3 (1_{\{e\}})  \iso \Ind \{e\} S_3 (1_{\{e\}}) \iso 1_{S_3} \oplus \Alt_{S_3} \oplus V^2   
\]
and $ \Ind A_3 S_3 (1_{A_3}) \iso 1_{S_3} \oplus \Alt_{S_3} $ it follows that
$ \Ind A_3 S_3 (V') \iso V $ for either non-trivial 1-dimensional representation $ V' $ of~$ A_3 $.
Applying this to the situation of Example~\ref{S3example} we see that
$ \L s  V {\Q} = \L s V' q $.
Similarly, if $ H $ is any subgroup of $ S_3 $ of order~2, then 
$ \Ind H S_3 (1_H) \iso 1_{S_3} \oplus V $.
In the situation of Example~\ref{S3example} this means that $ \z_{c}(s) = \z_{\Q}(s) \L s V {\Q} $.
\end{example}

\begin{example}\label{D8L}
If $ D_8 $ is as in Example~\ref{D8example}, $ \V a b $ for $ a,b=\pm1 $
the 1-dimensional representation of $ D_8 $ where $ \s $ acts as multiplication by
$ a $ and $ \t $ by $ b $, and $ V $ the irreducible 2-dimensional representation, then for any subgroup $ H $
$$
\Ind H D_8 \circ \Ind {\{e\}} H (1_{\{e\}}) \iso (\oplus_{a,b=\pm1} \V a b ) \oplus V^2
\,.
$$
If we take $ H = \langle \t, \s^2 \rangle $ then
$ \Ind {\{e\}} H (1_{\{e\}}) \iso \oplus_{a,b=\pm1} \W a b $
with $ \W a b $ the 1-dimensional representation of $ H $ where $ \s^2 $ acts as multiplication by 
$ a $ and $ \t $ by~$ b $.
Because $ \Ind H D_8 (\W 1 b ) \iso \oplus_{a=\pm1} \V a b $ we have $ \Ind H D_8 (\W -1 1 ) \iso \Ind H D_8 (\W -1 -1 ) \iso V $,
so with notation as in Figure~\ref{D8figure} we get $ \L s V {\Q}  = \L s {\W -1 1 } q_1 = \L s {\W -1 -1 } q_1 $.

Similarly, $ \Ind {\langle\t\rangle} D_8 (1_{\langle\t\rangle}) \iso ( \oplus_{a=\pm1} V_{a,1} ) \oplus V $
so that $ \z_{k_1}(s) = \z_{q_1}(s) \L s V {\Q} $ as well.
\end{example}

\begin{remark}
If $ H=\langle\s\rangle $ in Example~\ref{D8L}
then $ \Ind \{e\} H (1_{\{e\}}) \iso \oplus_{m=1}^4 U_{i^m} $ with
$ U_a $ the 1-dimensional representation of $ \langle\s\rangle $ where $ \s $ acts as multiplication
by~$ a $.  Now $ \Ind H D_8 (U_a) \iso \V a 1 \oplus \V a -1  $ when $ a = \pm 1 $
so that $ \Ind H D_8 (U_a) \iso V $ when $ a \ne \pm 1 $, and 
in Example~\ref{D8example} we also have $ \L s V {\Q}  = \L s U_i q_2 = \L s U_{-i} q_2 $.
\end{remark}

\begin{example}\label{sexticL}
If $ G = S_3 \times \Z/3\Z $ and $ H = A_3 \times \Z/3\Z $ then
$$
\Ind H G \circ \Ind {\{(e,0)\}} H (1_{\{(e,0)\}})
\iso W \oplus W_1^2 \oplus W_{\z_3}^2 \oplus W_{\z_3^2}^2 
$$
with $ W $ the direct sum of the 1-dimensional representations of $ G $
and the $ W_a $ as in Example~\ref{sexticexample}.
Write $  \Ind {\{(e,0)\}} H (1_{\{(e,0)\}}) \iso \oplus_{a,b} V_{a,b} $ where,
for $ a $ and $ b $ cubic roots of unity, $ V_{a,b} $ is the 1-dimensional representation of $ H $
on which $ (\s^m,n) $ acts as multiplication by $ a^m b^n $.
Since $ \Ind H G (\oplus_b V_{1,b}) \iso W $ we find by considering the action of $ (e,1) $
that $ \Ind H G (V_{a,b}) \iso W_b $  when $ a \ne 1 $.  In particular, in the notation of Example~\ref{sexticexample},
we have $ \L s W_b {\Q} = \L s V_{a,b} q $ since $ q = k^H $.

It is easy to see that the coefficients of $ p^{-s} $ in the Euler factors of
$ \L s W_{\z_3} {\Q} $ and $ \L s W_{\z_3^2} {\Q} $ are conjugate, but since not all of them can be real
it is not possible to express those functions in terms of $ \z $-functions of subfields of~$ k $.
\end{example}

\section{Examples}\label{examplesection}

In this section we describe the evidence for Conjecture~\ref{motiveconjecture} in the case of 
primitive idempotents $ \pi $ in $ E[\Gal(k/\Q)] $ as in Examples~\ref{S3example}, \ref{D8example} and \ref{sexticexample},
with $ E=\Q $ in the first two and $ E=\Q(\z_3) $ in the last.
In the first two cases part~(1) of the conjecture is known by Proposition~\ref{grossprop},
but in the last case our evidence for this part is numerical; the same holds for parts~(2) through (4) in all cases.

We would like to stress that the numerical verification of part~(4) of the conjecture actually proves
this part for all cases that we considered since we can check that a element of $ E\tensor\Qp $ is a unit
of that ring by computing a close enough approximation of it.

For the numerical calculations of the \padic regulator, $ R_{n,p}(M_\pi^E) $,
 we assumed that the syntomic regulator map on the subgroup of $ \Kn k $ described in Zagier's
Conjecture~\ref{zagierconjecture} 
is given by the $ \tLmod n ^\s $ (as in Remark~\ref{synregremark}).
Out of the many possible choices for $ \Lmod n (z) $ as described in ~\eqref{Pnpdef}
we used
\[
\Lmod n (z) = \Li_{n,p}(z) - \log_p^{n-1}(z) \Li_1(z) / n!
            = \Li_{n,p}(z) + \log_p^{n-1}(z) \log_p(1-z) / n!
\]
since it is relatively simple from a computational point of view.  The calculations of the $ \Lmod n (z) $
for the branch with $ \log_p(p)=0 $
were done in versions 2.11-7, 2.12-19 and 2.12-21 of \cite{magma},
using an implementation of the algorithm described in \cite{BdJ06}.

Note that the primitive idempotents in Examples~\ref{S3example} and \ref{D8example} are unique up to
conjugation in $ \Q[G] $ since they all correspond to $ \pi_2 $.  But in Example~\ref{sexticexample}
replacing $ \z_3 $ with $ \z_3^2 $ changes a primitive idempotent corresponding to $ \pi_3 $ into
one corresponding to $ \pi_4 $, and conversely.  Because the same holds for the 
classical and \padic \Lfunctions{} (with values in
$ E \tensor_\Q \C $ and $ E \tensor_\Q \Qp $ respectively), as well as the Euler factor for $ p $, the
conjectures for primitive idempotents corresponding to $ \pi_3 $ and $ \pi_4 $ are actually equivalent.

In Examples~\ref{S3computation},~\ref{D8realcomputation} and~\ref{D8CMcomputation} we have $ E=\Q $ so
that we identify $ E\tensor_\Q \C $ with $ \C $, etc., and use notation as in Remark~\ref{idenrem}.

\def\vphant{\vphantom{p_{p_p}}\vphantom{b^{b^b}}}
\def\vphantalt{\vphantom{p_{p_p}}\vphantom{b^b}}
\def\rat[#1]{e(#1,M_\pi)}
\def\ratE[#1]{e(#1,M_\pi^E)}
\def\preg[#1]{R_{#1,p}(M_\pi)/D(M_\pi)^{1/2,p}}
\def\pregE[#1]{R_{#1,p}(M_\pi^E)/D(M_\pi^E)^{1/2,p}}
\newcount\ncount
\def\makenumber#1{\ncount=0\loop\advance\ncount by 1\ifnum\ncount<#1\repeat\advance\ncount by -2 -\number\ncount}
\def\pL[#1]{\Lp p {#1} {\chi_\pi\omega_p^{\makenumber#1}} {\Q} }
\def\pLE[#1]{\Lp p {#1} {\chi_\pi\hskip-1pt\tensor\omega_p^{\makenumber#1}} {\Q} }
\def\quot[#1]{\rat[#1]\cdot \Eul_p(#1,M_\pi)\cdot(\preg[#1])\cdot\pL[#1]^{-1}}
\def\quotE[#1]{\ratE[#1]\cdot \Eul_p(#1,M_\pi^E)\cdot(\pregE[#1])\cdot\pLE[#1]^{-1}}

\long\def\begintable#1\endtable{\begin{table}[htb]#1\end{table}}
\def\beginarray#1\endarray{{\smaller\smaller$$\begin{array}{cc}#1\end{array}$$}}
\def\firstintableheadline[#1]{\vphant \hbox to 10 pt {\hfil $p$\hfil}&\hbox to 330 pt{\hfill$\preg[#1]$\hfill}\cr}
\def\secondintableheadline[#1]{\vphant \hbox to 10 pt {\hfil $p$\hfil}&\hbox to 330 pt{\hfill$\pL[#1]$\hfill}\cr}
\def\thirdintableheadline[#1]{\vphant \hbox to 10 pt {\hfil $p$\hfil}&\hbox to 330 pt{\hfill$\quot[#1]$\hfill}\cr}
\def\firstintableheadlineE[#1]{\vphant \hbox to 10 pt {\hfil $p$\hfil}&\hbox to 330 pt{\hfill$\pregE[#1]$\hfill}\cr}
\def\secondintableheadlineE[#1]{\vphant \hbox to 10 pt {\hfil $p$\hfil}&\hbox to 330 pt{\hfill$\pLE[#1]$\hfill}\cr}
\def\thirdintableheadlineE[#1]{\vphant \hbox to 10 pt {\hfil $p$\hfil}&\hbox to 330 pt{\hfill$\quotE[#1]$\hfill}\cr}
\def\capt[#1]{\caption{\smaller#1}}
\def\labl[#1]{\label{#1}\vskip-20pt}

\font\tablefont=cmcsc9

\begintable
\begin{center}
{\tablefont Table~1.}
{\smaller\smaller Splitting field of $
x^3
 - x^2
 - 3 x
 + 1
$, where $ C(\chi_\pi) = 2^{2}\cdot37$.}
\end{center}
\beginarray
\hline\multicolumn{2}{l}{\vphant n=3\hfill\rat[3]^{-1}=2^{-2}\cdot3^{-2}\cdot37^{2}}\cr
\hline
\firstintableheadline[3]
\hline
\vphantalt 2
&
(1.00010101001000000010000101111100100010101110100110010110110011000000)_{2}\times2^{0}
\cr
\vphantalt 3
&
(1.010012121202211021200001102121121100201201100121200211021212)_{3}\times3^{4}
\cr
\vphantalt 5
&
(4.2320143330214023104113344110103131003140)_{5}\times5^{6}
\cr
\vphantalt 7
&
(6.354304301412412415450326016336635)_{7}\times7^{6}
\cr
\vphantalt 11
&
(2.62161235A928A3423563A7888A)_{11}\times11^{6}
\cr
\hline
\secondintableheadline[3]
\hline
\vphantalt 2
&
(1.00000010100111000100001100100011111010110000101001101110011000101000)_{2}\times2^{2}
\cr
\vphantalt 3
&
(1.1202122002221002110112100012002001211111011101010220)_{3}\times3^{0}
\cr
\vphantalt 5
&
(4.11040024440232442233024131040014140)_{5}\times5^{0}
\cr
\vphantalt 7
&
(5.23516363226501261362543533110)_{7}\times7^{0}
\cr
\vphantalt 11
&
(A.9542728A692401225487A278)_{11}\times11^{0}
\cr
\hline
\hline\multicolumn{2}{l}{\vphant n=5\hfill\rat[5]^{-1}=2^{4}\cdot3^{2}\cdot5^{-2}\cdot37^{4}}\cr
\hline
\firstintableheadline[5]
\hline
\vphantalt 2
&
(1.00001010110110011101000011101001000001100000100010000011011110011101)_{2}\times2^{6}
\cr
\vphantalt 3
&
(1.200111022020102011202110020222210102001020101122110122022)_{3}\times3^{12}
\cr
\vphantalt 5
&
(3.214210342400224303122204142313310423002)_{5}\times5^{8}
\cr
\vphantalt 7
&
(1.2103004405355450123124426604231)_{7}\times7^{10}
\cr
\vphantalt 11
&
(1.2069AA607609242834655465)_{11}\times11^{10}
\cr
\hline
\secondintableheadline[5]
\hline
\vphantalt 2
&
(1.00000101111001110001111000010111110100101010001010001001100111100011)_{2}\times2^{2}
\cr
\vphantalt 3
&
(1.2022122002112102001110211120022002211212022201001210)_{3}\times3^{0}
\cr
\vphantalt 5
&
(3.34404422113121141420412040113130342)_{5}\times5^{0}
\cr
\vphantalt 7
&
(4.66001565513622316645262615350)_{7}\times7^{0}
\cr
\vphantalt 11
&
(1.1A26918431A09200099A81A7)_{11}\times11^{0}
\cr
\hline
\endarray
\smallskip
\begin{center}
{\tablefont Table~2.}
{\smaller\smaller Splitting field of $
x^3
 - 6 x
 - 2
$, where $ C(\chi_\pi) = 2^{2}\cdot3^{3}\cdot7$.}
\end{center}
\beginarray
\hline\multicolumn{2}{l}{\vphant n=3\hfill\rat[3]^{-1}=-2^{-2}\cdot3^{4}\cdot7^{2}}\cr
\hline
\firstintableheadline[3]
\hline
\vphantalt 2
&
(1.10000001110000110101010101111101011001000101111010100001000101111010)_{2}\times2^{0}
\cr
\vphantalt 3
&
(2.1122020011220001210112202100011222122000102121211020022012)_{3}\times3^{4}
\cr
\vphantalt 5
&
(3.4232302100220101433222440213432224202143)_{5}\times5^{6}
\cr
\vphantalt 7
&
(4.2544363600510010013523314315122)_{7}\times7^{7}
\cr
\vphantalt 11
&
(5.481569541633A6875525A00911)_{11}\times11^{6}
\cr
\hline
\secondintableheadline[3]
\hline
\vphantalt 2
&
(1.01111101110011001100110111100111100010100000100110111100101110110010)_{2}\times2^{2}
\cr
\vphantalt 3
&
(1.1211000101210110112021221022022211101101102021011220)_{3}\times3^{0}
\cr
\vphantalt 5
&
(2.13000020122302011221412322140400441)_{5}\times5^{0}
\cr
\vphantalt 7
&
(4.55510152416245405156545556000)_{7}\times7^{2}
\cr
\vphantalt 11
&
(1.14A21577033A19227344806A)_{11}\times11^{0}
\cr
\hline
\hline\multicolumn{2}{l}{\vphant n=5\hfill\rat[5]^{-1}=-2^{4}\cdot3^{14}\cdot5^{-2}\cdot7^{4}}\cr
\hline
\firstintableheadline[5]
\hline
\vphantalt 2
&
(1.10101101100000101100101100100110101101111111000110011100110000011111)_{2}\times2^{6}
\cr
\vphantalt 3
&
(2.120122222011210012110212011212210211022212221020022201)_{3}\times3^{14}
\cr
\vphantalt 5
&
(4.040210142044334334412003042243311012202)_{5}\times5^{8}
\cr
\vphantalt 7
&
(3.2403415331415026502116243110353)_{7}\times7^{9}
\cr
\vphantalt 11
&
(7.160964096503643537518492)_{11}\times11^{10}
\cr
\hline
\secondintableheadline[5]
\hline
\vphantalt 2
&
(1.01010111000010111011111101010101111011001011010111110100101100010111)_{2}\times2^{2}
\cr
\vphantalt 3
&
(1.2002000121100001212101001020022211022201122100111000)_{3}\times3^{0}
\cr
\vphantalt 5
&
(4.44311102333442121212402410333144343)_{5}\times5^{0}
\cr
\vphantalt 7
&
(3.41526030363465662002262442000)_{7}\times7^{0}
\cr
\vphantalt 11
&
(2.687349321A27130485753530)_{11}\times11^{0}
\cr
\hline
\endarray
\endtable

\begintable
\begin{center}
{\tablefont Table~3.}
{\smaller\smaller Splitting field of $
x^3
 - 4 x
 - 1
$, where $ C(\chi_\pi) = 229$.}
\end{center}
\beginarray
\hline\multicolumn{2}{l}{\vphant n=3\hfill\rat[3]^{-1}=-2^{-9}\cdot3^{-2}\cdot229^{2}}\cr
\hline
\firstintableheadline[3]
\hline
\vphantalt 2
&
(1.01110100101010000001111100101001000110110101000001001011100111100010)_{2}\times2^{0}
\cr
\vphantalt 3
&
(2.012012211010121011021101220021012212100110002220200211021221)_{3}\times3^{4}
\cr
\vphantalt 5
&
(4.3433241433301132322003123241441021322342)_{5}\times5^{6}
\cr
\vphantalt 7
&
(6.262401634562505353633061004510116)_{7}\times7^{6}
\cr
\vphantalt 11
&
(3.6717A20401481085265342054)_{11}\times11^{7}
\cr
\hline
\secondintableheadline[3]
\hline
\vphantalt 2
&
(1.01100001011001101110100010111011001000000110000010000100110111111101)_{2}\times2^{3}
\cr
\vphantalt 3
&
(2.2120211202111122100220011211121110002020101200021220)_{3}\times3^{0}
\cr
\vphantalt 5
&
(3.03414332241433441404241012442104412)_{5}\times5^{0}
\cr
\vphantalt 7
&
(3.24033044132523034433054506150)_{7}\times7^{0}
\cr
\vphantalt 11
&
(9.21A61A480A10801760893710)_{11}\times11^{1}
\cr
\hline
\hline\multicolumn{2}{l}{\vphant n=5\hfill\rat[5]^{-1}=-2^{-8}\cdot3^{6}\cdot11\cdot229^{4}}\cr
\hline
\firstintableheadline[5]
\hline
\vphantalt 2
&
(1.11110110101010011011010001011100111001101011001101110111000000010001)_{2}\times2^{5}
\cr
\vphantalt 3
&
(2.222112021012010001201220011021022211112022112001020100)_{3}\times3^{16}
\cr
\vphantalt 5
&
(4.3414223144333221124000410221113420102)_{5}\times5^{10}
\cr
\vphantalt 7
&
(2.5612036605646266635264251053110)_{7}\times7^{10}
\cr
\vphantalt 11
&
(8.5A9220278400758A0193AA6)_{11}\times11^{11}
\cr
\hline
\secondintableheadline[5]
\hline
\vphantalt 2
&
(1.01010000111010100010011010011100110010001101000101101100101100100011)_{2}\times2^{3}
\cr
\vphantalt 3
&
(2.1121000022222000210222202101202022011021200021010220)_{3}\times3^{0}
\cr
\vphantalt 5
&
(4.42011132342110332341130324344440134)_{5}\times5^{0}
\cr
\vphantalt 7
&
(1.25032245164656402236525435500)_{7}\times7^{0}
\cr
\vphantalt 11
&
(5.375705346AA0810959728302)_{11}\times11^{0}
\cr
\hline
\endarray
\smallskip
\begin{center}
{\tablefont Table~4.}
{\smaller\smaller Splitting field of $
x^3
 - 6 x^2
 + 2
$, where $ C(\chi_\pi) = 2^{2}\cdot3^{4}\cdot5$.}
\end{center}
\beginarray
\hline\multicolumn{2}{l}{\vphant n=3\hfill\rat[3]^{-1}=2^{-2}\cdot3^{6}\cdot5}\cr
\hline
\firstintableheadline[3]
\hline
\vphantalt 2
&
(1.11110010110000001011011110010111010011100101111011011000111111110001)_{2}\times2^{0}
\cr
\vphantalt 3
&
(2.120020100000012102010221000201210120200202222221111002100)_{3}\times3^{6}
\cr
\vphantalt 5
&
(3.3432222440014423344440441124111231400414)_{5}\times5^{4}
\cr
\vphantalt 7
&
(5.611651266622656555064166304513241)_{7}\times7^{6}
\cr
\vphantalt 11
&
(2.805A8344265760A69A2490828A)_{11}\times11^{6}
\cr
\hline
\secondintableheadline[3]
\hline
\vphantalt 2
&
(1.10100000100110111110111110101101001001110101011111100111101100001111)_{2}\times2^{2}
\cr
\vphantalt 3
&
(1.1120211121012011012202202010100211210110122222201210)_{3}\times3^{0}
\cr
\vphantalt 5
&
(2.13242313233011201133333143412022134)_{5}\times5^{0}
\cr
\vphantalt 7
&
(3.00123442142613541412624063510)_{7}\times7^{0}
\cr
\vphantalt 11
&
(2.2725A17844929880A8412281)_{11}\times11^{0}
\cr
\hline
\hline\multicolumn{2}{l}{\vphant n=5\hfill\rat[5]^{-1}=-2^{5}\cdot3^{18}\cdot5^{2}\cdot324762301}\cr
\hline
\firstintableheadline[5]
\hline
\vphantalt 2
&
(1.00010100011011001110000011110010000100001101100111101000101001001110)_{2}\times2^{7}
\cr
\vphantalt 3
&
(1.0200220011011111211112101212120212201200221122102001)_{3}\times3^{18}
\cr
\vphantalt 5
&
(1.321100221113142104412032003401031323110)_{5}\times5^{7}
\cr
\vphantalt 7
&
(3.5553063665305405650105042513546)_{7}\times7^{10}
\cr
\vphantalt 11
&
(A.568607725325302A5663A746)_{11}\times11^{10}
\cr
\hline
\secondintableheadline[5]
\hline
\vphantalt 2
&
(1.10111101010011001111001000001111111000110110101100010100000000001010)_{2}\times2^{2}
\cr
\vphantalt 3
&
(1.0002101112212121121001002201100122111211210111102110)_{3}\times3^{0}
\cr
\vphantalt 5
&
(2.30334130213120322410240113023441000)_{5}\times5^{0}
\cr
\vphantalt 7
&
(4.62042101122332630656130062420)_{7}\times7^{0}
\cr
\vphantalt 11
&
(5.491153185997A48326862A46)_{11}\times11^{0}
\cr
\hline
\endarray
\endtable

\begintable
\begin{center}
{\tablefont Table~5.}
{\smaller\smaller Splitting field of $
x^4
 - x^3
 - 3 x^2
 + x
 + 1
$, where $ C(\chi_\pi) = 5\cdot29$.}
\end{center}
\beginarray
\hline\multicolumn{2}{l}{\vphant n=3\hfill\rat[3]^{-1}=2^{-6}\cdot3^{-1}\cdot29^{2}}\cr
\hline
\firstintableheadline[3]
\hline
\vphantalt 2
&
(1.11111010111100011100011000110000111000001111111101101000111001001001)_{2}\times2^{2}
\cr
\vphantalt 3
&
(1.201220211120010110010200120221012101101121011201212202210212)_{3}\times3^{5}
\cr
\vphantalt 5
&
(2.21221421224100432231134432324224211113143)_{5}\times5^{3}
\cr
\vphantalt 7
&
(2.260524146626624605001621634435032)_{7}\times7^{6}
\cr
\vphantalt 11
&
(2.74697855444977862531A25026)_{11}\times11^{6}
\cr
\hline
\secondintableheadline[3]
\hline
\vphantalt 2
&
(1.01010010100111000101001000111001101011100000001100110101111010001100)_{2}\times2^{2}
\cr
\vphantalt 3
&
(1.1011210202202201201020222100121100002100020220220021)_{3}\times3^{0}
\cr
\vphantalt 5
&
(4.30210324423422014142420103340302330)_{5}\times5^{0}
\cr
\vphantalt 7
&
(1.16355063143251352506540201034)_{7}\times7^{0}
\cr
\vphantalt 11
&
(9.9A7720A981020A9A32071790)_{11}\times11^{0}
\cr
\hline
\hline\multicolumn{2}{l}{\vphant n=5\hfill\rat[5]^{-1}=2^{-4}\cdot3^{2}\cdot5^{2}\cdot29^{4}}\cr
\hline
\firstintableheadline[5]
\hline
\vphantalt 2
&
(1.01101101001101100011000000000011101110110010001010001100111111101001)_{2}\times2^{8}
\cr
\vphantalt 3
&
(1.0002011102111201211120012012011221112210212111112002202122)_{3}\times3^{12}
\cr
\vphantalt 5
&
(1.030204001130134123130241412102302410321)_{5}\times5^{7}
\cr
\vphantalt 7
&
(5.6213331545344431613600254300115)_{7}\times7^{10}
\cr
\vphantalt 11
&
(9.98414459608789A37989AA86)_{11}\times11^{10}
\cr
\hline
\secondintableheadline[5]
\hline
\vphantalt 2
&
(1.01111100101100101110101001110010110011011010001011101111000000111001)_{2}\times2^{2}
\cr
\vphantalt 3
&
(1.1111012010001010110000200000021011010010220212100211)_{3}\times3^{0}
\cr
\vphantalt 5
&
(4.01341421411140002024312030010110013)_{5}\times5^{0}
\cr
\vphantalt 7
&
(4.10032134632553222245313526015)_{7}\times7^{0}
\cr
\vphantalt 11
&
(8.6637671480756A9916284377)_{11}\times11^{0}
\cr
\hline
\endarray
\smallskip
\begin{center}
{\tablefont Table~6.}
{\smaller\smaller Splitting field of $
x^4
 - 2 x^3
 - 3 x^2
 + 4 x
 + 1
$, where $ C(\chi_\pi) = 2^{2}\cdot3^{2}\cdot11$.}
\end{center}
\beginarray
\hline\multicolumn{2}{l}{\vphant n=3\hfill\rat[3]^{-1}=2^{-2}\cdot3^{2}\cdot11^{2}}\cr
\hline
\firstintableheadline[3]
\hline
\vphantalt 2
&
(1.11001111111100110011011000011011110100000010001100010000111000001110)_{2}\times2^{5}
\cr
\vphantalt 3
&
(1.101222112200111120002122200210022002021121200101112022212121)_{3}\times3^{2}
\cr
\vphantalt 5
&
(1.3213042203212141302102134014142431100344)_{5}\times5^{6}
\cr
\vphantalt 7
&
(1.444434205533111532005253505615656)_{7}\times7^{6}
\cr
\vphantalt 11
&
(9.61613555149153993505670160)_{11}\times11^{5}
\cr
\hline
\secondintableheadline[3]
\hline
\vphantalt 2
&
(1.11110101101001001100111001001100001010010111011100011010101000010001)_{2}\times2^{4}
\cr
\vphantalt 3
&
(1.1221211120021010011021120022021022220202012220120020)_{3}\times3^{0}
\cr
\vphantalt 5
&
(1.03123432421241311210024333243134200)_{5}\times5^{0}
\cr
\vphantalt 7
&
(1.20515104501535164355503412253)_{7}\times7^{0}
\cr
\vphantalt 11
&
(4.A19905A97667183A61256312)_{11}\times11^{0}
\cr
\hline
\hline\multicolumn{2}{l}{\vphant n=5\hfill\rat[5]^{-1}=2^{7}\cdot3^{10}\cdot5^{-2}\cdot11^{7}\cdot151\cdot1389251}\cr
\hline
\firstintableheadline[5]
\hline
\vphantalt 2
&
(1.10101000100011001001110000000110011010000011010001100100101100100011)_{2}\times2^{16}
\cr
\vphantalt 3
&
(2.00020011020011210121200021121000100102202022102012021120)_{3}\times3^{10}
\cr
\vphantalt 5
&
(3.043124340101142041302443131432222411101)_{5}\times5^{8}
\cr
\vphantalt 7
&
(3.1016042121151553442503432550426)_{7}\times7^{10}
\cr
\vphantalt 11
&
(4.A244AA713986547114383)_{11}\times11^{13}
\cr
\hline
\secondintableheadline[5]
\hline
\vphantalt 2
&
(1.01011001111000100100100110110111001110010101101110101100000001001101)_{2}\times2^{4}
\cr
\vphantalt 3
&
(1.2200111010220200001122120000000211211201112212022102)_{3}\times3^{0}
\cr
\vphantalt 5
&
(4.31111224412403210010123414333414223)_{5}\times5^{0}
\cr
\vphantalt 7
&
(2.63463651641424006321031014001)_{7}\times7^{0}
\cr
\vphantalt 11
&
(2.054143779893294A023A2760)_{11}\times11^{1}
\cr
\hline
\endarray
\endtable

\begintable
\begin{center}
{\tablefont Table~7.}
{\smaller\smaller Splitting field of $
x^4
 - 6 x^2
 + 6
$, where $ C(\chi_\pi) = 2^{7}\cdot3^{2}$.}
\end{center}
\beginarray
\hline\multicolumn{2}{l}{\vphant n=3\hfill\rat[3]^{-1}=-2^{8}\cdot3^{2}}\cr
\hline
\firstintableheadline[3]
\hline
\vphantalt 2
&
(1.11111110110101110110001001001000100111010011011111000101011110110110)_{2}\times2^{10}
\cr
\vphantalt 3
&
(2.022112012020121000012010222021210101120000211220112222110100)_{3}\times3^{2}
\cr
\vphantalt 5
&
(2.041234240331110012313440303434111221114)_{5}\times5^{7}
\cr
\vphantalt 7
&
(4.126126412601443066506634321456330)_{7}\times7^{6}
\cr
\vphantalt 11
&
(1.508753A40028580A83265A6397)_{11}\times11^{6}
\cr
\hline
\secondintableheadline[3]
\hline
\vphantalt 2
&
(1.00111001111010100110000110000111111000001101110001010000100110100110)_{2}\times2^{2}
\cr
\vphantalt 3
&
(1.1102020120122110111012211100112220010211120222001022)_{3}\times3^{0}
\cr
\vphantalt 5
&
(3.24102423421030322330312322041041320)_{5}\times5^{1}
\cr
\vphantalt 7
&
(3.53421553111552403452502135351)_{7}\times7^{0}
\cr
\vphantalt 11
&
(9.993168A73395501045A926A6)_{11}\times11^{0}
\cr
\hline
\hline\multicolumn{2}{l}{\vphant n=5\hfill\rat[5]^{-1}=-2^{25}\cdot3^{10}\cdot5^{-1}\cdot11\cdot37\cdot180097}\cr
\hline
\firstintableheadline[5]
\hline
\vphantalt 2
&
(1.01000011110101000111110110010101100000101100100001110011010100110111)_{2}\times2^{27}
\cr
\vphantalt 3
&
(1.20110210022121202120201021000120210110011111202222221202)_{3}\times3^{10}
\cr
\vphantalt 5
&
(4.13321010222213321143100244321412404344)_{5}\times5^{9}
\cr
\vphantalt 7
&
(1.6542260045104311315015635541014)_{7}\times7^{10}
\cr
\vphantalt 11
&
(4.35786968692892127976707)_{11}\times11^{11}
\cr
\hline
\secondintableheadline[5]
\hline
\vphantalt 2
&
(1.00100100101111111001101001011101011111110110100110100000111111010001)_{2}\times2^{2}
\cr
\vphantalt 3
&
(1.0101222010121001211100001210022220120122121200102122)_{3}\times3^{0}
\cr
\vphantalt 5
&
(2.44334240213144033340330341233100231)_{5}\times5^{0}
\cr
\vphantalt 7
&
(2.43613231343555620664002413300)_{7}\times7^{0}
\cr
\vphantalt 11
&
(A.957125591551731675713855)_{11}\times11^{0}
\cr
\hline
\endarray
\smallskip
\begin{center}
{\tablefont Table~8.}
{\smaller\smaller Splitting field of $
x^4
 - 6 x^2
 - 4 x
 + 2
$, where $ C(\chi_\pi) = 2^{5}\cdot17$.}
\end{center}
\beginarray
\hline\multicolumn{2}{l}{\vphant n=3\hfill\rat[3]^{-1}=-2^{2}\cdot3^{-2}\cdot17^{2}}\cr
\hline
\firstintableheadline[3]
\hline
\vphantalt 2
&
(1.10111000110111011100110001011111110001011110100001100011100001101001)_{2}\times2^{7}
\cr
\vphantalt 3
&
(1.1001122110011102222210022111010102000101002122001210120111212)_{3}\times3^{4}
\cr
\vphantalt 5
&
(2.1330133431400143111041343243044132341000)_{5}\times5^{6}
\cr
\vphantalt 7
&
(1.251366021640656215554304534633666)_{7}\times7^{6}
\cr
\vphantalt 11
&
(6.8781153571182A7A3A52578590)_{11}\times11^{6}
\cr
\hline
\secondintableheadline[3]
\hline
\vphantalt 2
&
(1.01100101000011011110000111011100100110110110111011000100000101100110)_{2}\times2^{5}
\cr
\vphantalt 3
&
(2.2111111121010122120102220010120102202210101102021112)_{3}\times3^{0}
\cr
\vphantalt 5
&
(2.03010131243321342244201240202313243)_{5}\times5^{0}
\cr
\vphantalt 7
&
(2.45653546525345623002012601006)_{7}\times7^{0}
\cr
\vphantalt 11
&
(1.5054539309942A8686113521)_{11}\times11^{0}
\cr
\hline
\hline\multicolumn{2}{l}{\vphant n=5\hfill\rat[5]^{-1}=-2^{14}\cdot3^{2}\cdot5^{-2}\cdot17^{4}}\cr
\hline
\firstintableheadline[5]
\hline
\vphantalt 2
&
(1.01101001100111100011100000011101100011011111111011011010000000100111)_{2}\times2^{19}
\cr
\vphantalt 3
&
(1.0022021200111200200110022100110101220100021201020020021001)_{3}\times3^{12}
\cr
\vphantalt 5
&
(4.41122221322441142131444023001440141102)_{5}\times5^{9}
\cr
\vphantalt 7
&
(4.6425656621000232634542266362653)_{7}\times7^{10}
\cr
\vphantalt 11
&
(5.30669526A542413A15559106)_{11}\times11^{10}
\cr
\hline
\secondintableheadline[5]
\hline
\vphantalt 2
&
(1.10000111001100100001010000010001000011100111101011000000110111100101)_{2}\times2^{5}
\cr
\vphantalt 3
&
(2.1100001122200222122120000002011201000121010120122102)_{3}\times3^{0}
\cr
\vphantalt 5
&
(1.00313434131411434403030110103232430)_{5}\times5^{1}
\cr
\vphantalt 7
&
(4.51404541504052312152404613255)_{7}\times7^{0}
\cr
\vphantalt 11
&
(2.A43315061300727A39739622)_{11}\times11^{0}
\cr
\hline
\endarray
\endtable

\begintable
\begin{center}
{\tablefont Table~9.}
{\smaller\smaller Splitting field of $
x^4
 - 2 x^3
 + 5 x^2
 - 4 x
 + 2
$, where $ C(\chi_\pi) = 2^{3}\cdot17$.}
\end{center}
\beginarray
\hline\multicolumn{2}{l}{\vphant n=2\hfill\rat[2]^{-1}=-17}\cr
\hline
\firstintableheadline[2]
\hline
\vphantalt 2
&
(1.11110000011111000100011010100111010001000101111111000010111001000001)_{2}\times2^{8}
\cr
\vphantalt 3
&
(1.2001010102220112101111101220000000012020211022222012112010122)_{3}\times3^{4}
\cr
\vphantalt 5
&
(2.01203020434342403202140220421203111323421)_{5}\times5^{4}
\cr
\vphantalt 7
&
(4.1140305454665411023106103225364300)_{7}\times7^{4}
\cr
\vphantalt 11
&
(9.074667649A54466861880828A41)_{11}\times11^{4}
\cr
\hline
\secondintableheadline[2]
\hline
\vphantalt 2
&
(1.10011100010001001110011010011101110110111100101101100011101101010011)_{2}\times2^{6}
\cr
\vphantalt 3
&
(1.2211021011200210000101010021221011201210202220121211)_{3}\times3^{0}
\cr
\vphantalt 5
&
(4.33003403101112402103034043124123020)_{5}\times5^{0}
\cr
\vphantalt 7
&
(6.55655254400302602216266140200)_{7}\times7^{0}
\cr
\vphantalt 11
&
(4.84968AAA0466783629843316)_{11}\times11^{0}
\cr
\hline
\hline\multicolumn{2}{l}{\vphant n=4\hfill\rat[4]^{-1}=2^{4}\cdot3^{2}\cdot17^{3}}\cr
\hline
\firstintableheadline[4]
\hline
\vphantalt 2
&
(1.01111010101011110100100110010011000001110001101011101000010011010010)_{2}\times2^{15}
\cr
\vphantalt 3
&
(2.02020202222222212011121000222200020002221002201220201002211)_{3}\times3^{10}
\cr
\vphantalt 5
&
(3.121141434441031100040041101320123332342)_{5}\times5^{8}
\cr
\vphantalt 7
&
(4.30103246365105165401451263010636)_{7}\times7^{8}
\cr
\vphantalt 11
&
(9.8864321414815928A83132426)_{11}\times11^{8}
\cr
\hline
\secondintableheadline[4]
\hline
\vphantalt 2
&
(1.10101000111010101100110000110001001011001001001001010011001010011101)_{2}\times2^{7}
\cr
\vphantalt 3
&
(1.0220110000101201002201220101112220112012211122211111)_{3}\times3^{0}
\cr
\vphantalt 5
&
(4.24124003322330023340322023204112410)_{5}\times5^{0}
\cr
\vphantalt 7
&
(1.34631136124161206563156313232)_{7}\times7^{0}
\cr
\vphantalt 11
&
(6.484018160860423737043861)_{11}\times11^{0}
\cr
\hline
\endarray
\smallskip
\begin{center}
{\tablefont Table~10.}
{\smaller\smaller Splitting field of $
x^4
 - x^3
 + 3 x^2
 - 2 x
 + 4
$, where $ C(\chi_\pi) = 5\cdot41$.}
\end{center}
\beginarray
\hline\multicolumn{2}{l}{\vphant n=2\hfill\rat[2]^{-1}=2^{-3}\cdot5\cdot41}\cr
\hline
\firstintableheadline[2]
\hline
\vphantalt 2
&
(1.00100101101010001111110110001011010111110011100101011001110011011000)_{2}\times2^{6}
\cr
\vphantalt 3
&
(2.0021220111220210202122211101102111121010022202210012111112202)_{3}\times3^{4}
\cr
\vphantalt 5
&
(4.2413230212134320403404223440242144143010)_{5}\times5^{4}
\cr
\vphantalt 7
&
(5.6622665100264660251241453165220553)_{7}\times7^{4}
\cr
\vphantalt 11
&
(2.0A380682536A50856633A5963AA)_{11}\times11^{4}
\cr
\hline
\secondintableheadline[2]
\hline
\vphantalt 2
&
(1.10001100100010010110110100011001000110001101100111111011111110000000)_{2}\times2^{5}
\cr
\vphantalt 3
&
(1.0000112110000221101100101100220112122121100200111220)_{3}\times3^{0}
\cr
\vphantalt 5
&
(3.31322324200110323444242424141204000)_{5}\times5^{1}
\cr
\vphantalt 7
&
(6.20642211231410105634103112315)_{7}\times7^{0}
\cr
\vphantalt 11
&
(4.8488929AA248277415458863)_{11}\times11^{0}
\cr
\hline
\hline\multicolumn{2}{l}{\vphant n=4\hfill\rat[4]^{-1}=2^{-5}\cdot3^{2}\cdot5^{3}\cdot41^{3}}\cr
\hline
\firstintableheadline[4]
\hline
\vphantalt 2
&
(1.01101011011011001001110001011011001001100001000100001111001011100010)_{2}\times2^{9}
\cr
\vphantalt 3
&
(2.12201221002111122111110002102122122100002111212000210102220)_{3}\times3^{10}
\cr
\vphantalt 5
&
(3.03433213241214232102131324212010224400)_{5}\times5^{8}
\cr
\vphantalt 7
&
(3.21341424531605465104510422214051)_{7}\times7^{8}
\cr
\vphantalt 11
&
(A.180567AA0545A59652909A250)_{11}\times11^{8}
\cr
\hline
\secondintableheadline[4]
\hline
\vphantalt 2
&
(1.11100001111100101110011000101000001100001100100001000011000001000111)_{2}\times2^{6}
\cr
\vphantalt 3
&
(1.2110210100202111000001112120121020120200000202111102)_{3}\times3^{0}
\cr
\vphantalt 5
&
(1.22414201312142130131112334430123200)_{5}\times5^{1}
\cr
\vphantalt 7
&
(6.66200605216643655434600606423)_{7}\times7^{0}
\cr
\vphantalt 11
&
(3.5059939560A2AA91106844A7)_{11}\times11^{0}
\cr
\hline
\endarray
\endtable

\begintable
\begin{center}
{\tablefont Table~11.}
{\smaller\smaller Splitting field of $
x^4
 - 2 x^3
 + 6 x^2
 - 5 x
 + 2
$, where $ C(\chi_\pi) = 13\cdot17$.}
\end{center}
\beginarray
\hline\multicolumn{2}{l}{\vphant n=2\hfill\rat[2]^{-1}=2^{-3}\cdot13\cdot17}\cr
\hline
\firstintableheadline[2]
\hline
\vphantalt 2
&
(1.11100111011011110110000110110011010010001100010000011011000000100011)_{2}\times2^{7}
\cr
\vphantalt 3
&
(1.001110110101022221010202101001220021110111221112100200010011)_{3}\times3^{5}
\cr
\vphantalt 5
&
(4.12002342224442011223100142021013114003331)_{5}\times5^{4}
\cr
\vphantalt 7
&
(6.0065051533250322312551556340446112)_{7}\times7^{4}
\cr
\vphantalt 11
&
(9.58602382412A681650341722097)_{11}\times11^{4}
\cr
\hline
\secondintableheadline[2]
\hline
\vphantalt 2
&
(1.01010001000000000001000000101011111111001100001000001100111111100100)_{2}\times2^{6}
\cr
\vphantalt 3
&
(2.0212122200002021202200211122100200001001012201202000)_{3}\times3^{1}
\cr
\vphantalt 5
&
(2.12300004023110442141103112001044200)_{5}\times5^{0}
\cr
\vphantalt 7
&
(5.65035644055466446205106662353)_{7}\times7^{0}
\cr
\vphantalt 11
&
(6.38803A213738A731776470A4)_{11}\times11^{0}
\cr
\hline
\hline\multicolumn{2}{l}{\vphant n=4\hfill\rat[4]^{-1}=-2^{-5}\cdot3^{2}\cdot13^{3}\cdot17^{3}}\cr
\hline
\firstintableheadline[4]
\hline
\vphantalt 2
&
(1.01010110010100100110100011110011010001011010001010001001110111101010)_{2}\times2^{10}
\cr
\vphantalt 3
&
(1.1111010110110112222112002001211122102011200022210022220110)_{3}\times3^{11}
\cr
\vphantalt 5
&
(2.300041000440134413002331411231121430312)_{5}\times5^{8}
\cr
\vphantalt 7
&
(2.63246241360463115315231424221502)_{7}\times7^{8}
\cr
\vphantalt 11
&
(A.49299953329A1125214A939A7)_{11}\times11^{8}
\cr
\hline
\secondintableheadline[4]
\hline
\vphantalt 2
&
(1.00100111010100001110011000011111110110110001001100110010111000011110)_{2}\times2^{7}
\cr
\vphantalt 3
&
(1.2001110201210120022112022021221001220010010202222220)_{3}\times3^{1}
\cr
\vphantalt 5
&
(4.31111034231223303224440300102101300)_{5}\times5^{0}
\cr
\vphantalt 7
&
(3.23250645010664211233114156026)_{7}\times7^{0}
\cr
\vphantalt 11
&
(6.722A82921527466612856704)_{11}\times11^{0}
\cr
\hline
\endarray
\smallskip
\begin{center}
{\tablefont Table~12.}
{\smaller\smaller Splitting field of $
x^4
 + 3 x^2
 - 6 x
 + 6
$, where $ C(\chi_\pi) = 2^{2}\cdot3^{2}\cdot11$.}
\end{center}
\beginarray
\hline\multicolumn{2}{l}{\vphant n=2\hfill\rat[2]^{-1}=-2^{-1}\cdot3^{2}\cdot11}\cr
\hline
\firstintableheadline[2]
\hline
\vphantalt 2
&
(1.00011111110011010111101100110001100000010001000010010011001000110100)_{2}\times2^{6}
\cr
\vphantalt 3
&
(2.001121000222211121002102201101101100211202000001120001201101)_{3}\times3^{2}
\cr
\vphantalt 5
&
(4.32211042044110431112434230011130322040334)_{5}\times5^{4}
\cr
\vphantalt 7
&
(6.460003406606016165433631101363300)_{7}\times7^{5}
\cr
\vphantalt 11
&
(4.355421A54525023347874A49229)_{11}\times11^{3}
\cr
\hline
\secondintableheadline[2]
\hline
\vphantalt 2
&
(1.11101000100101000011100111111110011101011000001101111010111100101100)_{2}\times2^{5}
\cr
\vphantalt 3
&
(1.2022010222110011100120101011011011020010020011221012)_{3}\times3^{0}
\cr
\vphantalt 5
&
(3.22012020201133014223242323022213400)_{5}\times5^{0}
\cr
\vphantalt 7
&
(2.43314511556206561022451602200)_{7}\times7^{1}
\cr
\vphantalt 11
&
(7.514814537963976728A46001)_{11}\times11^{0}
\cr
\hline
\hline\multicolumn{2}{l}{\vphant n=4\hfill\rat[4]^{-1}=-2\cdot3^{8}\cdot11^{3}}\cr
\hline
\firstintableheadline[4]
\hline
\vphantalt 2
&
(1.00011101001000000111110111000101001111110011011101001101001001110101)_{2}\times2^{11}
\cr
\vphantalt 3
&
(2.0222010120001012000101112011222120000210000112021000101001)_{3}\times3^{8}
\cr
\vphantalt 5
&
(4.200024100322411004323111441041141020122)_{5}\times5^{8}
\cr
\vphantalt 7
&
(5.12536440256103300364331512021150)_{7}\times7^{8}
\cr
\vphantalt 11
&
(6.6966231827671A91984A434A7)_{11}\times11^{7}
\cr
\hline
\secondintableheadline[4]
\hline
\vphantalt 2
&
(1.10111000101101010100110110100100101011000111100000010010100110001101)_{2}\times2^{6}
\cr
\vphantalt 3
&
(1.0221110122220100011212020012222021100000010022110011)_{3}\times3^{0}
\cr
\vphantalt 5
&
(3.44200104443132341021214030231112200)_{5}\times5^{0}
\cr
\vphantalt 7
&
(4.44003540500226200364534021346)_{7}\times7^{0}
\cr
\vphantalt 11
&
(5.766584290652809A28501199)_{11}\times11^{0}
\cr
\hline
\endarray
\endtable

\begintable
\begin{center}
{\tablefont Table~13.}
{\smaller\smaller Splitting field of $
x^6
 - x^5
 - 6 x^4
 + 7 x^3
 + 4 x^2
 - 5 x
 + 1
$, where $ C(\chi_\pi\tensor\id) = 2^{2}\cdot5\cdot19$.}
\end{center}
\beginarray
\hline\multicolumn{2}{l}{\vphant n=3\hfill\ratE[3]^{-1}=3^{-2}\cdot19^{2}}\cr
\hline
\firstintableheadlineE[3]
\hline
\vphantalt 2
&
1\tensor ((1.010101011100101010)_{2}\times2^{4}) + \z_3 \tensor ((1.111100111000011101)_{2}\times2^{2})
\cr
\vphantalt 3
&
1\tensor ((2.121122121210110121)_{3}\times3^{4}) + \z_3 \tensor ((1.100202121120121202)_{3}\times3^{4})
\cr
\vphantalt 5
&
1\tensor ((4.402020001241122012)_{5}\times5^{3}) + \z_3 \tensor ((1.402102204113114001)_{5}\times5^{4})
\cr
\vphantalt 7
&
1\tensor ((2.1404154326422162)_{7}\times7^{6}) + \z_3 \tensor ((5.1444362630564332)_{7}\times7^{6})
\cr
\vphantalt 11
&
1\tensor ((6.9026829596822)_{11}\times11^{6}) + \z_3 \tensor ((8.048680429402A)_{11}\times11^{6})
\cr
\hline
\secondintableheadlineE[3]
\hline
\vphantalt 2
&
1\tensor ((1.010111010101101000)_{2}\times2^{4}) + \z_3 \tensor ((1.111111110001001011)_{2}\times2^{2})
\cr
\vphantalt 3
&
1\tensor ((1.120012101022112222)_{3}\times3^{0}) + \z_3 \tensor ((2.212202101110022110)_{3}\times3^{0})
\cr
\vphantalt 5
&
1\tensor ((1.222321004420021231)_{5}\times5^{0}) + \z_3 \tensor ((1.123432442320044212)_{5}\times5^{0})
\cr
\vphantalt 7
&
1\tensor ((6.215644553551135401)_{7}\times7^{0}) + \z_3 \tensor ((1.642236344041330514)_{7}\times7^{0})
\cr
\vphantalt 11
&
1\tensor ((6.6A0036761338492361)_{11}\times11^{0}) + \z_3 \tensor ((8.749011AA3694540145)_{11}\times11^{0})
\cr
\hline
\hline\multicolumn{2}{l}{\vphant n=5\hfill\ratE[5]^{-1}=2^{7}\cdot3^{2}\cdot5^{2}\cdot19^{4}\cdot(37754+43673\z_3)}\cr
\hline
\firstintableheadlineE[5]
\hline
\vphantalt 2
&
1\tensor ((1.000011100001001101)_{2}\times2^{9}) + \z_3 \tensor ((1.101010111001110110)_{2}\times2^{9})
\cr
\vphantalt 3
&
1\tensor ((1.110110120010021222)_{3}\times3^{12}) + \z_3 \tensor ((2.020220021122220210)_{3}\times3^{12})
\cr
\vphantalt 5
&
1\tensor ((3.020013014242432200)_{5}\times5^{7}) + \z_3 \tensor ((1.220134043203242432)_{5}\times5^{7})
\cr
\vphantalt 7
&
1\tensor ((4.15141551064450)_{7}\times7^{10}) + \z_3 \tensor ((6.53120060244502)_{7}\times7^{10})
\cr
\vphantalt 11
&
1\tensor ((4.285500134)_{11}\times11^{12}) + \z_3 \tensor ((2.A3934821906)_{11}\times11^{10})
\cr
\hline
\secondintableheadlineE[5]
\hline
\vphantalt 2
&
1\tensor ((1.110010010001000101)_{2}\times2^{4}) + \z_3 \tensor ((1.111111100101010110)_{2}\times2^{2})
\cr
\vphantalt 3
&
1\tensor ((1.211021010101222221)_{3}\times3^{0}) + \z_3 \tensor ((2.122212200211101011)_{3}\times3^{0})
\cr
\vphantalt 5
&
1\tensor ((1.401222401123304240)_{5}\times5^{0}) + \z_3 \tensor ((1.432224122043203014)_{5}\times5^{0})
\cr
\vphantalt 7
&
1\tensor ((2.441442356342000541)_{7}\times7^{0}) + \z_3 \tensor ((3.664531413660053644)_{7}\times7^{0})
\cr
\vphantalt 11
&
1\tensor ((6.3A4A10A60052302518)_{11}\times11^{0}) + \z_3 \tensor ((4.7AA8843572652457A4)_{11}\times11^{0})
\cr
\hline
\endarray
\smallskip
\begin{center}
{\tablefont Table~14.}
{\smaller\smaller Splitting field of $
x^6
 - x^5
 - 8 x^4
 - x^3
 + 12 x^2
 + 7 x
 + 1
$, where $ C(\chi_\pi\tensor\id) = 2^{2}\cdot5\cdot31$.}
\end{center}
\beginarray
\hline\multicolumn{2}{l}{\vphant n=3\hfill\ratE[3]^{-1}= 2^{-2}\cdot3^{-2}\cdot31^{2}\cdot \z_3}\cr
\hline
\firstintableheadlineE[3]
\hline
\vphantalt 2
&
1\tensor ((1.100100010111101000)_{2}\times2^{2}) + \z_3 \tensor ((1.111010100101111010)_{2}\times2^{0})
\cr
\vphantalt 3
&
1\tensor ((1.221220200011201020)_{3}\times3^{4}) + \z_3 \tensor ((2.020111011200022220)_{3}\times3^{4})
\cr
\vphantalt 5
&
1\tensor ((1.141130342002111122)_{5}\times5^{3}) + \z_3 \tensor ((4.342032011212400441)_{5}\times5^{3})
\cr
\vphantalt 7
&
1\tensor ((4.4505556355566640)_{7}\times7^{6}) + \z_3 \tensor ((3.0610023664215014)_{7}\times7^{6})
\cr
\vphantalt 11
&
1\tensor ((7.3700473650A19)_{11}\times11^{6}) + \z_3 \tensor ((1.3659A52396880)_{11}\times11^{6})
\cr
\hline
\secondintableheadlineE[3]
\hline
\vphantalt 2
&
1\tensor ((1.101110001100101000)_{2}\times2^{2}) + \z_3 \tensor ((1.010100111010111101)_{2}\times2^{4})
\cr
\vphantalt 3
&
1\tensor ((2.111112222201000112)_{3}\times3^{0}) + \z_3 \tensor ((1.120112022212020122)_{3}\times3^{0})
\cr
\vphantalt 5
&
1\tensor ((4.321142113321234330)_{5}\times5^{0}) + \z_3 \tensor ((1.130134403242024422)_{5}\times5^{0})
\cr
\vphantalt 7
&
1\tensor ((5.530006300632330553)_{7}\times7^{0}) + \z_3 \tensor ((2.200601212403230506)_{7}\times7^{0})
\cr
\vphantalt 11
&
1\tensor ((2.00762927347972431A)_{11}\times11^{0}) + \z_3 \tensor ((A.191A1780051997A650)_{11}\times11^{0})
\cr
\hline
\hline\multicolumn{2}{l}{\vphant n=5\hfill\ratE[5]^{-1}=-2^{6}\cdot3^{3}\cdot5^{2}\cdot31^{4}\cdot(2+\z_3)\cdot(2-\z_3)\cdot(7-15\z_3)\cdot(54+31\z_3)}\cr
\hline
\firstintableheadlineE[5]
\hline
\vphantalt 2
&
1\tensor ((1.011100011000010010)_{2}\times2^{8}) + \z_3 \tensor ((1.011100000000101011)_{2}\times2^{8})
\cr
\vphantalt 3
&
1\tensor ((2.202202022111202222)_{3}\times3^{14}) + \z_3 \tensor ((1.112002110200111001)_{3}\times3^{15})
\cr
\vphantalt 5
&
1\tensor ((3.141433024141141322)_{5}\times5^{8}) + \z_3 \tensor ((2.002014030014311244)_{5}\times5^{7})
\cr
\vphantalt 7
&
1\tensor ((6.50135302661002)_{7}\times7^{10}) + \z_3 \tensor ((4.04356114536113)_{7}\times7^{10})
\cr
\vphantalt 11
&
1\tensor ((5.A3306A68921)_{11}\times11^{10}) + \z_3 \tensor ((3.85A6650344A)_{11}\times11^{10})
\cr
\hline
\secondintableheadlineE[5]
\hline
\vphantalt 2
&
1\tensor ((1.100110110011011000)_{2}\times2^{2}) + \z_3 \tensor ((1.110110100001011101)_{2}\times2^{4})
\cr
\vphantalt 3
&
1\tensor ((2.201101211210020202)_{3}\times3^{0}) + \z_3 \tensor ((1.212021202011010201)_{3}\times3^{0})
\cr
\vphantalt 5
&
1\tensor ((4.011013012410044000)_{5}\times5^{0}) + \z_3 \tensor ((1.431333301302323324)_{5}\times5^{0})
\cr
\vphantalt 7
&
1\tensor ((2.130552365243450401)_{7}\times7^{0}) + \z_3 \tensor ((2.601221155363314062)_{7}\times7^{0})
\cr
\vphantalt 11
&
1\tensor ((9.53274A487A14818064)_{11}\times11^{0}) + \z_3 \tensor ((9.89A625650825167861)_{11}\times11^{0})
\cr
\hline
\endarray
\endtable

\begintable
\begin{center}
{\tablefont Table~15.}
{\smaller\smaller Splitting field of $
x^6
 - x^5
 - 8 x^4
 + 5 x^3
 + 19 x^2
 - 4 x
 - 11
$, where $ C(\chi_\pi\tensor\id) = 5\cdot139$.}
\end{center}
\beginarray
\hline\multicolumn{2}{l}{\vphant n=3\hfill\ratE[3]^{-1}=2^{-5}\cdot3^{-2}\cdot139^{2}}\cr
\hline
\firstintableheadlineE[3]
\hline
\vphantalt 2
&
1\tensor ((1.000110100111100100)_{2}\times2^{4}) + \z_3 \tensor ((1.000011011101110010)_{2}\times2^{5})
\cr
\vphantalt 3
&
1\tensor ((2.102212010212021001)_{3}\times3^{4}) + \z_3 \tensor ((2.122111112211000220)_{3}\times3^{4})
\cr
\vphantalt 5
&
1\tensor ((1.400211310114320013)_{5}\times5^{3}) + \z_3 \tensor ((1.434230430334022423)_{5}\times5^{3})
\cr
\vphantalt 7
&
1\tensor ((2.2523154003106061)_{7}\times7^{6}) + \z_3 \tensor ((2.3200021260361110)_{7}\times7^{6})
\cr
\vphantalt 11
&
1\tensor ((9.4472801A61A44)_{11}\times11^{6}) + \z_3 \tensor ((1.7232221727046)_{11}\times11^{6})
\cr
\hline
\secondintableheadlineE[3]
\hline
\vphantalt 2
&
1\tensor ((1.000100101010000101)_{2}\times2^{4}) + \z_3 \tensor ((1.000010100000000110)_{2}\times2^{3})
\cr
\vphantalt 3
&
1\tensor ((2.111012101221111220)_{3}\times3^{2}) + \z_3 \tensor ((1.010011111201220212)_{3}\times3^{0})
\cr
\vphantalt 5
&
1\tensor ((3.240243311024013404)_{5}\times5^{0}) + \z_3 \tensor ((4.142342004244233331)_{5}\times5^{2})
\cr
\vphantalt 7
&
1\tensor ((6.543343430631523116)_{7}\times7^{1}) + \z_3 \tensor ((2.342260265363120032)_{7}\times7^{0})
\cr
\vphantalt 11
&
1\tensor ((4.8521241A647A883282)_{11}\times11^{0}) + \z_3 \tensor ((1.46220220356163464A)_{11}\times11^{0})
\cr
\hline
\hline\multicolumn{2}{l}{\vphant n=5\hfill\ratE[5]^{-1}=2^{-4}\cdot3^{2}\cdot5^{2}\cdot139^{4}\cdot(4 - 5 \z_3)}\cr
\hline
\firstintableheadlineE[5]
\hline
\vphantalt 2
&
1\tensor ((1.000100111100110101)_{2}\times2^{10}) + \z_3 \tensor ((1.011111001011110010)_{2}\times2^{9})
\cr
\vphantalt 3
&
1\tensor ((1.202001002112012221)_{3}\times3^{13}) + \z_3 \tensor ((1.021221211222021020)_{3}\times3^{12})
\cr
\vphantalt 5
&
1\tensor ((3.404332041431011004)_{5}\times5^{7}) + \z_3 \tensor ((3.213210304101144004)_{5}\times5^{7})
\cr
\vphantalt 7
&
1\tensor ((4.11031052300040)_{7}\times7^{10}) + \z_3 \tensor ((4.03411453000015)_{7}\times7^{10})
\cr
\vphantalt 11
&
1\tensor ((6.8523188534)_{11}\times11^{11}) + \z_3 \tensor ((5.2A39A447198)_{11}\times11^{10})
\cr
\hline
\secondintableheadlineE[5]
\hline
\vphantalt 2
&
1\tensor ((1.101001011001000000)_{2}\times2^{4}) + \z_3 \tensor ((1.000100011110101010)_{2}\times2^{3})
\cr
\vphantalt 3
&
1\tensor ((1.200212010102010102)_{3}\times3^{1}) + \z_3 \tensor ((1.122011112111021002)_{3}\times3^{0})
\cr
\vphantalt 5
&
1\tensor ((3.124113213204101323)_{5}\times5^{0}) + \z_3 \tensor ((2.114244212313224030)_{5}\times5^{3})
\cr
\vphantalt 7
&
1\tensor ((6.430024313466045624)_{7}\times7^{0}) + \z_3 \tensor ((5.065261416212661354)_{7}\times7^{0})
\cr
\vphantalt 11
&
1\tensor ((5.A88931A0106A0172A3)_{11}\times11^{0}) + \z_3 \tensor ((3.A4151A485571789978)_{11}\times11^{0})
\cr
\hline
\endarray
\smallskip
\begin{center}
{\tablefont Table~16.}
{\smaller\smaller Splitting field of $
x^6
 - x^5
 - 8 x^4
 + 11 x^3
 + 2 x^2
 - 5 x
 + 1
$, where $ C(\chi_\pi\tensor\id) = 2^{2}\cdot3^{2}\cdot13$.}
\end{center}
\beginarray
\hline\multicolumn{2}{l}{\vphant n=3\hfill\ratE[3]^{-1}= -2^{-2}\cdot3^{2}\cdot13^{2}\cdot \z_3}\cr
\hline
\firstintableheadlineE[3]
\hline
\vphantalt 2
&
1\tensor ((1.101100110101111100)_{2}\times2^{2}) + \z_3 \tensor ((1.001000001010110111)_{2}\times2^{0})
\cr
\vphantalt 3
&
1\tensor ((2.021221201100120120)_{3}\times3^{5}) + \z_3 \tensor ((2.020012210111121212)_{3}\times3^{7})
\cr
\vphantalt 5
&
1\tensor ((3.000431031211204434)_{5}\times5^{6}) + \z_3 \tensor ((1.012314421212444244)_{5}\times5^{6})
\cr
\vphantalt 7
&
1\tensor ((3.0426445363200511)_{7}\times7^{6}) + \z_3 \tensor ((1.1035504615611621)_{7}\times7^{6})
\cr
\vphantalt 11
&
1\tensor ((8.9799A36158285)_{11}\times11^{6}) + \z_3 \tensor ((9.229923AAA0734)_{11}\times11^{6})
\cr
\hline
\secondintableheadlineE[3]
\hline
\vphantalt 2
&
1\tensor ((1.100100111001100011)_{2}\times2^{2}) + \z_3 \tensor ((1.101001000011101110)_{2}\times2^{4})
\cr
\vphantalt 3
&
1\tensor ((1.022110110201121122)_{3}\times3^{2}) + \z_3 \tensor ((2.102110200120211001)_{3}\times3^{0})
\cr
\vphantalt 5
&
1\tensor ((2.424401443414003124)_{5}\times5^{0}) + \z_3 \tensor ((4.440233034214103134)_{5}\times5^{0})
\cr
\vphantalt 7
&
1\tensor ((5.421026022055251515)_{7}\times7^{0}) + \z_3 \tensor ((4.510213122044500212)_{7}\times7^{0})
\cr
\vphantalt 11
&
1\tensor ((A.187435006878255608)_{11}\times11^{0}) + \z_3 \tensor ((6.47020389393AA07A54)_{11}\times11^{0})
\cr
\hline
\hline\multicolumn{2}{l}{\vphant n=5\hfill\ratE[5]^{-1}=\text{ see Example~\ref{sexticcomputation}}}\cr
\hline
\firstintableheadlineE[5]
\hline
\vphantalt 2
&
1\tensor ((1.001111100010111111)_{2}\times2^{7}) + \z_3 \tensor ((1.100110001101001011)_{2}\times2^{8})
\cr
\vphantalt 3
&
1\tensor ((1.100222220001012001)_{3}\times3^{15}) + \z_3 \tensor ((1.201100102012122102)_{3}\times3^{15})
\cr
\vphantalt 5
&
1\tensor ((1.302021433120040233)_{5}\times5^{8}) + \z_3 \tensor ((1.414010322010320220)_{5}\times5^{8})
\cr
\vphantalt 7
&
1\tensor ((6.54022035336014)_{7}\times7^{10}) + \z_3 \tensor ((1.31434160250212)_{7}\times7^{10})
\cr
\vphantalt 11
&
1\tensor ((9.5A062420171)_{11}\times11^{10}) + \z_3 \tensor ((A.873665790A8)_{11}\times11^{10})
\cr
\hline
\secondintableheadlineE[5]
\hline
\vphantalt 2
&
1\tensor ((1.101111101101011100)_{2}\times2^{2}) + \z_3 \tensor ((1.001110000011000011)_{2}\times2^{4})
\cr
\vphantalt 3
&
1\tensor ((2.001012111110022002)_{3}\times3^{1}) + \z_3 \tensor ((2.011221210002201001)_{3}\times3^{0})
\cr
\vphantalt 5
&
1\tensor ((3.241120123211444141)_{5}\times5^{1}) + \z_3 \tensor ((2.103430110032401010)_{5}\times5^{0})
\cr
\vphantalt 7
&
1\tensor ((5.454116060262152231)_{7}\times7^{0}) + \z_3 \tensor ((5.160242506230053303)_{7}\times7^{0})
\cr
\vphantalt 11
&
1\tensor ((9.9950088486A6479846)_{11}\times11^{0}) + \z_3 \tensor ((A.178492061888054088)_{11}\times11^{0})
\cr
\hline
\endarray
\endtable

\begin{example}\label{S3computation}
Let $ k $ be a totally real $ S_3 $-extension of $ \Q $.  Use notation as in Case~1 of Example~\ref{S3example}, and
let $ \pi $ be the primitive idempotent corresponding to $ \pi_2 $ of that example.  Since $ E=\Q $ we shall write
$ M_\pi $ instead of $ M_\pi^E $ in all notation referring to Conjecture~\ref{motiveconjecture}.

Since $ \pi \Kn k = \pi_2 \Kn c $ by Remark~\ref{galoisdescent}
we computed the $K$-theory of $ c $ using the (numerical) methods of Section~\ref{practical} and tried to recognize the number
\begin{equation}\label{enM}
e(n,M_\pi)= \L n {\chi_\pi} {\Q} \cdot D(M_\pi)^{1/2,\infty} / R_{n,\infty}(M_\pi) 
\end{equation}
for $ n = 3 $ and 5 as an element of $ \Q^* $ by employing the same methods as described in Section~\ref{practical}
for number fields, where $ \L n {\chi_\pi} {\Q} $ was computed as $ \z_c(n)/\z_\Q(n) $ (see Example~\ref{S3L}) using
Pari-GP~\cite{PARI2}.
We succeeded in all cases, and then verified numerically for $ p=2 $, 3, 5, 7 and 11 if
$ \Lp p n {\chi_\pi\omega_p^{1-n}} {\Q} \ne 0 $ and 
\begin{equation}\label{quotientexp}
e(n,M_\pi) \cdot \Eul_p(n, \chi_\pi, \Q) \cdot \dfrac{R_{n,p}(M_\pi)}{D(M_\pi)^{1/2,p}} 
\cdot \Lp p n {\chi_\pi\omega_p^{1-n}} {\Q} ^{-1}
\end{equation}
in $ \Q_p^* $ was equal to~1.

We put our results for the four totally real $ S_3 $-extensions $ k $ of $ \Q $ with smallest discriminants
in Tables~1 through~4.
Note that by Lemma~\ref{Qpquot} $ R_{n,p}(M_\pi) / D(M_\pi)^{1/2,p} $ is in $ \Q_p $
and the same holds for $  \Lp p n {\chi_\pi\omega_p^{1-n}} {\Q} $ by Lemma~\ref{pLvalueslemma},
so that (approximations of) those numbers can be easily represented in our tables.
We have denoted an element $ (a_0+a_1p+a_2p^2+\cdots) \times p^s $  in $ \Q_p^* $ 
with all $ a_j $ in $ \{0,\dots,p-1\} $ and $ a_0\ne0 $ by $ (a_0.a_1a_2\cdots)_p \times p^s $, writing
$ A $ to represent 10 when $ p = 11 $.

For the \padic regulators we computed each $ \Lmod n (z) $ up to $ O(p^{L(p)}) $ with $ p^{L(p)} $
approximately equal to $ 10^{30} $.
The relevant part of $ \preg[n] $ has been given in the tables whenever it fitted.
For the values of the \padic \Lfunction we can only prove that the
relative error is at most $  1+O(p^{M(p)}) $, 
where $ M(2) = 72$, $ M(3) = 47 $,  $ M(5) = 32 $, $ M(7) = 26 $ and $ M(11) = 22$.
This proves that $ \Lp p n {\chi_\pi\omega_p^{1-n}} {\Q} $ does not vanish;
together with the verification that~\eqref{enM} equals 1 numerically this proves part~(4)
of Conjecture~\ref{motiveconjecture}.
The value that we found for~\eqref{quotientexp} was $ 1+O(p^{N(p)}) $ with $ N(2) \ge 82 $, 
$ N(3) \ge 52 $, $ N(5) \ge 36  $, $ N(7) \ge 27 $ and $ N(11) \ge 24 $, giving
numerical evidence for parts~(2) and~(3) of the conjecture, but also 
suggesting that the relative precision of the value of the \padic \Lfunction
is slightly higher than we can prove, justifying the higher precision given in the tables.

Similar calculations were carried out for the next four such $ S_3 $ extensions
and the same primes, with very similar results.

\end{example}

The other examples proceed mostly along the same lines as Example~\ref{S3computation}.
The values of $ \L n {\chi_\pi} {\Q} $ in Examples~\ref{D8realcomputation} and~~\ref{D8CMcomputation} below
were computed using the formula at the end of Example~\ref{D8L}, but in Example~\ref{sexticcomputation}
the values of $ \L n \chi_\pi^E {\Q} $, with  $ E = \Q(\z_3) $ for $ \z_3 $ a primitive cube root of unity,
were computed using the algorithms described in \cite{Dok1} and its associated program.
In Examples~\ref{D8realcomputation} and \ref{D8CMcomputation} we determined
if $ e(n,\chi_\pi) $ (as in~\eqref{enM}) was in $ \Q^* $ using the methods of Section~\ref{practical}, and
for Example~\ref{sexticcomputation} we wrote $ e(n,M_\pi^E)^{-1} = 1 \tensor a + \z_3 \tensor b $ in
$ E \tensor_\Q \R $ and then recognized $ a $ and $ b $ as integers divided by products of some relatively small primes.
In each case we then verified for $ p=2,\dots,11 $ if~\eqref{quotientexp} 
(or its equivalent) was equal to 1 in $ \Qp $ or, for Example~\ref{sexticcomputation},
in~$ \Qp(\z_3)\tensor_\Q\Q = \Qp(\z_3) $.  The precision for $ \Lmod n (z) $ and the \padic \Lfunctions
in Examples~\ref{D8realcomputation} and~\ref{D8CMcomputation} was as in Example~\ref{S3computation}.
In Example~\ref{sexticcomputation} all $ \Lmod n (z) $ were computed up to $ O(p^{L(p)}) $ with $ p^{L(p)} $
approximately $ 10^{16} $ since in this case the check if the equivalent of ~\eqref{quotientexp} is~1
is done in $ \Q(\z_3) \tensor_\Q \Qp $ so two coefficients are checked;
the values of the \padic \Lfunction in this case, also in $ \Q(\z_3)\tensor_\Q \Qp $, were determined
up to multiplication by $ 1\tensor(1+O(p^{M(p)})) $ with $ M(p) $ as in Example~\ref{S3computation}.

\begin{example}\label{D8realcomputation}
Let $ k $ be a totally real $ D_8 $-extension of $ \Q $ as in Case~1 of Example~\ref{D8example}, 
and use notation as in that example.  We again write $ M_\pi $ instead of $ M_\pi^E $ since $ E=\Q $ and consider
odd $ n \ge 2 $.
Because $ \pi \Kn k = \pi_2 \Kn k_1 $ we computed $ \Kn k_1 $ numerically using the methods of Section~\ref{ksection}.
Our results with $ n=3 $ and 5 for the first four such extensions when ordered according to the discriminant
are in Tables~5 through~8.
The value that we found for~\eqref{quotientexp} in this case was $ 1+O(p^{N(p)}) $ where
$ N(2) \ge 76 $, $ N(3) \ge 52 $, $ N(5) \ge 34 $, $ N(7) \ge 29 $ and $ N(11) \ge 22 $.
Such calculations were also carried out for the next four such extensions, with similar results.
\end{example}

\begin{example}\label{D8CMcomputation}
Let $ k $ be a CM Galois extension of $ \Q $ with $ \Gal(k/\Q) \iso D_8 $ as in Case~2 of Example~\ref{D8example}, 
and use notation as in that example.  We again write $ M_\pi $ instead of $ M_\pi^E $ since $ E=\Q $, and consider
even $  n \ge 2 $.
Again $ \pi \Kn k = \pi_2 \Kn k_1 $, so we computed $ \Kn k_1 $ numerically using the methods of Section~\ref{ksection}.
We put our results at $ n=2 $ and 4 for the first four such extensions when ordered by the value of the
(positive) discriminant in Tables~9 through~12.
The value that we found for~\eqref{quotientexp} in this case was $ 1+O(p^{N(p)}) $ where
$ N(2) \ge 79 $, $ N(3) \ge 53 $, $ N(5) \ge 35 $, $ N(7) \ge 27 $ and $ N(11) \ge 25 $.
Such calculations were also carried out for the next four such extensions, with similar results.
\end{example}

\begin{example}\label{sexticcomputation}
Let $ k $ be a totally real $ S_3 \times \Z/3\Z $ extension of $ \Q $, as in Case~1 of Example~\ref{sexticexample}.
We use notation as in that example, taking $ E=\Q(\z_3) $ for a primitive cubic root of unity $ \z_3 $.
A different identification of $ \Gal(k/\Q) $ with $ S_3 \times \Z/3\Z $
might interchange $ \pi_3 $ and $ \pi_4 $ as well as
corresponding primitive idempotents, but, as we noted before, the validity of the conjecture for either is equivalent.

Since $ \pi \KnE k = \pi_3 \KnE s_1 $ we computed $ \Kn s_1 $ as before.
Also, in order to speed up the calculations of the regulator
we note that the action of $ (e,1) $ and $ (\s^2,0) $ on $ s_1 = k^{\langle (\s,1) \rangle } $
is the same, so that $ \pi = \pi_3 \pi_{H_1} $ and
\begin{equation*}
(1-\pi_1) ((e,0) + \z_3 (\s,0) + \z_3^2 (\s^2,0))/3 
=
((e,0) + \z_3 (\s,0) + \z_3^2 (\s^2,0))/3 
\end{equation*}
give the same result when applied to $ \KnE s_1 $.
We put our results for the first four $ s_1 $, when ordered according to their discriminant, for $ n=3 $ and 5 in Tables~13 through~16.
In Table~16 the number $ \ratE[5]^{-1} $ is given by
\[
-2^5\cdot3^{10}\cdot5^{-2}\cdot13^{4}\cdot(7+3\z_3)\cdot(68-43\z_3)\cdot(1547202-1603487\z_3)\cdot(a+b\z_3)
\]
with
\begin{alignat*}{1}
a & = 8482255892311139091186217686741233 
\\
\intertext{and}
b & = 3527440500058817421018094757947617
\,,
\end{alignat*}
presumably because we experimentally found a subgroup of $ \Kn s_1 $ that is not close to actually being
$ K_{2n-1}(s_1) $.
In all cases $ \Lp p n \chi_\pi\tensor\omega_p^{1-n} {\Q} $ is a unit in $ E\tensor_\Q\Qp $, and
when we write the equivalent of~\eqref{quotientexp} in the form $ 1 \tensor \a  + \z_3 \tensor \b $ then
we find that $ \a-1 $ and $ \b $ are $ O(p^{N(p)}) $ with
$ N(2) \ge 50 $, $ N(3) \ge 24 $, $ N(5) \ge 20 $, $ N(7) \ge 14 $ and $ N(11) \ge 12 $, with in this case
the precision bounded by that of the~$ R_{n,p}(M_\pi^E) $.
Since $ e(n,M_\pi^E) $, $ D(M_\pi^E)^{1/2,p} $ and $ \Eul_p(n, \chi_\pi\tensor\id,\Q) $ are units in
$ E\tensor_\Q\Qpbar $ this also shows that $ R_{n,p}(M_\pi^E) $ is a unit in $ E\tensor_\Q F $, as conjectured
in Conjecture~\ref{motiveconjecture}(4).

Such calculations were also carried out for the next four such $ s_1 $ with similar results.
\end{example}

\begin{remark}
By Proposition~\ref{abelprop} and Remark~\ref{consequencesremark}(1),
the conjecture for $ \pi $ in Example~\ref{S3computation} is equivalent to Conjecture~\ref{pborel}
for the field $ k_1 $ in the example.  The same statement holds for $ \pi $ and $ k_1 $ in
Example~\ref{D8realcomputation} below.
\end{remark}

\begin{remark}\label{constantremark}
In the tables we have also included the value of the  constant
$ C(\chi_\pi\tensor\id) = C(\chi_\pi^\vee\tensor\id) $
in the functional
equation~\eqref{newfuneq} of $ \L s {\chi_\pi\tensor\id} {\Q} $ to make possible a comparison between its
prime factors and those of $ e(n,M_\pi) $.
Namely, for the $ D(M_\pi^E)^{1/2,\infty} $ used in the calculations,
the $ e_n $ in~\eqref{edef} (which in our case
are independent of $ n $) are as follows.

\begin{center}
$$
\begin{array}{cc@{\qquad\quad}cc@{\qquad\quad}cc@{\qquad\quad}cc}
\hline
\vphant\hbox{Table} & e_n & \hbox{Table} & e_n &\hbox{Table} & e_n & \hbox{Table} & e_n \cr
\hline
\vphant
1 & 3 & 5 & 4 & 9 & 4 & 13 & 1\tensor 9 +\z_3\tensor 9
\cr
\vphant
2 & 3 & 6 & 4 & 10 & 4 & 14 & 1\tensor 9 +\z_3\tensor 9
\cr
\vphant
3 & 3 & 7 & 4 & 11 & \hbox to 20.8 pt{$-4$\hfill} & 15 & 1\tensor 0 +\z_3\tensor 9
\cr
\vphant
4 & 3 & 8 & 4 & 12 & 4 & 16 & 1\tensor 0 +\z_3\tensor (-9)
\cr
\hline
\end{array}
$$
\end{center}

\bigskip\noindent
As can be seen from these values and Tables~1-16, in the number
\begin{equation*}
\delta_n e_n ((n-1)!/2)^m C(\chi_\pi^\vee\tensor\id)^{n-1} e(n,M_\pi^E) 
\end{equation*}
of~\eqref{newe} the larger prime factors occurring in $ e(n,M_\pi^E) $ are often
cancelled by those of $ C(\chi_\pi^\vee\tensor\id)^{n-1} $.
This suggests that the formulation of Beilinson's
conjecture at $ s=1-n $ normally involves simpler prime factors in $ E^* $ than the formulation at $ s=n $
in Conjecture~\ref{motiveconjecture}(1), with the remaining complicated factors in our examples
quite possibly due to an awkward choice of basis for $ \KMn M_\pi^E $.
\end{remark}


\begin{thebibliography}{10}

\bibitem{Bar78}
D.~Barsky.
\newblock Fonctions z\^eta $p$-adiques d'une classe de rayon des corps de
  nombres totalement r\'eels.
\newblock In {\em Groupe d'Etude d'Analyse Ultram\'etrique (5e ann\'ee:
  1977/78)}, pages Exp. No. 16, 23. Secr\'etariat Math., Paris, 1978.

\bibitem{bei85}
A.~A. Beilinson.
\newblock Higher regulators and values of ${L}$-functions.
\newblock {\em J. Sov. Math.}, 30:2036--2070, 1985.

\bibitem{Bes98a}
A.~Besser.
\newblock Syntomic regulators and $p$-adic integration {I}: rigid syntomic
  regulators.
\newblock {\em Israel Journal of Math.}, 120:291--334, 2000.

\bibitem{Bes98b}
A.~Besser.
\newblock Syntomic regulators and $p$-adic integration {II}: {$K_2$} of curves.
\newblock {\em Israel Journal of Math.}, 120:335--360, 2000.

\bibitem{BdJ03}
A.~Besser and R.~de~Jeu.
\newblock The syntomic regulator for the $ {K} $-theory of fields.
\newblock {\em Annales Scientifiques de l'\'Ecole Normale Sup\'erieure},
  36(6):867--924, 2003.

\bibitem{BdJ06}
A.~Besser and R.~de~Jeu.
\newblock $ {L}i^{(p)} $-service? {A}n algorithm for computing $ p $-adic
  polylogarithms.
\newblock To appear in {\it Mathematics of Computation}.

\bibitem{Bes-Den99}
A.~Besser and C.~Deninger.
\newblock $p$-adic {M}ahler measures.
\newblock {\em J. Reine Angew. Math.}, 517:19--50, 1999.

\bibitem{Beu06}
F.~Beukers.
\newblock Irrationality of some $p$-adic ${L}$-values, 2006.
\newblock Preprint available from http://arxiv.org/abs/math.NT/0603277.

\bibitem{bl00}
S.~Bloch.
\newblock Higher regulators, algebraic ${K}$-theory, and zeta functions of
  elliptic curves.
\newblock Manuscript from 1978 (``Irvine notes"). Published as volume 11 of CRM
  Monographs Series by American Mathematical Society.

\bibitem{Borel74}
A.~Borel.
\newblock Stable real cohomology of arithmetic groups.
\newblock {\em Ann. {S}ci. {E}{N}{S}}, 4:235--272, 1974.

\bibitem{Borel77}
A.~Borel.
\newblock Cohomologie de {${\rm SL}\sb{n}$} et valeurs de fonctions z\^eta aux
  points entiers.
\newblock {\em Ann. Scuola Norm. Sup. Pisa Cl. Sci. (4)}, 4(4):613--636, 1977.
\newblock Errata at vol. 7, p. 373 (1980).

\bibitem{Bur02}
J.~I. Burgos~Gil.
\newblock {\em The regulators of {B}eilinson and {B}orel}, volume~15 of {\em
  CRM Monograph Series}.
\newblock Amer. Math. Soc., Providence, RI, 2002.

\bibitem{Cal05}
F.~Calegari.
\newblock Irrationality of certain {$p$}-adic periods for small {$p$}.
\newblock {\em Int. Math. Res. Not.}, (20):1235--1249, 2005.

\bibitem{Cas-Nog78}
P.~Cassou-Nogu{\`e}s.
\newblock {$p$}-adic {$L$}-functions for totally real number field.
\newblock In {\em Proceedings of the Conference on $p$-adic Analysis (Nijmegen,
  1978)}, volume 7806 of {\em Report}, pages 24--37. Katholieke Univ. Nijmegen,
  1978.

\bibitem{Cas-Nog79}
P.~Cassou-Nogu{\`e}s.
\newblock Valeurs aux entiers n\'egatifs des fonctions z\^eta et fonctions
  z\^eta {$p$}-adiques.
\newblock {\em Invent. Math.}, 51(1):29--59, 1979.

\bibitem{Col82}
R.~Coleman.
\newblock Dilogarithms, regulators, and $p$-adic {$L$}-functions.
\newblock {\em Invent. Math.}, 69:171--208, 1982.

\bibitem{Col-Gro89}
R.~Coleman and B.~Gross.
\newblock $p$-adic heights on curves.
\newblock In {\em Algebraic number theory}, pages 73--81. Academic Press,
  Boston, MA, 1989.

\bibitem{Colm88}
P.~Colmez.
\newblock R\'esidu en {$s=1$} des fonctions z\^eta {$p$}-adiques.
\newblock {\em Invent. Math.}, 91(2):371--389, 1988.

\bibitem{dJ95}
R.~de~Jeu.
\newblock Zagier's conjecture and wedge complexes in algebraic ${K}$-theory.
\newblock {\em Compositio Mathematica}, 96:197--247, 1995.

\bibitem{Del77}
P.~Deligne.
\newblock Valeurs de fonctions {$L$} et p\'eriodes d'int\'egrales.
\newblock In {\em Automorphic forms, representations and $L$-functions (Proc.
  Sympos. Pure Math., Oregon State Univ., Corvallis, Ore., 1977), Part 2},
  Proc. Sympos. Pure Math., XXXIII, pages 313--346. Amer. Math. Soc.,
  Providence, R.I., 1979.
\newblock With an appendix by N. Koblitz and A. Ogus.

\bibitem{del:imd}
P.~Deligne.
\newblock Interpr{\'e}tation motivique de la conjecture de {Z}agier reliant
  polylogarithms et r{\'e}gulateurs.
\newblock Preprint, 1991.

\bibitem{Del-Rib}
P.~Deligne and K.~Ribet.
\newblock Values of abelian {$L$}-functions at negative integers over totally
  real fields.
\newblock {\em Invent. Math.}, 59(3):227--286, 1980.

\bibitem{DenSch91}
C.~Deninger and A.~J. Scholl.
\newblock The {B}e\u\i linson conjectures.
\newblock In {\em $L$-functions and arithmetic (Durham, 1989)}, volume 153 of
  {\em London Math. Soc. Lecture Note Ser.}, pages 173--209. Cambridge Univ.
  Press, Cambridge, 1991.

\bibitem{Dok1}
T.~Dokchitser.
\newblock Computing special values of motivic ${L}$-functions.
\newblock {\em Experiment. Math.}, 13(2):137--149, 2004.

\bibitem{gonXpam}
A.~B. Goncharov.
\newblock Polylogarithms and motivic {G}alois groups.
\newblock In {\em Motives (Seattle, WA, 1991)}, volume~55 of {\em Proc. Sympos.
  Pure Math.}, pages 43--96. Amer. Math. Soc., Providence, RI, 1994.

\bibitem{gon:goc}
A.~B. Goncharov.
\newblock Geometry of configurations, polylogarithms, and motivic cohomology.
\newblock {\em Adv. Math.}, 114(2):197--318, 1995.

\bibitem{G3}
A.~B. Goncharov.
\newblock Deninger's conjecture of {$L$}-functions of elliptic curves at
  {$s=3$}.
\newblock {\em J. Math. Sci.}, 81(3):2631--2656, 1996.
\newblock Algebraic geometry, 4.

\bibitem{greenberg83}
R.~Greenberg.
\newblock On {$p$}-adic {A}rtin {$L$}-functions.
\newblock {\em Nagoya Math. J.}, 89:77--87, 1983.

\bibitem{Gro90}
M.~Gros.
\newblock R\'egulateurs syntomiques et valeurs de fonctions ${L}$ $p$-adiques.
  {I}.
\newblock {\em Invent. Math.}, 99(2):293--320, 1990.
\newblock With an appendix by Masato Kurihara.

\bibitem{Gro94}
M.~Gros.
\newblock R\'egulateurs syntomiques et valeurs de fonctions ${L}\;p$-adiques.
  {I}{I}.
\newblock {\em Invent. Math.}, 115(1):61--79, 1994.

\bibitem{Jan88b}
U.~Jannsen.
\newblock Deligne homology, {H}odge-{${\mathcal D}$}-conjecture, and motives.
\newblock In {\em Be\u\i linson's conjectures on special values of
  $L$-functions}, volume~4 of {\em Perspect. Math.}, pages 305--372. Academic
  Press, Boston, MA, 1988.

\bibitem{Kat76}
N.~Katz.
\newblock {$p$}-adic interpolation of real analytic {E}isenstein series.
\newblock {\em Ann. of Math. (2)}, 104(3):459--571, 1976.

\bibitem{Kat78}
N.~Katz.
\newblock {$p$}-adic {$L$}-functions for {CM} fields.
\newblock {\em Invent. Math.}, 49(3):199--297, 1978.

\bibitem{Kat94}
N.~Katz.
\newblock Review of {$l$}-adic cohomology.
\newblock In {\em Motives (Seattle, WA, 1991)}, volume~55 of {\em Proc. Sympos.
  Pure Math.}, pages 21--30. Amer. Math. Soc., Providence, RI, 1994.

\bibitem{Kub-Leo}
T.~Kubota and H.-W. Leopoldt.
\newblock Eine {$p$}-adische {T}heorie der {Z}etawerte. {I}. {E}inf\"uhrung der
  {$p$}-adischen {D}irichletschen {$L$}-{F}unktionen.
\newblock {\em J. Reine Angew. Math.}, 214/215:328--339, 1964.

\bibitem{Lan70}
S.~Lang.
\newblock {\em Algebraic number theory}.
\newblock Addison-Wesley Publishing Co., Inc., Reading, Mass.-London-Don Mills,
  Ont., 1970.

\bibitem{MTT86}
B.~Mazur, J.~Tate, and J.~Teitelbaum.
\newblock On {$p$}-adic analogues of the conjectures of {B}irch and
{S}winnerton-{D}yer.
\newblock {\em Invent. Math.}, 84(1):1--48, 1986.

\bibitem{MW84}
B.~Mazur and A.~Wiles.
\newblock Class fields of abelian extensions of {${\bf Q}$}.
\newblock {\em Invent. Math.}, 76(2):179--330, 1984.
		
\bibitem{Neu88}
J.~Neukirch.
\newblock The {B}e\u\i linson conjecture for algebraic number fields.
\newblock In {\em Be\u\i linson's conjectures on special values of
  $L$-functions}, volume~4 of {\em Perspect. Math.}, pages 193--247. Academic
  Press, Boston, MA, 1988.

\bibitem{Neu99}
J.~Neukirch.
\newblock {\em Algebraic number theory}, volume 322 of {\em Grundlehren der
  Mathematischen Wissenschaften [Fundamental Principles of Mathematical
  Sciences]}.
\newblock Springer-Verlag, Berlin, 1999.
\newblock Translated from the 1992 German original and with a note by Norbert
  Schappacher, With a foreword by G. Harder.

\bibitem{Niz97}
W.~Nizio\l.
\newblock On the image of $p$-adic regulators.
\newblock {\em Invent. Math.}, 127:375--400, 1997.

\bibitem{Peri96}
B.~Perrin-Riou.
\newblock Fonctions ${L}$ $p$-adiques des repr\'esentations $p$-adiques.
\newblock {\em Ast\'erisque}, (229):198pp, 1995.

\bibitem{qui73b}
D.~Quillen.
\newblock Finite generation of the groups ${K}\sb{i}$ of rings of algebraic
  integers.
\newblock In {\em Algebraic ${K}$-theory {1}}, volume 341 of {\em Lecture
  {N}otes in {M}athematics}, pages 179--198. Springer {V}erlag, Berlin, 1973.

\bibitem{Rob06}
X.-F. Roblot.
\newblock Computing values of {$p$}-adic {$L$}-functions of totally real number
  fields.
\newblock In preparation.

\bibitem{Sol-Rob}
X.-F. Roblot and D.~Solomon.
\newblock Verifying a {$p$}-adic abelian {S}tark conjecture at {$s=1$}.
\newblock {\em J. Number Theory}, 107(1):168--206, 2004.

\bibitem{schn88}
P.~Schneider.
\newblock Introduction to the {B}eilinson {C}onjectures.
\newblock In {\em Beilinson's {C}onjectures on {S}pecial {V}alues of
  $L$-Functions}, pages 1--35. Academic Press, Boston, MA, 1988.

\bibitem{Scho98}
A.~J. Scholl.
\newblock An introduction to {K}ato's {E}uler systems.
\newblock In {\em Galois representations in arithmetic algebraic geometry
  (Durham, 1996)}, pages 379--460. Cambridge Univ. Press, Cambridge, 1998.

\bibitem{Ser73}
J.-P. Serre.
\newblock Formes modulaires et fonctions z\^eta {$p$}-adiques.
\newblock In {\em Modular functions of one variable, III (Proc. Internat.
  Summer School, Univ. Antwerp, 1972)}, pages 191--268. Lecture Notes in Math.,
  Vol. 350. Springer, Berlin, 1973.

\bibitem{Shin76}
T.~Shintani.
\newblock On evaluation of zeta functions of totally real algebraic number
  fields at non-positive integers.
\newblock {\em J. Fac. Sci. Univ. Tokyo Sect. IA Math.}, 23(2):393--417, 1976.

\bibitem{sou81}
C.~Soul{\'e}.
\newblock On higher {$p$}-adic regulators.
\newblock In {\em Algebraic $K$-theory, Evanston 1980 (Proc. Conf.,
  Northwestern Univ., Evanston, Ill., 1980)}, volume 854 of {\em Lecture Notes
  in Math.}, pages 372--401. Springer Verlag, Berlin, 1981.

\bibitem{susXkof}
A.~A. Suslin.
\newblock {$K\sb 3$} of a field, and the {B}loch group.
\newblock {\em Trudy Mat. Inst. Steklov.}, 183:180--199, 229, 1990.
\newblock Translated in Proc.\ Steklov Inst.\ Math.\ {\bf 1991}, no.\ 4,
  217--239, Galois theory, rings, algebraic groups and their applications
  (Russian).

\bibitem{magma}
{The Magma~group}, Sydney.
\newblock {\em {Magma}}.
\newblock Available from \url{http://magma.maths.usyd.edu.au/}.

\bibitem{PARI2}
{The PARI~Group}, Bordeaux.
\newblock {\em {PARI/GP, version {\tt 2.2.11}}}, 2005.
\newblock Available from \url{http://pari.math.u-bordeaux.fr/}.

\bibitem{Zag91}
D.~Zagier.
\newblock Polylogarithms, {D}edekind {Z}eta {F}unctions and the {A}lgebraic
  ${K}$-theory of {F}ields.
\newblock In {\em Arithmetic algebraic geometry (Texel, 1989)}, pages 391--430.
  Birkh\"auser Boston, Boston, MA, 1991.

\end{thebibliography}
\end{document}